\documentclass[3p,times,preprint]{elsarticle}

\usepackage{capt-of}
\usepackage{nicefrac}

\usepackage{amssymb}

\usepackage{dsfont}
\usepackage{amsmath,amssymb,amsfonts, comment,tikz}
\usepackage{inputenc}

\usepackage{multirow}
\usepackage{amsthm}
\usepackage{epsfig} 
\usepackage{bm}
\usepackage{graphicx}
\usepackage{textcomp}
\usepackage{xcolor}
\usepackage{hyperref}
\usepackage{cleveref}

\usepackage{enumerate}
\usepackage{ifthen}
\usepackage{fontsize}

\usepackage[ruled,vlined]{algorithm2e}

\usepackage{fancyhdr}
\usepackage{lastpage}
\def\BibTeX{{\rm B\kern-.05em{\sc i\kern-.025em b}\kern-.08em T\kern-.1667em\lower.7ex\hbox{E}\kern-.125emX}}

\newcommand{\revg}[1]{{\color{black}#1}}
\newcommand{\reva}[1]{{\color{black}#1}}
\newcommand{\revr}[1]{{\color{black}#1}}

\usepackage[figuresright]{rotating}
\usepackage{subfigure}
\usepackage{enumerate}
\usepackage{ulem}

\usepackage{graphicx}
\usepackage{enumitem}
\usepackage{xcolor}
\usepackage{nicefrac}

\usepackage{bm}
\usepackage{wrapfig}
\usepackage{amsmath}
\usepackage{amssymb}

\usepackage{hyperref}

\newtheorem{theorem}{Theorem}
\newtheorem{lem}{Lemma}
\newdefinition{rmk}{Remark}
\newtheorem{corollary}{Corollary}
\newproof{pf}{Proof}
\newproof{pfthrmone}{\textbf{Proof of Theorem \ref{thm1}}}
\newproof{pfthrmtwo}{\textbf{Proof of Theorem \ref{thrm_beta_ODE_prop}}}

\newproof{pfthrmthree}{\textbf{Proof of Theorem \ref{thrm_BP_to_fake}}}
\newproof{pfthrmfour}{\textbf{Proof of Theorem \ref{thrm_unique_att}}}
\newproof{pfthrmfive}{\textbf{Proof of Theorem \ref{thrm_opt}}}
\newproof{pflimitswarning}{\textbf{Proof of Corollary \ref{cor_limits_warning}}}
\newproof{pfthrmsix}{\textbf{Proof of Theorem \ref{thrm_Q_old}}}
\newproof{pfcorollaryQ}{\textbf{Proof of Corollary \ref{cor_Q_old}}}
\newproof{pfthrmseven}{\textbf{Proof of Theorem \ref{thrm_BP_to_fake_ea}}}
\newproof{pfcorollaryeaWM}{\textbf{Proof of Theorem \ref{corollary_ea_wm}}}
\newproof{pfcorollaryexWM}{\textbf{Proof of Corollary \ref{corollary_ex_wm}}}
\newproof{pfcorollaryehWM}{\textbf{Proof of Theorem \ref{corollary_eh_WM}}}
\newproof{pfcorbetaona}{\textbf{Proof of Corollary \ref{cor_beta_o_na}}}

\begin{document}

\begin{frontmatter}



\title{Robust fake-post detection against real-coloring adversaries}


\author{Khushboo Agarwal\corref{cor1}\fnref{label2}}
\ead{agarwal.khushboo@iitb.ac.in}

\author{Veeraruna Kavitha}
\ead{vkavitha@iitb.ac.in}
\fntext[label2]{The work of first author is partially supported by Prime Minister’s Research Fellowship (PMRF), India.}
\cortext[cor1]{Corresponding author}

\address{IEOR, IIT Bombay, Powai, Mumbai, 400076, Maharashtra, India}

\begin{abstract}
The viral propagation of fake posts on online social networks (OSNs) has become an alarming concern. The paper aims to design control mechanisms for fake post detection while negligibly affecting the propagation of real posts. Towards this, a warning mechanism based on crowd-signals was recently proposed, where all users actively declare the post as real or fake. In this paper, we consider a more realistic framework where users exhibit different adversarial or non-cooperative behaviour: (i) they can independently decide whether to provide their response, (ii) they can choose not to consider the warning signal while providing the response, and (iii) they can be real-coloring adversaries who deliberately declare any post as real. To analyze the post-propagation process in this complex system, we propose and study a new branching process, namely total-current population-dependent branching process with multiple death types. At first, we compare and show that the existing warning mechanism significantly under-performs in the presence of adversaries. Then, we design new mechanisms which remarkably perform better than the existing mechanism by cleverly eliminating the influence of the responses of the adversaries. \revr{Finally, we propose another enhanced mechanism which assumes minimal knowledge about the user-specific parameters.} The theoretical results are validated using Monte-Carlo simulations.
\end{abstract}




\end{frontmatter}

\newcommand{\tcprocess}{total-current population-dependent BP }
\newcommand{\tcprocessnospace}{total-current population-dependent BP}
\newcommand{\hide}[1]{}

\newcommand{\psiL}{\psi^c_{\mbox{\tiny o}}}
\newcommand{\psiaL}{\psi^a_{\mbox{\tiny o}}}

\newcommand{\bL}{\beta^c_{\mbox{\tiny o}}}
\newcommand{\baL}{\beta^a_{\mbox{\tiny o}}}

\newcommand{\tL}{\theta^c_{\mbox{\tiny o}}}
\newcommand{\taL}{\theta^a_{\mbox{\tiny o}}}
\newcommand{\sa}{z}
\newcommand{\Sa}{Z}
\newcommand{\N}{\mathcal{N}}
\newcommand{\I}{I_{\theta/\psi}}
\newcommand{\eop}{\hfill{$\square$}}
\newcommand{\nto}{\nrightarrow}
\newcommand{\Om}{\Phi}
\newcommand{\om}{\phi}
 \newcommand{\ups}{{\mbox{\small ${\Upsilon}$}}}
 \newcommand{\Ups}{{\mathbf \Upsilon}}
\newcommand{\Bin}{{\cal B}}
\newcommand{\tp }{\tau^+}
\newcommand{\tm}{\tau^-}
\newcommand{\up}{\uparrow}
\newcommand{\offs}{\Gamma}
\newcommand{\down}{\hspace{.15mm}\downarrow}
\newcommand{\beq}{\begin{eqnarray*}}
\newcommand{\eeq}{\end{eqnarray*}}
\newcommand{\Cx}{C_x}
\newcommand{\Cy}{C_y}
\newcommand{\cx}{c_x}
\newcommand{\cy}{c_y}
\newcommand{\cM}{{\cal E}}

\newcommand{\SampS}{{\tiny \mathbb{S}}}

\newcommand{\betax}{{\overline c^x}}
\newcommand{\betay}{{\overline c^y}}
\newcommand{\bax}{{\overline a^x}}
\newcommand{\Ax}{T_x}
\newcommand{\Ay}{T_y}
\newcommand{\ax}{t_x}
\newcommand{\ay}{t_y}
\newcommand{\wm}{\widetilde{e}}
\newcommand{\cA}{{\cal A}}
\newcommand{\cR}{{\cal S}}
\newcommand{\cRo}{{\cal R}_e}
\newcommand{\tc}{\theta^c}
\newcommand{\ta}{\theta^a}
\newcommand{\pc}{\psi^c}
\newcommand{\pa}{\psi^a}
\newcommand{\Tc}{\Theta^c}
\newcommand{\Ta}{\Theta^a}
\newcommand{\Pc}{\Psi^c}
\newcommand{\Pa}{\Psi^a}
\newcommand{\ba}{\beta^a}
\newcommand{\Bc}{\mathrm{B}^c}
\newcommand{\Ba}{\mathrm{B}^a}
\newcommand{\Beta}{\mathrm{B}}
\newcommand{\dist}{ d_{st}}
\newcommand{\st}{{\cal T}}
\newcommand{\bstar}{\beta^*}
\newcommand{\minf}{m^\infty}
\newcommand{\ueps}{\overline{\varepsilon} }
\newcommand{\leps}{\underline{\varepsilon}}
\newcommand{\polya}{\mbox{P\'{o}lya} }
\newcommand{\q}{\mathbf{q}}
\newcommand{\newbeta}{f_{\beta}^\infty(\beta)}
\newcommand{\overomega}{\overline{\omega}(\beta)}
\newcommand{\ga}{\mathbf{g}}
\newcommand{\gna}{\bm{\varrho}}

\newcommand{\cS}{\mathcal{D}_b}
\newcommand{\cD}{\mathcal{D}}
\newcommand{\cB}{\mathcal{B}}

\newcommand{\betana}{\beta^o_{\mbox{{\footnotesize na}}}}
\newcommand{\bpam}{e}
\newcommand{\bpaum}{\overline{e}}
\newcommand{\bpalm}{\underline{e}}
\newcommand{\propum}{\overline{m}}
\newcommand{\proplm}{\underline{m}}
\newcommand{\propwm}{m^\infty}
\newcommand{\propm}{m}

\newcommand{\G}{A_{\tiny{\mathbb{G}}}}

\newcommand\overlinebelow[1]{\stackunder[1.2pt]{$#1$}{\rule{1.2ex}{.075ex}}}
\def\lc{\left\lceil}   
\def\rc{\right\rceil}

\section{Introduction}
\label{sec_intro}
The prevalence of online social networks (OSNs), like Facebook or Twitter, is unprecedented today. A variety of content is available on the OSNs for users to consume, which can either be for education, entertainment, advertisement or awareness purposes, among many more. 
Users also read news on such platforms instead of using classical mediums like newspapers.

One of the reasons for such high usage of OSNs is the ease with which users can access or share information. Further, there is no instant check to ensure that the shared post is authentic. On one hand, this freedom allows users to express their views freely, but at the same time, it provides users with the flexibility to post fake content - the one that contains fabricated (mis)information that propagates through OSNs like authentic posts (see \cite{lazer2018science}). Once a post is shared on the OSN with an initial set of users, called seed users, the post can be further shared repeatedly by the recipients of the post to the extent of getting viral (the copies of the post grow significantly with time), or the post can get extinct in the initial phase (\cite{ranbir2019decomposable, agarwal2021co, agarwal2022saturated, van2010viral}). 

Now, there are several reasons for a fake-post to get viral. Authors in \cite{talwar2019people} theorize that users may share any information obtained from their reliable
source, or they can share any exciting post to seek their peers’
attention and have a sense of belonging. Also, users share posts that match their beliefs to continue using social media (due to its perceived usefulness). 
There have been many instances in the past where fake-posts have proven to be fatal, and the most controversial of all is the 2016 US Presidential elections (\cite{allcott2017social}). Thus, studies on the generation, propagation, detection, and control of fake posts are the need of the hour. In this paper, we focus on the detection aspect of fake-posts.

Machine learning or deep learning is one of the commonly used approaches for fake-post detection (see \cite{feng2022misreporting, sharma2019combating, ruchansky2017csi, ahmed2021detecting}). However, as argued in \cite{ahmed2021detecting}, such algorithms often face difficulty in obtaining training datasets in certain languages, and it gets difficult to determine the actuality using only the content (\cite{sharma2019combating}). 
Another approach used for fake-post identification is using crowd-signals. The basic idea is to allow users to declare any post as real or fake, and then leverage user responses to identify the actuality of the post. Such an approach is being used by Facebook\footnote{\url{https://www.facebook.com/help/1753719584844061}}, where any user can report any post on the OSN. They can also provide specific reasons for reporting the post. When a post is reported, it is reviewed by third-party fact-checking organizations and is removed if it is against their policies. 
However, until the post is reviewed, the users on the OSN can view it without any warning. 

In \cite{kapsikar2020controlling}, the authors design a warning-based mechanism to control fake-posts using crowd-signals. The idea is to leverage users' fake/real responses (tags) to the post and generate a warning signal for future recipients. Since the real-time warning signal/status of the post is continuously displayed to the users, this approach of using crowd-signals is different and should be more effective than that of Facebook. The objective is to ensure the maximal correct identification of the fake-post, while maintaining the proportion of fake-tags for the real-post within a given threshold. The paper assumes that each user participates in the tagging process. 

In this paper, we consider a more realistic framework. Firstly, we assume that not all users would be willing to tag. Secondly, if a user tags, it can consider the warning signal provided by the OSN; or it can tag without viewing the warning. And lastly, the users can be adversarial- these users always assign the real-tag to any post.

For such a system, we compare and show that the warning mechanism in \cite{kapsikar2020controlling} is insufficient. With just $1\%$ (with $2\%$) adversaries in the system, while everyone else tagging exactly as in \cite{kapsikar2020controlling}, we observed that the performance decreases approximately by $10\%$ (nearly $18.2\%$). This observation highlights the need for mechanisms which are robust against adversaries. We precisely achieve the same in this paper. 

The new warning mechanisms are designed by cleverly eliminating the effect of adversarial users. We derive a one-dimensional ordinary differential equation (ODE) that captures the performance of any such general warning mechanism, and utilizing that ODE, we design the new warning mechanisms as well as illustrate the improved performance guarantees theoretically. 

\revr{We have also presented Monte-Carlo simulation-based exhaustive numerical study to confirm our theoretical findings. The performance is expressed in two ways: (i) quality of service (QoS) which measures the proportion of fake-tags for the fake-post, and (ii) improved QoS (i-QoS) which represents the proportions only from non-adversarial users. The second metric i-QoS provides better interpretation for the performance of warning mechanisms, as actions of adversarial users can not be controlled. Note that, accordingly the threshold with respect to the real-post also changes, to consider the responses only from non-adversarial users. 

According to the parameters in \cite{kapsikar2020controlling}, the non-adversarial users are assumed to be smart (i.e., have high intrinsic ability to identify the actuality of the posts). Thus, no warning mechanism can accentuate their ability beyond a limit -- we observe minor improvements in QoS of $2.66\%$ and $5.34\%$ with $1\%$ and $2\%$ of adversary respectively; these numbers translate to $98.64\%$ and $98.63\%$ of i-QoS under new mechanisms as compared to $95.8\%$ and $92.53\%$ with the mechanism as in \cite{kapsikar2020controlling}.

In another instance, where users are less informed and more likely to wrongly recognize the posts (as is the case in reality), significant improvements are noticed even for a larger fraction of adversaries. Under newly proposed mechanism, the QoS is $52.89\%$ (i-QoS is $80.86\%$), which is only $45.31\%$ (i-QoS is only $45.31\%$) under old mechanism, when $32.5\%$ of adversarial users are involved. In fact, this performance is achieved with minimal knowledge about users sensitivity to the warning, and their behavioural type. }
 
The warning dynamics are modelled using a new variant of branching processes (BPs).
This paper also contributes towards total-current population-dependent two-type branching processes with population dependent death rates and also considers a variety of unnatural deaths. In particular, we derive all possible limits and limiting behaviours of the population sizes as time progresses.

\textbf{Related Literature for Branching processes with unnatural deaths:}  
The literature on BPs has previously investigated unnatural deaths in a restricted setting. The BP analyzed in \cite{BPwithinteraction} is population-independent, while the authors in \cite{etheridge2013conditioning} consider unnatural deaths due to competition, modelled using a quadratic function of population size. The BP with pairwise interaction in \cite{ojeda2020branching} models natural births and deaths, along with additional births and deaths occurring due to cooperation and competition. Further, the birth and death rates in \cite{ojeda2020branching} are proportional to current population sizes. Our work provides a much more generalised framework where the interactions are not limited to cooperation or competition. Further, the birth and death rate functions can additionally depend on the total and current population-sizes. 

\section{Problem description}\label{sec_prob_desc}
Consider an OSN with a large user base like Facebook or Twitter. Any post, $u$ on the OSN can be either fake ($u = F$) or real ($u = R$). The OSN aims to identify the actuality of the post. In \cite{kapsikar2020controlling}, the authors have proposed a warning mechanism where the recipients of the post themselves are guided in such a way that it leads to correct identification. We first study its robustness against adversarial users and then propose improved mechanisms.

We begin by describing the system and the warning mechanism of \cite{kapsikar2020controlling}. The posts are stored in a last-in-show-at-top structure named timeline for each user. The users are provided a warning for each post, and are asked to assign a tag (fake or real) to it. Whenever a user views the post on its timeline, it guesses the actuality of the post, assigns the tag as real or fake accordingly and then forwards the same to its friends. This results in more unread copies of the post tagged as fake or real. 
The process continues when another user with the post on its timeline visits the OSN. The warning mechanism relies on the tags provided by the users and is updated with each new tag. 

We will now introduce a few notations and then discuss the propagation and tagging dynamics of the post. Let the fake and real tagged copies of the $u$-post be denoted as $x$-type and $y$-type, respectively. Further, let $\Cx(t)$ and $\Cy(t)$ be the number of users who have received the post tagged as fake and real, respectively but have not yet read/shared it; thus, these are the number of unread copies of the post with fake or real tag. The total number of users who have received the post tagged as fake or real are represented by $\Ax(t)$ and $\Ay(t)$ respectively; these are read plus unread copies of the post. Let $\Om(t) := (\Cx(t), \Cy(t), \Ax(t), \Ay(t))$ be the tuple of number of copies at time $t$.

\begin{wrapfigure}{r}{6.2cm}
    \centering
    \includegraphics[trim = {1.2cm 0cm 7.2cm 2cm}, clip, scale = 0.35]{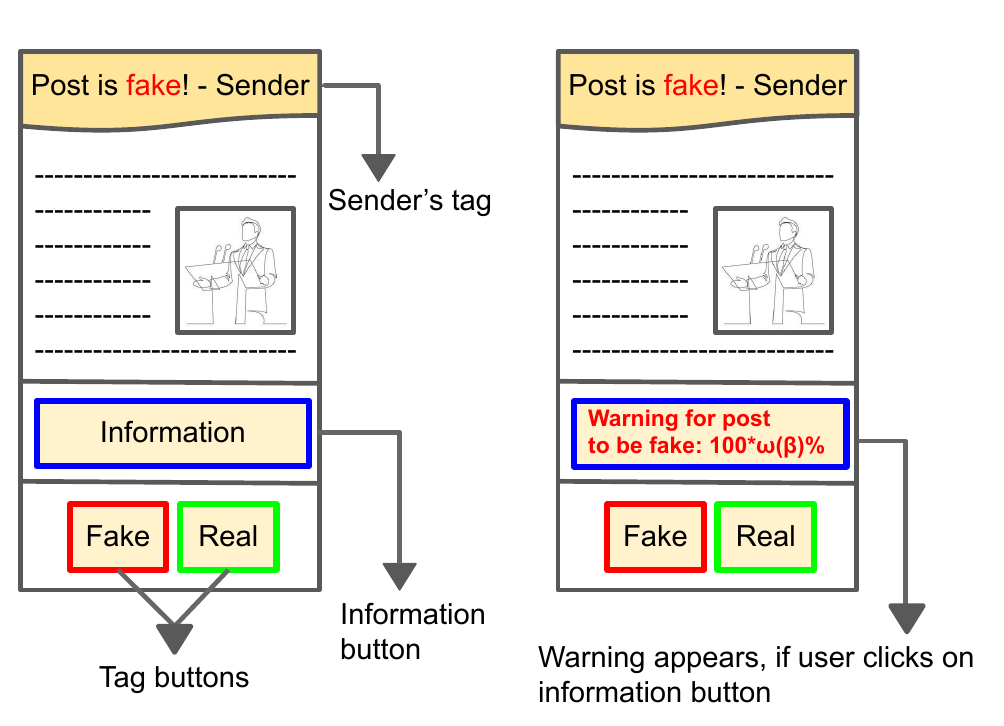}
    \caption{Design of the post}
    \label{fig_post_design}
\end{wrapfigure}

Each post contains two pieces of information: first, the sender's tag and second, the warning by the OSN, which is available at the click of a button (see Figure \ref{fig_post_design}). Users can exhibit different behaviours about utilising the provided information. For example, some users may prefer to read the warning before tagging, while others may not. Therefore, motivated by \cite{agarwal2023single}, we broadly divide user behaviour into four categories. 

\subsection{Warning-ignoring (wi) users} These are the users who tag the post only based on the sender's tag and their intrinsic ability to judge the post's actuality, not the warning. They prefer to invest less time in the system. Let $\tau$ be the time when a wi-user (with an unread copy of the post) reads it. At this time, the user will tag and then share the post with its friends. Let $I_{x, wi}(\Om(\tau^-))$ and $I_{y,wi}(\Om(\tau^-))$ be the indicator that the $wi$-user with fake or real-tagged copy of the post tags it as fake.

If the sender has tagged the post as fake, then the recipient tags the post as fake or real with probability (w.p.) $p_x^u \in (0,1)$ and $1-p_x^u$ respectively. Similarly, let $p_y^u \in (0, 1)$ be the probability of fake-tagging the post, received with real-tag. Therefore:
\begin{align}\label{eqn_tag_wi}
    P(I_{x,wi}(\Om(\tau^-)) = 1|{\cal G}_\tau) &= p_x^u \mbox{ and } P(I_{y, wi}(\Om(\tau^-)) = 1|{\cal G}_\tau) = p_y^u, 
\end{align}where $\mathcal{G}_t$ is the sigma-algebra generated by $\{\Om(t'); t'\le t\}$. Naturally, the users get more suspicious about the post when received with fake-tag. Thus, we assume $p_x^u > p_y^u$ for any $u \in \{R, F\}$.

As said before, the user forwards the post to some/all of its friends after tagging. The number of shares depends on how attractive the post is, which we measure by $\eta^u \in (0,1)$. As argued in \cite{sui2023falsehood}, the design of fake-posts is deceptive and more appealing; therefore, we assume $\eta^F > \eta^R$. 

Let ${\mathcal F}$ be the number of friends of a typical user of the OSN and assume that  ${\cal F}$ is independent and identically distributed across various users. Let $\tau^+$ and $\tau^-$ be the usual limits, e.g., $\Cx(\tau^-) := \lim_{t \uparrow \tau} \Cx(t)$. When a wi-user receives a post with fake-tag and shares it with fake-tag, it generates $\xi_{xx, wi}$ number of fake-tagged copies. Similarly, when it tags the post as real, it shares to $\xi_{xy, wi}$ friends. Define $\xi_{yx, wi}$ and $\xi_{yy, wi}$ in a similar manner. We assume ($k$ is some constant):
\begin{eqnarray}
\label{eqn_shares_wi}
 \xi_{ix, wi}(\Om(\tau^-))  = \xi_{iy, wi}(\Om(\tau^-)) \sim Bin\left(\mathcal{F},\eta^u + \frac{k}{(\Sa(\tau^-))^2}\right) \mbox{ for } i \in \{x, y\},
\end{eqnarray}$Bin(\cdot, \cdot)$ denotes a binomial random variable; many times, users receive the post more than once, however, they may not be interested in it again - thus, the new effective shares in \eqref{eqn_shares_wi} reduces with the total copies/shares of the post generated so far, i.e., $\Sa(\tau^-) := \Ax(\tau^-) + \Ay(\tau^-)$, for example as in \eqref{eqn_shares_wi}. The distribution considered in \eqref{eqn_shares_wi} is a specific example; however, our analysis can extend to any total-current shares-dependent sharing-distribution that satisfies assumption \ref{a2_prop} (see Section \ref{sec_prop_BP}).



\subsection{Warning-seeking (ws) users}
These users  click on the warning button also - they incorporate the sender's tag, their innate capacity and the warning provided by the OSN to decide the tag.

Say a ws-user views the fake-tagged post at time $\tau$. Let $\omega_\tau$ be the warning at this time. Then, as in \cite{kapsikar2020controlling}, we assume that such user tags the post as fake (real) w.p. $\min\{\alpha_x^u \omega_\tau, 1\}$ (respectively, $1-\min\{\alpha_x^u \omega_\tau, 1\}$) before sharing; here, $\alpha_x^u > 0$ is the sensitivity parameter to the warning when the post is received with fake-tag. Similarly, if the post received by the ws-user has a real-tag, then it tags the post as fake or real w.p. $\min\{\alpha_y^u \omega_\tau, 1\}$ and  $1-\min\{\alpha_y^u \omega_\tau, 1\}$, respectively, where $\alpha_y^u > 0$ is the sensitivity parameter when the post is received with real-tag. Thus, we have:
\begin{align}\label{eqn_tag_ws}
    P(I_{x, ws}(\Om(\tau^-))=1|{\cal G}_\tau) = \min\{\alpha_x^u \omega_\tau, 1\} \mbox{ and }  P(I_{y, ws}(\Om(\tau^-))=1|{\cal G}_\tau) = \min\{\alpha_y^u \omega_\tau, 1\}.
\end{align}
The \textit{sensitivity parameters are indicative of the user's intrinsic ability to recognize the actuality of the post}. These parameterize \textit{warning-aided identification}, while  $p_F^u, p_R^u$ are the probabilities of \textit{un-aided identification}; both are characteristics of the users of the OSN. We thus assume a linear dependence between the two as in \cite{agarwal2023single}, i.e., we assume a $\rho \in (0,1)$ such that 
\begin{align}\label{eqn_prob_wi}
p_F^u = \alpha_x^u \rho \mbox{ and } p_R^u = \alpha_y^u \rho.
\end{align}

Now, similar to wi-users, a ws-user also shares the post with its friends. Using notations as in \eqref{eqn_shares_wi}, we have ($k$ is some constant):
\begin{eqnarray}
\label{eqn_shares_ws}
 \xi_{ix, ws}(\Om(\tau^-)) = \xi_{iy, ws}(\Om(\tau^-)) \sim Bin\left(\mathcal{F},\eta^u + \frac{k}{(\Sa(\tau^-))^2}\right) \mbox{ for } i \in \{x, y\}.
\end{eqnarray}

\subsection{Adversaries (a)} As is usually the case, there can be a small fraction of adversarial users on the OSN. These users aim to harm the efficacy of the system-generated warning by incorrectly tagging the post. Their agenda for doing so can be in self-interest or political. Often, such users do not have prior information about the actuality of the post, but to meet their objective they target to confuse the users about the actuality of the posts. Towards this, we consider that they always tag any post as real. In a way, such users are the ones who wish to color (tag) the posts as real, irrespective of the actuality of the posts. 

Let $I_{x, a}(\Om(\tau^-))$ and $I_{y, a}(\Om(\tau^-))$ be the indicator that an a-user with a fake or real-tagged copy of the post tags the post as fake, where $\tau$ is the time when an a-user views the post. Here, we have:
\begin{align}\label{eqn_tag_rc}
    P(I_{x, a}(\Om(\tau^-))=1|{\cal G}_\tau) = P(I_{y, a}(\Om(\tau^-))=1|{\cal G}_\tau) = 0. 
\end{align}

An adversarial user shares the post with a real-tag to its friends with probability $\eta_a \in (0,1)$, irrespective of the attractiveness of the post. Therefore, we have ($k$ is some constant):
\begin{align}\label{eqn_shares_adv}
    \xi_{ix, a}(\Om(\tau^-)) \equiv 0 \mbox{ and } \xi_{iy, a}(\Om(\tau^-)) \sim Bin\left(\mathcal{F},\eta_a  + \frac{k}{(\Sa(\tau^-))^2}\right) \mbox{ for } i \in \{x, y\}.
\end{align}

\subsection{Non-participants (np)} In \cite{kapsikar2020controlling}, it is assumed that all users viewing the post share and tag it. In reality, there can be users named as non-participants who neither participate in the tagging process nor share the post. In other words, when they receive a copy of the post, they do not respond, which we capture as:
\begin{align}
    P(I_{i, np}(\Om(\tau^-))=1|{\cal G}_\tau) &= P(I_{i, np}(\Om(\tau^-))=1|{\cal G}_\tau) = 0, \label{eqn_tag_np}
\end{align}
and shares to none, i.e.,
\begin{align}
    \xi_{ix, np}(\Om(\tau^-)) &= \xi_{iy, np}(\Om(\tau^-)) \equiv 0, \mbox{ for } i \in \{x, y\}. \label{eqn_share_np}
\end{align}

\vspace{2mm}

\noindent \textbf{Number of shares:} Let ${\cal U} := \{\mbox{wi, ws, a, np}\}$ be the set of types of users in the system. Let $\mu_0, \mu_1, \mu_2, \mu_a$ be the respective proportions of np, wi, ws, a-users on the OSN such that $\mu_1 + \mu_2 + \mu_a + \mu_0 = 1$; \textit{we assume that the OSN knows these proportions}. Since our approach is based on crowd-signals, therefore, it is meaningful to assume that $\mu_2 \in (0, 1)$. Any user of the OSN visits it after a random time which is exponentially distributed with parameter $1$ (without loss of generality); this is a commonly made assumption in the literature (see, for example, \cite{agarwal2021co, dhounchak2023viral, van2010viral}). If required, one can model different users visiting the OSN at different rates, for example, a-users might visit more often; our framework can easily extend to such a case. Any user of $j$-type, after viewing the post with fake-tag ($i = x$) or real-tag ($i = y$), generates $\Gamma_{ix, j}$ and $\Gamma_{iy, j}$ number of new fake and real-tagged copies of the post respectively, where:
\begin{align}\label{eqn_final_shares}
\begin{aligned}
    \Gamma_{ix, j}(\Om(\tau^-)) &:= I_{i, j} (\Om(\tau^-)) \xi_{ix, j} (\Om(\tau^-)), \mbox{ and }\\
    \Gamma_{iy, j}(\Om(\tau^-)) &:= \bigg(1-I_{i, j}(\Om(\tau^-))\bigg) \xi_{iy, j}  (\Om(\tau^-)), \mbox{ for } i \in \{x, y\} \mbox{ and } j \in {\cal U}.
\end{aligned}
\end{align}

Next, we discuss some meaningful assumptions (inspired by \cite{kapsikar2020controlling}).

\vspace{2mm}
\noindent \textbf{Regime of parameters and assumptions:} The probability of a user fake-tagging any $u$-post is higher when the sender's tag is fake, thus, $\alpha_x^u > \alpha_y^u$, for $u \in \{R, F\}$. We assume that the users are more likely to tag fake-post as fake, as compared to tagging real-post as fake, irrespective of sender's tag, i.e., $\alpha_i^F > \alpha_i^R$, for each $i \in \{x, y\}$. \revg{Since the intent of a-users is to share the post rigorously, therefore, we assume $\eta_a > \eta^u$, for each $u$, only in the numerical experiments; the theoretical results follow even if $\eta_a \leq \eta^u$. }
Thus, in all, we assume the following:
\begin{align}\label{eqn_relation_parameters}
\begin{aligned}
    \alpha_x^u &> \alpha_y^u > 0 , \mbox{ for each }  u \in \{R, F\},   \alpha_i^F > \alpha_i^R \mbox{ for each } i \in \{x, y\},\\
    \eta_a &> \eta^F > \eta^R > 0,
    \mu_2 \in (0,1) \mbox{ and } \rho \in (0,1).
\end{aligned}
\end{align}
For the sake of clarity, we summarize all the notations which will be used consistently throughout the paper:
\begin{table}[htbp]
\centering
\scalebox{1}{
\begin{tabular}{|c|c|l|}
\hline
Sr. No. & Notation                        & \multicolumn{1}{c|}{Description}                                \\ \hline
1.      & ${\cal U} = \{\mbox{wi, ws, a, np}\}$                   & types of users: warning-ignoring, warning-seeking, adversarial, non-participating \\ \hline
2.      & $\mu_0, \mu_1, \mu_2, \mu_a$           & proportion of np, wi, ws and a-users \\ \hline
3.      & $u \in \{R, F\}$ & actuality of the post as real or fake respectively                                        \\ \hline
4.      & $\eta^u, \eta_a$                        & probability of a user/adversary sharing the post to its friend                                     \\ \hline
5.      & $x, y$                    & fake or real tag by the sender \\ \hline
6.      & $\alpha_x^u, \alpha_y^u$                    & sensitivity of a user towards the warning when received with fake or real tag  \\ \hline
\end{tabular}}
\caption{Summary of the notations} \label{table_notations}
\end{table}

\subsection{Warning Mechanism (WM) - system-generated warning}
\begin{figure}[http]
    \centering
    \includegraphics[trim = {0cm 3.3cm 0cm 0cm}, clip, scale = 0.5]{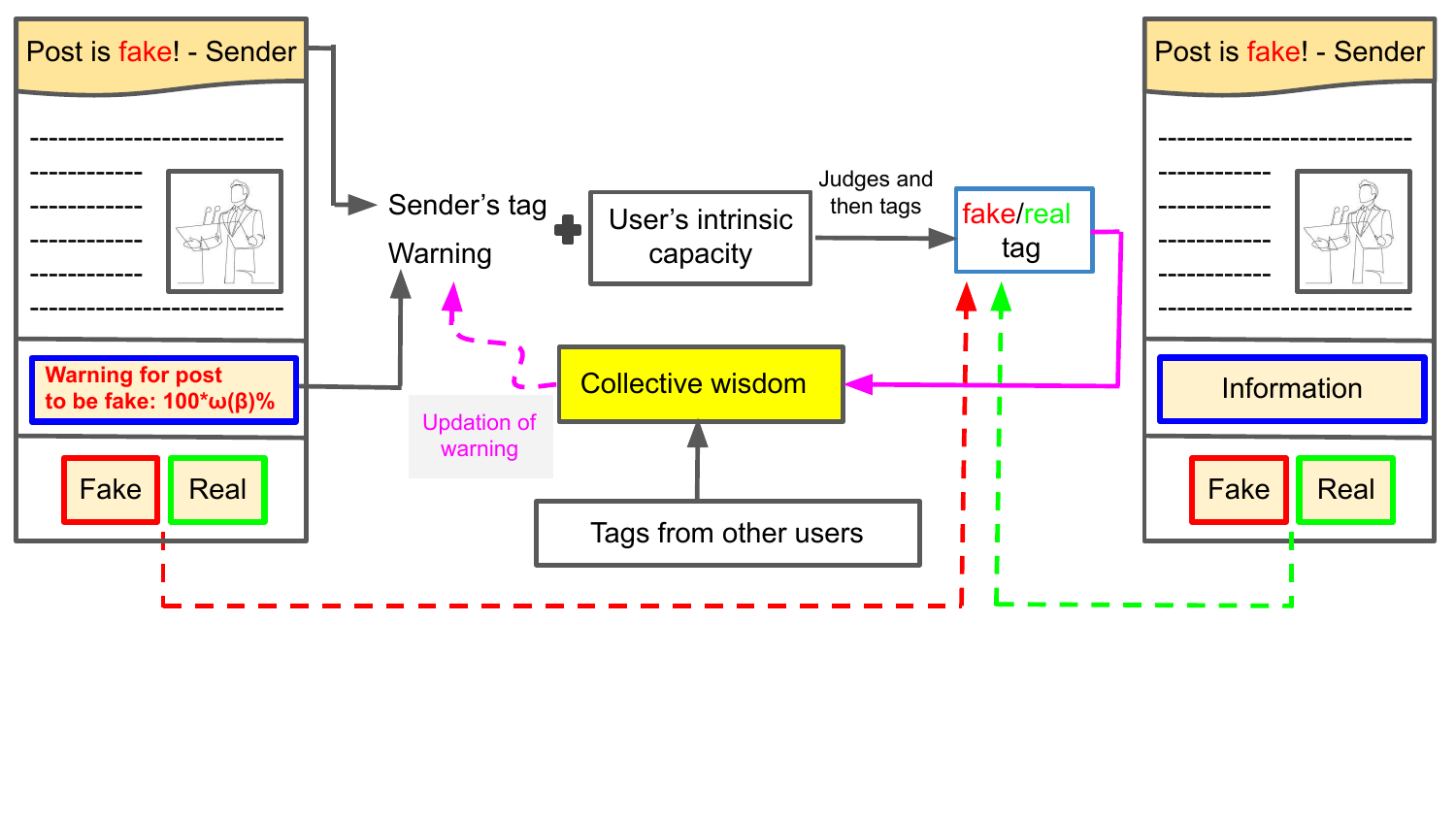}
    \caption{On the left, ws-user tags the post as fake. On the right, a-user tags the post as real, without checking the warning or sender's tag.}
    \label{fig_warning_mech}
\end{figure}

In \cite{kapsikar2020controlling}, the authors designed a warning mechanism (WM) by leveraging upon the responses of the users. They assumed all users are ws-users and did not consider the adversaries (i.e., $\mu_2 = 1$). The main idea behind the design of the mechanism is to exploit the collective wisdom of the users (via responses of all users), as depicted in Figure \ref{fig_warning_mech} (left side). 
The warning considered in \cite{kapsikar2020controlling} is:
\begin{equation} \label{eqn_warning}
 \omega_t = \bigg ( \frac{w\Cx(t)}{\Cx(t)+b\Cy(t)} + \gamma \bigg) = \bigg ( \frac{w\Beta(t)}{\Beta(t)+b(1 - \Beta(t))} + \gamma \bigg), \mbox{ where } \Beta(t) := \frac{\Cx(t)}{\Cx(t) + \Cy(t)}
\end{equation}represents the relative fraction of (unread) fake-tagged copies at time $t$; $w$ and $b$ are the control parameters; $\gamma > 0$ is the parameter which captures the prior knowledge OSN has about the post via some fact-check mechanism. Here, $w \in [0, \overline{w}]$ for $\overline{w} := \frac{1}{\alpha_x^F}-\gamma$. This ensures that a ws-user tags the fake-tagged copy of the post as fake with probability $\min\{\alpha_x^u \omega(\beta), 1\} =  \alpha_x^u \omega(\beta)$ for any $\beta \in [0,1]$, when the warning is as in \eqref{eqn_warning}; thus, $\min\{\alpha_y^u \omega(\beta), 1\} = \alpha_y^u \omega(\beta)$ (since $\alpha_y^u < \alpha_x^u$, see \eqref{eqn_relation_parameters}). Further, the parameter $b \in [0, \infty)$. The warning in \eqref{eqn_warning} is generated individually for each post.

In this paper, we are considering a variety of user behavior. Therefore, the warning is now influenced by the responses of users who ignore the warning while tagging or are purposely providing incorrect tags. 
In Figure \ref{fig_warning_mech}, we depict that the warning is updated by the response (fake) of the ws-user (left side of the figure) and also by that of a-user (right side of the figure). Similarly, one can visualize how a warning gets updated when a wi-user tags. This suggests that the warning \eqref{eqn_warning} needs to be studied for our complex and more realistic system.

It is clear from the discussion so far that the end goal of the OSN is to nudge users towards the correct identification of the posts. \revg{Let $B^u(t)$ represents the proportion of fake-tags, given that the actuality of the post is $u \in \{R, F\}$.} Then, similar to \cite{kapsikar2020controlling}, we aim to optimally choose $w, b$:
\begin{itemize}
    \item to maximize the proportion of fake-tags for the fake-post, $\max \lim_{t\to \infty} \Beta^F(t)$,  and
    \item to ensure that the proportion of fake-tags for the real-post, $\lim_{t\to \infty} \Beta^R(t)$, is at most $\delta$, for some $\delta \in (0,1)$.
\end{itemize}
\revr{The above objective is well defined if the limits in the above exist and are unique almost surely. By Theorem \ref{thrm_BP_to_fake} stated in section \ref{sec_warning_old}, we prove that the limits indeed exist (but need not be unique) for any general warning mechanism. Hence, define ${\cal L}^u := \left\{\lim_{t\to \infty} \Beta^u(t)\right\}$ as the set of all possible limits for $u$-post, across all sample paths, and consider the following optimization problem:
\begin{align}\label{eqn_gen_opt}
    \max_{w, b} \inf({\cal L}^F) \mbox{ subject to } \sup({\cal L}^R) \leq \delta.
\end{align}}

Further, we shall investigate the following two questions:
\begin{enumerate}
    \item How does the optimal WM in \eqref{eqn_warning} perform in the presence of wi-users and a-users? 
    \item If the performance degrades, can we design improved WMs which are robust against adversaries?
\end{enumerate}

\subsection{Warning dynamics and Branching process}
It is clear that when a user tags the post as fake, the fake number of copies (represented by $x$)  gets updated; otherwise, the real ($y$) number of copies gets updated. Further, the user who receives the post can be one among the type $i$, for $i \in {\cal U}$, w.p. given by the proportion of the type it belongs to; for example, the recipient can be a wi-user w.p. $\mu_1$. As discussed in \eqref{eqn_shares_wi}, \eqref{eqn_shares_ws}, \eqref{eqn_shares_adv} and \eqref{eqn_share_np}, the distribution of the number of shares depends on the type of the user who received the post.

Let $\tau$ be the time when a type-$i$ user views the post on its timeline with a fake-tag. Then, the number of fake-tagged and real-tagged copies of the underlying post evolves at time $\tau$ as follows:
\begin{equation}\label{eqn_transition_fake_tag}
\begin{aligned}
\Cx(\tau^+) & = \Cx(\tau^-) - 1 + \Gamma_{xx, i}(\Om(\tau^-)), 
\Cy(\tau^+) = \Cy(\tau^-) + \Gamma_{xy, i}(\Om(\tau^-)),\\ 
\Ax(\tau^+) & = \Ax(\tau^-) + \Gamma_{xx, i}(\Om(\tau^-)), \mbox{ and }
\Ay(\tau^+) = \Ay(\tau^-) + \Gamma_{xy, i}(\Om(\tau^-)).
\end{aligned}
\end{equation}
We argued before that once a user reads a post, it is seldom interested in the same post again; thus, the current (unread) number of fake-tagged copies decreases by $1$. Similarly, when a type-$i$ user who received the post with the real-tag views the post,  the system evolves as:
\begin{equation}\label{eqn_transition_real_tag}
\begin{aligned}
\Cx(\tau^+) & = \Cx(\tau^-) + \Gamma_{yx, i}(\Om(\tau^-)), 
\Cy(\tau^+) = \Cy(\tau^-) - 1 + \Gamma_{yy, i}(\Om(\tau^-)),\\
\Ax(\tau^+) & = \Ax(\tau^-) + \Gamma_{yx, i}(\Om(\tau^-)), \mbox{ and }
\Ay(\tau^+) = \Ay(\tau^-) + \Gamma_{yy, i}(\Om(\tau^-)).
\end{aligned}
\end{equation} 
We shall briefly call the above warning-mechanism aided dynamics as \underline{warning dynamics}. At this point, it is important to state that the dynamics described above can be modelled as a continuous-time total-current population-dependent branching process (TC-BP) of \cite{agarwal2021new}, except for varying death rates. We will discuss how such correspondences can be made in Section \ref{sec_warning_old}; in particular, we will see that the viewing of the post can be modelled as a death in an appropriate TC-BP and hence, will have different death-types and rates owing to different types of users. However, we first analyze the TC-BPs with multiple death types in the next section using ODE based stochastic approximation technique, which will be instrumental for our study.

\vspace{2mm}
\noindent \textbf{Informal outline for design of WMs:}
\revr{
We consider any general warning mechanism $\omega(\beta)$, which depends on the proportion of fake-tags ($\beta$) provided by the previous recipients of the post. The limiting behaviour of the warning-guided post-propagation process is analyzed using the ODE derived via the analysis of the underlying BP. In particular, we will show that the analysis of a one-dimensional ODE suffices to study the limits of the underlying process; of course, the limits depend upon the warning mechanism utilized. 
The main idea is to reverse-engineer:  consider the design of the warning mechanism (to achieve the desired output), based  on the anticipated attractors of the one-dimensional ODE. We will follow this approach in section \ref{sec_warning_old} and thereafter, where we bring our attention back to the control of fake-post propagation over OSNs.}

\section{Total-Current population-dependent Branching Process (TC-BP) with multiple death types} \label{sec_prop_BP}

Consider two types of populations, namely $x$ and $y$-types, and  
let $c_{x, 0}$ and $c_{y, 0}$ be their respective initial sizes. An individual can either die naturally, or it may die  differently due to unnatural circumstances. We refer any death which is not natural as `unnatural death'\footnote{
In biological systems, unnatural deaths may occur due to exposition to a virus, competition with other species, etc. We discuss unnatural deaths for the application at hand in section \ref{sec_warning_old}.}. Let $D_i := \{0, 1, \dots, d_i\}$ be the set of variety of deaths for $i$-type individual, where $d_i \in [0, \infty)$. Here, $d = 0$ represents the natural death and $d \in D_i-\{0\}$ represents an unnatural death; $D_x$ need not equal $D_y$ as some circumstances may affect only one population. We shall briefly refer to the death of variety $d$ as $d$-death.

Now, given that the interest of this paper is in controlling the fake post propagation over OSNs, our focus is on the time-asymptotic proportion of the population (fake-tags). Therefore, it is sufficient to study the embedded process (discrete-time chain defined at death instances) of the continuous-time Markov process. In \cite{agarwal2021new}, the authors analysed the TC-BP using stochastic approximation based approach, where only natural deaths occur. In this section, we will follow the same approach to incorporate different varieties of deaths. We begin by introducing few notations which are exactly as in \cite{agarwal2021new}, however are re-written here for the ease of reading.

Let $\tau_n$ be the time at which $n$-th individual dies. Consider any $n \geq 1$. Let $\Om_n := (C_{x, n}, C_{y, n}, T_{x, n}, T_{y, n})$, where $C_{x, n}, C_{y, n}$ represent the \textit{current population} and $T_{x, n}, T_{y, n}$ are the \textit{total population} sizes immediately after $\tau_n$, e.g., $C_{x, n} = \Cx(\tau_n^+)$. Let $S_n :=  C_{x, n} + C_{y, n}$ be the sum current population, again immediately after $\tau_n$.  Let $\om = (\cx, \cy, \ax, \ay)$ be a realisation of the random vector $\Om$. \revg{Any individual can die naturally or unnaturally. We assume that the time till $d$-death of an $i$-type individual is exponentially distributed with parameter $\lambda_{i, d} \in (0, \infty)$. An individual in the population will die according to the first death (variety) event that occurs. By memoryless property, after any given instance of time (e.g., $\tau_n$), the death-time of any $i$-type individual in the population is again exponentially distributed with parameter $\sum_{d} \lambda_{i, d}$, and hence the first death in the two populations is exponentially distributed with parameter $\left(\sum_{d} \lambda_{x, d} + \sum_{d} \lambda_{y, d}\right)$. We further assume that the parameter $\lambda_{i, d}$ depends on the population-size, i.e., $\lambda_{i, d}(\om_n)$, conditioned on $\om_n$, for each $i \in \{x, y\}$. Observe that we have population-dependency even for the natural deaths, in contrast to the classical models studying only population-independent natural deaths (see, for example, \cite{klebaner1993population, athreya2012classical, jagers1969proportions}).}


The current population can get extinct, and thus let $\nu_e := \inf \{n : S_n = 0\}$ be the extinction epoch, with the usual convention that $\nu_e = \infty$, when $S_n > 0$ for all $n$. \textit{For the sake of completion, define $\Om_n := \Om_{\nu_e}$ and $\tau_{n} :=\tau_{\nu_e}$, for all $n \geq \nu_e$, when $\nu_e < \infty$.} \revg{We refer the sample paths in which $\nu_e = \infty$ as the non-extinction paths, and the complementary ones as the extinction paths.} \textit{Define $\Beta_n :=  C_{x, n}/S_n$ as the proportion of $x$-type population among current population}. Let $\beta = \cx/(\cx+\cy)$ be a realisation of $\Beta$. 

\hide{
Define $\Om(t) := (\Cx(t), \Cy(t), \Ax(t), \Ay(t))$, where $\Cx(t), \Cy(t)$ represent the \textit{current population} and $\Ax(t), \Ay(t)$ are the \textit{total population} sizes at time $t$. Observe 
$(\Ax(0), \Ay(0)) = (\cx_0, \cy_0)$.  Let $\om = (\cx, \cy, \ax, \ay)$ be a realisation of the random vector $\Om$. Now, conditioned on $\om$, say an $i$-type individual lives for an exponentially distributed with parameter $\lambda_{i, d}(\om) \in (0, \infty)$ before it $d$-dies, where $i \in \{x, y\}$. In classical BPs, such parameters are population-independent, while here, we consider them to depend not only on the population sizes, but also on the variety of death. 
Here, $\lambda_{i, 0}(\om)$ represents the lifetime parameter if populations were living (and reproducing) independently.

Define $\beta := \cx/(\cx + \cy)$ as the proportion of $x$-type individuals, conditioned on $\om$. }

\subsection{Evolution of embedded process}
In classical BPs, each individual lives for a random time which is exponentially distributed with a common parameter  (say) $\lambda > 0$. Thus, an individual to die at $n$-th epoch is of $x$-type w.p.\footnote{This happens due to the memory-less property of exponential distribution and as minimum of $k$ independent and identically distributed exponentially distributed random variables with parameter $\lambda$ is exponentially distributed with parameter $k\lambda$. } $\beta_n$, conditioned on $\Om_n = \om_n$. In similar lines, with the possibility of unnatural deaths, the probability that an $i$-type individual $d$-dies is given by:
\begin{align}\label{eqn_prob_d_death}
\begin{aligned}
\mbox{P($x$-type individual $d$-dies$|\om$)} &= \frac{\lambda_{x, d}(\om) \beta}{d(\om)} \mbox{ and} \\
\mbox{P($y$-type individual $d$-dies$|\om$)} &= \frac{\lambda_{y, d}(\om) (1-\beta)}{ d(\om)}, \mbox{ where }
d(\om) := \sum_{d \in D_x} \lambda_{x, d}(\om) \beta + \sum_{d \in D_y} \lambda_{y, d}(\om) (1-\beta).
\end{aligned}
\end{align}
In all, the overall probability that an $i$-type individual is the first to die after previous death instance, $\tau$, is given by:
\begin{align}\label{eqn_prob_death}
\begin{aligned}
\mbox{P($x$-type individual dies$|\om$)} &= \frac{\beta \sum_{d \in D_x} \lambda_{x, d}(\om)}{d(\om)} =: f_{\beta}(\om) \mbox{ and } \mbox{P($y$-type individual dies$|\om$)} = 1 - f_{\beta}(\om).
\end{aligned}
\end{align}

Say an individual of $i$-type dies at $n$-th epoch. Then, the current size (not the total size) of $i$-type reduces by $1$ due to death. Further, if it $d$-dies for $d \in D_i$, it produces $\offs_{ii, d}(\Om_{n-1})$ and $\offs_{ij, d}(\Om_{n-1})$ offspring of $i$-type and $j$-type ($j\neq i$) respectively, conditioned on the sigma algebra $\sigma\{\Om_{n-1}\}$, where $\offs_{ij, d} (\Om_{n-1})$ is an  integer-valued random variable. Basically, when $\Om_{n-1} = \om_{n-1}$, the random offspring are represented by $\offs_{ij, d}(\om_{n-1})$ for each $i, j$ and $d$. Thus, the embedded process immediately after an $i$-type individual $d$-dies at $n$-th epoch is given by:
\begin{equation}\label{evolve_x_up_time_prop}
\begin{aligned}
C^i_n &= C^i_{n-1}  + \offs_{ii, d}(\Om_{n-1}) - 1, \ \ \  A^i_n = A^i_{n-1}  + \offs_{ii, d}(\Om_{n-1}), 
\\
C^j_n &= C^j_{n-1} + \offs_{ij, d}(\Om_{n-1}), \ \ \ \ \ \ \ \ \ A^j_n = A^j_{n-1} + \offs_{ij, d}(\Om_{n-1}), \mbox{ for } i \neq j.
\end{aligned}
\end{equation}

Now, conditioned on $\om$, we assume the $\om$-dependent random offspring satisfy the following, which also ensures throughout super-criticality, a notion defined in \cite{agarwal2021new}:
\begin{enumerate}[label=\textbf{C.\arabic*}, ref=\textbf{C.\arabic*}]
    \item \label{a1_prop} There exist two integrable random variables $\overline{\offs}$ and $\underline{\offs}$ which bound the random offspring as: $0 \leq \underline{\offs} \leq \offs_{ix, d}(\om) + \offs_{iy, d}(\om) \leq \overline{\offs}$ almost surely (a.s.), for each $\om$, for each $d$. Also,  $E[\overline{\offs}^2] < \infty$ and $E[\underline{\offs}] > 1$. Further, $\offs_{ii, d}(\om) \geq 0$ a.s., for each $i, \om, d$. Furthermore, assume that $\inf_{\om} \lambda_{x, d}(\om) > 0$ for each $d \in D_x$ and $\inf_{\om} \lambda_{y, d}(\om) > 0$ for each $d \in D_y$.
\end{enumerate}

\subsection{Mean matrix}

Let \textit{$m_{ij, d} (\om) := E[ \offs_{ij, d} (\om) ]$ denote the expectation of the number of $j$-type offspring, when an $i$-type parent $d$-dies, conditioned on $\om$, for $i, j \in \{x, y\}$ and $d \in D_i$}. Further, define the mean matrix $M(\om) := [m_{ij}(\om)]_{i, j\in\{x, y\}}$ as given below:
\begin{eqnarray}\label{eqn_mean_matrix}
M(\om) :=
\left [  
\begin{array}{ll}
    \frac{\sum_{d \in D_x} \lambda_{x, d}(\om) m_{xx, d}(\om)}{\sum_{d \in D_x} \lambda_{x, d}(\om)}    & \frac{\sum_{d \in D_x} \lambda_{x, d}(\om) m_{xy, d}(\om)}{\sum_{d \in D_x} \lambda_{x, d}(\om)}   \\ \\
      \frac{\sum_{d \in D_y} \lambda_{y, d}(\om) m_{yx, d}(\om)}{\sum_{d \in D_y} \lambda_{y, d}(\om)}    & \frac{ \sum_{d \in D_y} \lambda_{y, d}(\om) m_{yy, d}(\om)}{\sum_{d \in D_y} \lambda_{y, d}(\om)}
\end{array}
\right ].  
\end{eqnarray}
Then,  for $j \in \{x, y\}$, we have (see \eqref{eqn_prob_d_death}, \eqref{eqn_prob_death} and \eqref{eqn_mean_matrix}):

\vspace{-4mm}
{\small
\begin{align}
\begin{aligned}
    E[j\mbox{-type offspring produced by an } x\mbox{-type parent}|\om] &= \sum_{d \in D_x} \frac{\lambda_{x, d}(\om) \beta}{d(\om)} m_{xj, d}(\om) = f_{\beta}(\om) m_{xj}(\om),\\
    E[j\mbox{-type offspring produced by a } y\mbox{-type parent}|\om] &= \sum_{d \in D_y} \frac{\lambda_{y, d}(\om) (1-\beta)}{d(\om)} m_{yj, d}(\om)= (1-f_{\beta}(\om)) m_{yj}(\om).
\end{aligned}
\end{align}}

As in \cite[Lemma 2, Appendix A]{agarwal2021new}, one can prove the dichotomy for the sum current population of TC-BP with multiple death types, as in the following:
\begin{lem}\label{lemma_dichotomy}
Assume \ref{a1_prop} and define $\underline{m} =: E[\underline{\offs}]$. Then, we have:
$$
P\left(\left\{\liminf_{n} S_n e^{-\underline{\lambda}(\underline{m}-1)\tau_n} > 0\right\} \cup \left\{\lim_{n \to \infty} S_n = 0\right\}\right) = 1, 
$$where $\underline{\lambda} := \inf_\om\{\lambda_{x, 0}(\om), \dots, \lambda^x_{d_x}(\om), \lambda_{y, 0}(\om), \dots, \lambda^y_{d_y}(\om)\} > 0$.
\end{lem}
\revg{Thus, the sum current population either gets extinct or in non-extinction paths, it explodes (i.e., it grows exponentially larger at rate $\underline{\lambda}(\underline{m}-1)$).
}

\subsection{Main Result}
We will now provide the first main result of the paper which determines the limit proportion, $\lim_{t \to \infty} \Bc(t)$ in non-extinction paths and additionally, provides the deterministic approximate trajectories for the underlying BP. The result follows in similar lines to \cite[Theorem 1]{agarwal2021new}, while accommodating some important changes for multiple deaths. As established in Lemma \ref{lemma_dichotomy}, the underlying BP can explode. In such a case, it is a common practice to scale the process appropriately that enables convergence to a finite limit (see, for example, \cite{agarwal2021new, kapsikar2020controlling}).



\hide{Let $\tau_n$ be the time at which $n$-th individual dies.
Consider any $n \geq 1$. Let $\Om_n := (C_{x, n}, C_{y, n}, T_{x, n}, T_{y, n})$ be the individual (current and total) populations and $S_n :=  C_{x, n} + C_{y, n}$ be the sum current population, both immediately after $\tau_n$, e.g., $C_{x, n} = \Cx(\tau_n^+)$. The current population can get extinct, and thus let $\nu_e := \inf \{n : S_n = 0\}$ be the extinction epoch,  
with the usual convention that $\nu_e = \infty$, when $S_n > 0$ for all $n$. \textit{For the sake of completion, define $\Om_n := \Om_{\nu_e}$ and $\tau_{n} :=\tau_{\nu_e}$, for all $n \geq \nu_e$, when $\nu_e < \infty$.}}

To this end, define the scaled ratios $\Pc_n := S_{n}/n$ and $\Tc_n := C_{x, n}/n$. Let $Z_n := T_{x, n} + T_{y, n}$ be the total population size immediately after $\tau_n$, and then analogously, define  $\Pa_n$ and $\Ta_n$ for the total population.  Let $\Ups_n := (\Pc_n, \Tc_n, \Pa_n, \Ta_n)$, and let $\Ups_0 := (s_0^c, c_{x, 0}, s_0^c, c_{x,0})$ denote the initial population, where $s_0^a = s_0^c := c_{x, 0} + c_{y, 0}$. \hide{\textit{Define $\Bc_n := \Tc_n/\Pc_n  = C_{x, n}/S_n$ as the proportion of $x$-type population among current population}. } Let $\ups := (\pc, \tc, \pa, \ta)$ be a realisation of $\Ups$. \hide{ and $\beta = \cx/(\cx+\cy) = \tc/\pc$ be that of $\Bc$. }

In \cite[\textbf{A.2}]{agarwal2021new}, the authors assumed that the total-current population-dependent  mean functions converge to proportion dependent mean functions, which can further be discontinuous. Similar to that, we assume that the resultant mean functions ($m_{ij}(\om)$, and not $m_{ij, d}(\om)$) converge to proportion-dependent mean functions at a certain rate. However, to accommodate for the variety of deaths, we assume that the lifetime parameters of the populations also become proportion-dependent asymptotically (at the same rate of convergence as that of the mean functions).
\begin{enumerate}[label=\textbf{C.\arabic*}, ref=\textbf{C.\arabic*}]
\setcounter{enumi}{1}
    \item Define $\beta(\ups) := \tc/\pc = \cx/s$. As sum current population,  $s \to \infty$:   \label{a2_prop}
    \begin{align*}
     |m_{ij}(\om) - \minf_{ij}(\beta(\ups))| &\leq \frac{1}{(s)^\alpha}, \mbox{ for each }i, j \in \{x, y\} \mbox{ and}\\
    |\lambda_{i, d}(\om) - \lambda_{i, d}^\infty(\beta(\ups))| &\leq \frac{1}{(s)^\alpha}, \mbox{ for each } d \in D_i \mbox{ for each } i \in \{x, y\},    \mbox{ for some } \reva{\alpha \ge 1}.
    \end{align*}
\end{enumerate}
Further, under \ref{a2_prop}, the function $f_\beta(\om)$ converges to $f_{\beta}^\infty(\beta)$ as given below (see \eqref{eqn_prob_death}):
\begin{align}\label{eqn_f_infty}
\begin{aligned}
|f_{\beta}(\om) - f_{\beta}^\infty(\beta) | &\leq \frac{1}{(s)^\alpha}, \mbox{ where } f_{\beta}^\infty(\beta) := \frac{\beta \sum_{d \in D_x} \lambda_{x, d}^\infty(\beta)}{d^\infty(\beta)} \mbox{ with}\\
d^\infty(\beta) &:= \beta \sum_{d \in D_x} \lambda_{x, d}^\infty(\beta) + (1-\beta) \sum_{d \in D_y} \lambda_{y, d}^\infty(\beta).
\end{aligned}
\end{align}
In all, under \ref{a1_prop}-\ref{a2_prop}, we analyze the ratios $\Ups_n$ using the solutions of the following ODE:
\begin{align}\label{eqn_ODE_prop}
\dot{\ups} &= \ga(\ups) = \mathbf{h}(\beta)1_{\{\pc > 0\}} - \ups, \mbox{ where } \mathbf{h}(\beta) := (h_{\psi}^{c},  h_{\theta}^{c},  h_{\psi}^{a},  h_{\theta}^{a}), \mbox{ with} \nonumber \\
        h_{\psi}^{c}(\beta) &= f_{\beta}^\infty(\beta) \bigg(\minf_{xx}(\beta) +     \minf_{xy}(\beta)\bigg) + \left(1-f_{\beta}^\infty(\beta)\right)\bigg(\minf_{yy}(\beta) + \minf_{yx}(\beta)\bigg) - 1,  \nonumber\\
        h_{\theta}^{c}(\beta)  &= f_{\beta}^\infty(\beta)\bigg(\minf_{xx}(\beta) - 1\bigg) + \left(1-f_{\beta}^\infty(\beta)\right) \minf_{yx}(\beta), \\
        h_{\psi}^{a}( \beta) &=  f_{\beta}^\infty(\beta) \bigg(\minf_{xx}(\beta) + \minf_{xy}(\beta) \bigg) + \left(1-f_{\beta}^\infty(\beta)\right)\bigg(\minf_{yy}(\beta) + \minf_{yx}(\beta) \bigg)   \mbox{ \normalsize  and}  \nonumber  \\
        h_{\theta}^{a}( \beta)  &= f_{\beta}^\infty(\beta)\minf_{xx}(\beta) + \left(1-f_{\beta}^\infty(\beta)\right) \minf_{yx}(\beta).  \nonumber
\end{align}
Now, exactly as in \cite[\textbf{A.3}]{agarwal2021new}, we assume the following (see \cite[Definition 1]{agarwal2021new} for the definition of extended solution):
\begin{enumerate}[label=\textbf{C.\arabic*}, ref=\textbf{C.\arabic*}]
\setcounter{enumi}{2}
    \item There exists a unique solution $\ups(\cdot)$ for ODE \eqref{eqn_ODE_prop} in the extended sense over any bounded interval.\label{a3_prop}
\end{enumerate}
As per \cite[Definition 2]{agarwal2021new}, let $\cA$ be the attractor set and $\cR$ be the saddle set with respect to the ODE \eqref{eqn_ODE_prop}. For systems modelling the BPs,  the following subset of the combined domain of attraction of $\cA$ and $\cR$ is relevant (recall the definition of ratios $\ups$):
\begin{align}\label{eqn_domain of attraction_prop}
    \cD &:= \{\ups \in (\mathbb{R}^+)^4: \tc \leq \pc \leq \pa, \ta \leq \pa \mbox{ and } \ups(t) \to {\cA}\cup \cR \mbox{ as } t \to \infty, \mbox{ if } \ups(0) = \ups\}.
\end{align}
Therefore, we will be interested in initial conditions $\ups(0) \in \cD_I$ for the ODE \eqref{eqn_ODE_prop}.

In \cite[Definition 4]{agarwal2021new}, we introduced a new notion of limiting behavior of the stochastic process, named `hovering around the saddle set' - here, the stochastic trajectory visits every neighborhood of $\cR$ infinitely often (i.o.), but also leaves some neighborhood of $\cR$ i.o. The main result in \cite{agarwal2021new} states that the random trajectory either converges to the attractor set or it converges to/hovers around a special kind of saddle set.  In particular, if any non-zero saddle point, $\ups^* \neq \mathbf{0}$, is attracted exponentially to $\cR$ along a particular affine sub-space, $\revg{{\mathbb S}(\ups^*) := \{\ups: \beta(\ups) = \beta(\ups^*)\}}$ and to $\cA$ in the remaining space, then such $\ups^*$ are named as (quasi) q-attractor in \cite{agarwal2021new}.  We have a similar result for the case with multiple deaths.

Similar to \cite{agarwal2021new}, 
under above definition, we finally assume the following:
\begin{enumerate}[label=\textbf{C.\arabic*}, ref=\textbf{C.\arabic*}]
\setcounter{enumi}{3}
    \item Let $\cA\cap\cD_I$ be the attractor set and each $\ups \in \cR\cap\cD_I$ be the q-attractor. Consider $\cD$ as in \eqref{eqn_domain of attraction_prop} and let  $\cS := \cD \cap \{\pa \leq b\}$, for some $b > 0$, be a compact subset of combined domain of attraction. Assume $p_{b} := P(\mathcal{V}) > 0$, where $\mathcal{V} :=  \{\omega : \Ups_n(\omega) \in \cS \mbox{ i.o.}\}$. \label{a4_prop}
\end{enumerate}
We have the following result:
\begin{theorem}\label{thm1}
  Under \ref{a1_prop}-\ref{a3_prop}, we have:
  \begin{enumerate}[label=(\roman*)]
        \item For every $T>0$, a.s.  there exists a sub-sequence $(n_l)$ such that:
            $$
            \sup_{k: t_k \in [t_{n_l}, t_{n_l} + T]} d(\Ups_k, \ups(t_k - t_{n_l})) \to 0  \mbox{ as } l \to \infty, \mbox{ where } t_n := \sum_{k=1}^n \frac{1}{k} \mbox{ and}
            $$
        $\ups(\cdot)$ is the extended solution of ODE \eqref{eqn_ODE_prop} which starts at $\ups(0) =
        \lim_{n_l \to \infty} \Ups_{n_l}$.
        \item Further, assume \ref{a4_prop}. Then, $P({\cal C}_1 \cup {\cal C}_2) \geq p_b$, where
        \begin{align*}
            \hspace{5mm} {\cal C}_1 &: =\{\Ups_n \to (\cA \cup \cR)\cap \cD_I \mbox{ as } n \to \infty\}, \mbox{ and }
                {\cal C}_2 := \{ \Ups_n \mbox{ hovers around } \cR \}. \hspace{4mm}  \mbox{ \eop} 
        \end{align*}
  \end{enumerate}
\end{theorem}
\textit{The proof of the above Theorem and all forthcoming results will be provided in \ref{appendix_journal2}.}

\subsection{Derivation of attractor and saddle sets}
It is evident from Theorem \ref{thm1} that the limit proportion, $\lim_{n \to \infty} \Beta_n$ can be deduced if one derives the attractor and saddle (specifically, q-attractor) sets. 
In \cite{agarwal2021new}, the authors proposed a procedure to derive these sets for the ODE \eqref{eqn_ODE_prop}, when only natural deaths occur. The main idea was to  exploit the dependence of limit mean functions on $\beta$ as in \ref{a2_prop} and finally, it is showed that the analysis of $\beta(\ups)$-ODE suffices. We extend the same approach for the new process with both natural and unnatural deaths. Towards this, one can derive the following limit $\beta$-ODE, using \eqref{eqn_ODE_prop}:
\begin{align}\label{eqn_beta_ODE_prop}
\begin{aligned}
\dot{\beta} &= \frac{1}{\pc} g_\beta(\beta)1_{\{\pc > 0\}},\mbox{ where}\\ 
g_\beta(\beta) &:= -f^\infty_{\beta}(\beta) m_{xy}^\infty(\beta) + (1-f^\infty_{\beta}(\beta))m_{yx}^\infty(\beta)   + \beta - f^\infty_{\beta}(\beta) \\
&\hspace{1cm}+ (1-\beta)f^\infty_{\beta} (\beta) (m_{xx}^\infty(\beta) + m_{xy}^\infty(\beta)) - \beta (1-f^\infty_{\beta}(\beta))  (m_{yy}^\infty(\beta) + m_{yx}^\infty(\beta))  .
    \end{aligned}
\end{align}
Similar to \cite{agarwal2021new}, we will also show that \textit{the asymptotic analysis of $\beta$ is independent of other components of $\ups$}. In particular, the result stated below shows that the analysis of the following one-dimensional ODE suffices:
\begin{align}\label{eqn_beta_ode_simple_prop}
    \dot{\beta} &= g_\beta(\beta).
\end{align}

\begin{theorem} \label{thrm_beta_ODE_prop}
\reva{Consider the interval $[0,1]$ such that $g_\beta(0) \geq 0$ and $g_\beta(1) \leq 0$.} Define ${\cal I} := \{x^*: g_\beta (x^*) = 0\}$ and say ${\cal I}= \{x_i^*:     1 \leq i \leq n\}$, for some $1 \leq n < \infty$. For each $i$, let there exist an open/closed/half-open non-empty interval around $x_i^* \in {\cal I}$, say ${\cal N}_i^*$, such that $\cup_{1\leq i\leq n} {\cal N}_i^* = [0,1]$ and ${\cal N}_i^* \cap {\cal N}_j^* = \emptyset$ for $i\neq j$. Define ${\cal N}_i^- := {\cal N}_i^*\cap[0, x_i^*)$ and 
${\cal N}_i^+ := {\cal N}_i^*\cap (x_i^*, 1]$. Let $g_\beta(x)$ be Lipschitz continuous on ${\cal N}_i^-$ and $ {\cal N}_i^+$ for each $i$:
\begin{enumerate}[label=(\roman*)]
    \item if $g_\beta(x) > 0$ for all $x \in {\cal N}_i^-$, $g_\beta(x) < 0$ for all $x \in {\cal N}_i^+$, then, $x_i^*$ is an attractor for ODE \eqref{eqn_beta_ode_simple_prop}; 
    \item   if $g_\beta(x) < 0$ for all $x \in {\cal N}_i^-$ and $g_\beta(x) > 0$ for all $x \in {\cal N}_i^+$, then, $x_i^*$ is a repeller for ODE \eqref{eqn_beta_ode_simple_prop}; 
    \item else if $g_\beta(x) > 0$ (or $g_\beta(x) < 0$) for all $x \in {\cal N}_i^-$ and $g_\beta(x) > 0$ (or $g_\beta(x) < 0$ respectively) for all $x \in {\cal N}_i^+$, then, $x_i^*$ is a saddle point for ODE \eqref{eqn_beta_ode_simple_prop}. 
\end{enumerate}
Further, ODE \eqref{eqn_ODE_prop} satisfies \ref{a3_prop}. Furthermore, the attractor and saddle sets in $\cD_I$ are respectively given by:
\begin{align*}
{\cal A} &:= \{\mathbf{h}(x^*): x^* \in {\cal I} \mbox{ is an attractor for the ODE \eqref{eqn_beta_ode_simple_prop}}\},\\
{\cal S} &:= \{\mathbf{h}(x^*): x^* \in {\cal I} \mbox{ is a repeller or saddle point for the ODE \eqref{eqn_beta_ode_simple_prop}}\} \cup \{\mathbf{0}\}, \mbox{ and}
\end{align*}entire $\cD_I$ is the combined domain of attraction for  \eqref{eqn_ODE_prop}. \eop
\end{theorem}
The above result provides the limiting behaviour of a one-dimensional ODE with possibly discontinuous right hand sides, that typically arises while studying our type of application. \revg{The condition $g_\beta(0) \geq 0$ and $g_\beta(1) \leq 0$ ensures that the interval $[0,1]$ is positive invariant for the ODE \eqref{eqn_beta_ode_simple_prop}.}  
It is important to note that in \cite[Theorem 2]{agarwal2021new}, the authors consider the function $g_\beta$ such that its zeroes could be either attractors or repellers for the ODE \eqref{eqn_beta_ode_simple_prop}. The above result is an extension of the former as here the zeroes of the function $g_\beta$ can be either attractors or repellers or saddle points for the ODE \eqref{eqn_beta_ode_simple_prop}. Such an extended result is required for the application at hand as we will see in the coming sections. 

\section{Modelling of warning dynamics using TC-BP with multiple deaths}\label{sec_warning_old}
We begin this section by demonstrating how the warning dynamics can be modelled using TC-BP with multiple deaths discussed in the previous section. Towards this, we model the copies with fake and real tags as the $x$ and $y$-type populations respectively. The time instance when a user views, tags and shares the post corresponds to the time of death of an individual in the BP. \revg{As seen in section \ref{sec_prob_desc}, in \eqref{eqn_tag_wi}-\eqref{eqn_final_shares}, the distribution of shares, types of shares, etc., depends on the type-$d$ of the user that reads the post with $d \in {\cal U}$. Thus, one can correspond each $d$-type user to a $d$-death because of the following details.  When a $d$-type user reads and shares the post, the said post becomes a read copy, resulting in a $d$-death. Further, clearly $D_x = D_y = {\cal U}$. At any given time, the proportions of the users of any type are given by $\mu_0, \mu_1, \mu_2$ and $\mu_a$, which also correspondingly represent the proportions of unread copies with np, wi, ws and a-users. Thus, one can easily infer that a type-$d$ user reads the post first among the existing unread copies, or in other words, $d$-type death occurs first with probability $\mu_d/(\mu_0+ \mu_1+\mu_2+\mu_a) = \mu_d$. 
Therefore, one can set the parameter of $d$-death as:}
\begin{align}\label{eqn_lambda}
    \lambda_{z,d} (\phi) := \mu_d  \mbox{ for all } \om, \mbox{ for each }  d \in {\cal U} \mbox{ and } z \in \{x, y\}.
\end{align}
Now, after viewing the post, if a ws-user with fake-tagged copy shares the post with fake-tag, then we say that the number of shares, $\Gamma_{xx, ws}$, corresponds to the number of $x$-type offspring produced by an $x$-type parent, when ws-death occurs. In general, the number of shares with fake and real-tag correspond to offspring of $x$ and $y$-type respectively, see \eqref{eqn_final_shares}; the number of shares (offspring) also depend upon $\Om(\cdot)$.

\revg{The underlying TC-BP with multiple death-types that models the warning dynamics \eqref{eqn_final_shares}  is exactly like the well-known irreducible BP, except for the inclusion of multiple death-types (see 
\cite{athreya2004branching}). In irreducible BPs, the extinction occurs only when both the population-types die; individual extinction of a population-type is not possible. The same is the case with our model. For example, say there are no unread copies with fake-tag, i.e., $\Cx(t) = 0$ at some time $t > 0$, while the system still has real-tagged unread copies ($\Cy(t) > 0$). Then, if at some time $t' > t$, a wi-user or ws-user reads and shares the post, then, with non-zero probability, it can tag the post as fake (see \eqref{eqn_tag_wi}, \eqref{eqn_tag_ws}). If so happens, then it will lead to new unread copies with fake-tag, i.e., $\Cx(t') > 0$. Thus, \textit{the number of fake-tagged copies can be regenerated even after they are not present on the OSN, as long as there are some unread copies of the post on the OSN}.}

Next, we provide the general framework for analyzing the warning dynamics with respect to any warning mechanism ($\omega$). \revg{Observe that when any real/fake post gets extinct, then it's effect is harmless.} Therefore, our focus shall only be on the non-extinction paths.

\subsection{Analysis of warning dynamics for general WM}
\revg{Consider a general warning mechanism defined using a continuous-function $\omega : [0,1] \to \mathbb{R}^+$ which depends only on the proportion of fake-tags $\beta$.} Further, consider any post with actuality, $u \in \{R, F\}$. Then, for each $u$, it is clear from the previous section that the analysis of the TC-BP with multiple-death types, and hence the warning dynamics, is driven by the limit mean matrix (see \eqref{eqn_mean_matrix} and \ref{a2_prop}). Thus, we first construct the limit mean matrix, $M^{\infty, u}(\beta) := [m_{ij}^{\infty, u}(\beta)]_{\{i, j \in\{x,y\}\}}$, as follows
(see \eqref{eqn_tag_wi}-\eqref{eqn_final_shares}):
\begin{align}\label{eqn_mean_matrix_red_coloring}
    M^{\infty, u}(\beta) = 
    \begin{bmatrix}
    \bigg(\mu_1\rho \alpha_x^u + \mu_2 \min\{\omega(\beta) \alpha_x^u, 1\}\bigg) m_f\eta^u  & \bigg(\mu_1(1-\alpha_x^u \rho) + \mu_2(1-\min\{\omega(\beta) \alpha_x^u, 1\}) \bigg) m_f\eta^u + \mu_a m_f \eta_a\\
    & \\ 
    \bigg(\mu_1\rho \alpha_y^u + \mu_2 \min\{\omega(\beta) \alpha_y^u, 1\}\bigg) m_f\eta^u & \bigg(\mu_1(1-\alpha_y^u \rho) + \mu_2(1-\min\{\omega(\beta) \alpha_y^u, 1\})\bigg) m_f\eta^u + \mu_a m_f \eta_a
    \end{bmatrix}.
\end{align}
Next, we will identify the attractor, repeller and saddle sets for the ODE \eqref{eqn_beta_ode_simple_prop} which will lead to the limits for the stochastic trajectory corresponding to the warning dynamics by using Theorem \ref{thrm_beta_ODE_prop} and Theorem \ref{thm1}. Towards this, observe that $d(\phi) = d^\infty(\beta) = 1$, as $\sum_{d } \lambda_{z,d} (\phi) = 1$ for any $\om$ and $z \in \{x, y\}$. This implies, $f^\infty_\beta (\beta) = \beta$ (see \eqref{eqn_f_infty}). Thus, by \eqref{eqn_beta_ODE_prop}, the function $g^u_\beta$ and the corresponding ODE \eqref{eqn_beta_ode_simple_prop} for the warning dynamics for both types of posts, $u \in \{R, F\}$, is given by:
\begin{align}\label{eqn_general_g_beta}
\dot{\beta^u} &= g_\beta^u(\beta) \mbox{ where }\\
    g^u_\beta(\beta) 
    :=  \bigg(-\beta \mu_2 - \beta \mu_1 (1-\alpha_x^u \rho) + (1-\beta) \mu_1 \rho  \alpha_y^u  + \mu_2 &\bigg(\beta \min\{\omega(\beta) \alpha_x^u, 1\} + (1-\beta)\min\{\omega(\beta) \alpha_y^u, 1\}\bigg) \bigg) m_f\eta^u - \beta \mu_a m_f \eta_a. \nonumber
\end{align}
Define $\cA^u_\beta$ as the set of attractors in $[0,1]$ and $\cR^u_\beta$ as the combined set of repellers and saddle points in $[0,1]$ for the above ODE. Then, we have the following result: 
\begin{theorem}\label{thrm_BP_to_fake}
Consider the warning dynamics as in \eqref{eqn_transition_fake_tag} and \eqref{eqn_transition_real_tag}. Let the distribution of number of friends, ${\cal F} \geq 0$ be such that $E[{\cal F}] \eta_R > 1$ and $E[{\cal F}^2] < \infty$. Then, the following statements are true for each~$u$, the actuality of post:
\begin{enumerate}[label=(\roman*)]
    \item the assumptions \ref{a2_prop} and \ref{a3_prop} hold for the ODE \eqref{eqn_ODE_prop}; hence Theorem \ref{thm1}(i) is true,
    \item the set $\cA^u_\beta \neq \emptyset$, and then $\Ups_n$ converges to $\cA^u \cup \cR^u$, as $n \to \infty$ or hovers around $\cR^u$ w.p. $1$, where $\cA^u = \{\mathbf{h}(\beta): \beta \in \cA^u_\beta\}$ and $\cR^u = \{\mathbf{0}\} \cup \{\mathbf{h}(\beta): \beta \in \cR^u_\beta\}$. 
    \item Further, any potential limit proportion corresponding to the warning dynamics, i.e., $\beta \in \cA^u_\beta \cup \cR^u_\beta$, can be bounded as below:
    \begin{align}\label{eqn_beta_bar}
        0 < \frac{\mu_1 \rho \alpha_y^u \eta^u }{q^u} =: \underline{\beta}^u &< \beta \leq \overline{\beta}^u := \frac{(\mu_2 + \mu_1 \rho \alpha_y^u) \eta^u}{q^u} \leq 1, \mbox{ where} 
    \end{align}
    the constant $q^u := \bigg( \mu_2 + \mu_1 (1 - (\alpha_x^u - \alpha_y^u) \rho)\bigg) \eta^u + \mu_a \eta_a$. \eop
\end{enumerate}
\end{theorem}
At first, observe that any warning mechanism $\omega$ only affects the likelihood of tagging the post as real or fake by a ws-user (see \eqref{eqn_tag_ws}). \revg{It does not affect the probability of a post getting viral or extinct as extinction depends on the sum current number of unread copies (i.e., sum current population in the BP)}. Now, given that our interest lies in non-extinction paths, the above Theorem gives a generalised result which holds for any warning dynamics. It is important to note that \revg{viral paths are possible only when the probability of non-extinction is non-zero; this is possible if $E[{\cal F}] \eta_R > 1$, as then the TC-BP with multiple deaths can be in throughout super-critical regime  (see Lemma \ref{lemma_dichotomy}).}

Theorem \ref{thrm_BP_to_fake}(i) implies that the warning dynamics can be approximated by the solution of the ODE \eqref{eqn_general_g_beta} over any finite-time window, where the limit mean functions are given by \eqref{eqn_mean_matrix_red_coloring}. \revg{The more important result for our context is the second part of the Theorem which} states that the stochastic trajectory $\Ups_n$ either converges to $\cA^u \cup \cR^u$ or hovers around $\cR^u$. The set \revg{$\cR^u$} contains $\mathbf{0}$ which represents the limiting behavior of the stochastic trajectory in the extinction paths. \textit{Thus, all the results henceforth will focus on deriving the limits \revg{which are not equal to $\mathbf{0}$, which in turn provide the limit proportion of fake-tags for the warning dynamics in non-extinction paths.}}

\revr{Further, Theorem \ref{thrm_BP_to_fake} provides the above limits using the zeroes $\{\beta^{\infty, u}\}$ of $g^u_\beta$ (see \eqref{eqn_general_g_beta}). Now, observe that the function $g^u_\beta$ and therefore the zeroes $\{\beta^{\infty, u}\}$ depends on $\chi$, where $\chi := \{\mu_1, \mu_2, \mu_a, b, w\}$ is the set of parameters.
For some warning mechanisms, the function $g^u_\beta$ can have multiple zeroes, and the warning dynamics can converge to one of them. Thus, one would want to ensure that the  maximum limit proportion of fake-tags for the real-post is within a given limit and optimise the minimum proportion of fake-tags for the fake-post. This aspect is considered in the optimization problem \eqref{eqn_opt_prob} given in the next section. 
In this context, we define the following Quality of Service (QoS) for any warning mechanism:
\begin{align}\label{eqn_qos}
Q := \inf\{\beta: {\beta}\in \cA^{F}_\beta \cup \cR^{F}_\beta \}.
\end{align}
Observe here that $Q = \inf({\cal L}^F)$, the objective function of \eqref{eqn_gen_opt} and is the almost sure lower bound on the limit proportion of fake-tags when the underlying post is fake. It measures the minimal extent to which a fake-post is identified by the users. From \eqref{eqn_beta_bar} of Theorem \ref{thrm_BP_to_fake}, $Q \in (\underline{\beta}^F, \overline{\beta}^F]$. We would see in the coming sections how (optimal) $Q$ varies with different warning mechanisms.}

Next, in Theorem \ref{thrm_unique_att}, we will derive some properties of $\{\beta^{\infty, u}\}$ with respect to each parameter in $\chi$\revg{, when $g_\beta^u$ has a unique zero}. This result will be instrumental in deriving important results in the coming sections. To keep it simple, we shall write $\beta^{\infty, u}(\kappa)$ and $g^u_\beta(\beta; \kappa)$ to show the dependency on the required parameter $\kappa$, an element of the tuple $\chi$. Towards this, we require the following difference term (note that $g^u_\beta(\beta^\infty(\kappa); \kappa) = 0$):
\begin{align}\label{eqn_nabla}
    \nabla^u(\kappa, \partial \kappa) := g^u_\beta(\beta^\infty(\kappa); \kappa+\partial \kappa)  -  g^u_\beta(\beta^\infty(\kappa); \kappa) = g^u_\beta(\beta^\infty(\kappa); \kappa+\partial \kappa).
\end{align}
\begin{theorem} \label{thrm_unique_att} Consider any warning mechanism, $\omega(\beta)$.
Let $\kappa$ be any parameter. Let $g^u_\beta(\beta; \kappa)$ be either a convex or concave or linear function of $\beta$ with a unique zero, $\beta^{\infty, u}(\kappa) \in (0,1)$, for each $u \in \{R, F\}$. Keeping all parameters in $\chi$ fixed, other than $\kappa$, if difference term $\nabla^u(\kappa, \partial \kappa)   > 0$ for some $\partial \kappa > 0$, then $\beta^{\infty, u}(\kappa + \partial \kappa) > \beta^{\infty, u}(\kappa)$. Else if $\nabla^u(\kappa, \delta\kappa) < 0$, then $\beta^{\infty, u}(\kappa + \partial \kappa) < \beta^{\infty, u}(\kappa)$. Else, $\beta^{\infty, u}(\kappa + \partial \kappa) = \beta^{\infty, u}(\kappa)$.
\eop
\end{theorem} 

Hereon, we will analyse the warning dynamics for some specific mechanisms. 

%
\section{Analysis of Extended Original WM (eo-WM)}
In this section, we will analyse the warning dynamics when the OSN provides the warning as in \eqref{eqn_warning}, which is originally proposed in \cite{kapsikar2020controlling}.  Recall that in \cite{kapsikar2020controlling}, the system has only ws-users who interact with the warning mechanism. Since we study the original mechanism \eqref{eqn_warning} under the influence of a variety of user behaviour, we refer to $\omega$ as \underline{extended original warning mechanism} (eo-WM) in our context.

Consider any post with actuality $u \in \{R, F\}$. Recall that we have $w \leq \overline{w} := \frac{1}{\alpha_x^F} - \gamma$, thus leading to $\alpha_j^u \omega(\beta) < 1$ for each $j \in \{x, y\}$ and for any $\beta \in [0,1]$. We begin the analysis by analyzing the ODE \eqref{eqn_general_g_beta} for the eo-WM. The $g_\beta^u$ defined in \eqref{eqn_general_g_beta} for the eo-WM, henceforth denoted as $g^{o, u}_\beta$, is as given below:
\begin{align}\label{eqn_beta_ODE_etac1}
\begin{aligned}
    g^{o, u}_\beta(\beta) &= -\beta \mu_2 m_f\eta^u - \beta \mu_1 (1-\alpha_x^u \rho) m_f\eta^u + (1-\beta) \mu_1 \rho  \alpha_y^u m_f\eta^u + \mu_2 \omega(\beta) \bigg(\beta \alpha_x^u + (1-\beta)\alpha_y^u\bigg) m_f\eta^u - \beta \mu_a m_f \eta_a.
\end{aligned}
\end{align}
Let $\cA^{o, u}_\beta \subset [0,1]$ be the corresponding attractor set and $\cR^{o, u}_\beta \subset [0,1]$ be the union of the corresponding repeller and saddle sets, i.e., with respect to ODE $\dot{\beta}^u = g^{o, u}_\beta(\beta)$. We study these sets in the following.
\begin{corollary}\label{corollary_ex_wm}
  There exists a unique zero, $\beta^{o, \infty, u}$, of $g^{o, u}_\beta$ in $[0,1]$. Further, $\beta^{o, \infty, u} \in (0,1)$, $\cA^{o, u}_\beta = \left\{\beta^{o, \infty, u}\right\}$ and $\cR^{o,u}_\beta = \emptyset$.  \eop
\end{corollary}
\revg{Thus, there is a unique attractor, $\beta^{o, \infty, u}$, of ODE \eqref{eqn_general_g_beta}. By Theorem \ref{thrm_BP_to_fake}, the stochastic trajectory $\Ups_n$ under eo-WM either converges to $\mathbf{h}(\beta^{o, \infty, u})$ or $\mathbf{0}$, or hovers around $\mathbf{0}$ almost surely. We re-iterate that our focus is on the non-extinction paths, and thus, the relevant proportion of fake-tags is unique and equals $\beta^{o, \infty, u}$. Further, by Theorem \ref{thrm_BP_to_fake}, for the given choice of $w, b$ and given $\mu_1, \mu_2, \mu_a \in \chi$, $\beta^{o, \infty, u} \in (\underline{\beta}^u, \overline{\beta}^u]$.}

We now consider the following robust optimization problem for the OSN discussed before:
\begin{align}\label{eqn_opt_prob}
\begin{aligned}
    \sup_{w \in \left[0,  \overline{w}\right], b \in [0, \infty)} \beta^{o, \infty, F}(w, b)
    \mbox{ subject to } \beta^{o, \infty, R}(w, b) \leq \delta, \mbox{ for some } \delta \in (0,1).
\end{aligned}
\end{align}
\revr{By uniqueness of the attractors in the non-extinction paths, the above constrained optimization problem optimizes the QoS defined in \eqref{eqn_qos}, $Q = \beta^{o, \infty, F}$ under eo-WM by choosing $w, b$, while ensuring that the unique zero for the real-post, $\beta^{o, \infty, R} \leq \delta$.}
The problem in \eqref{eqn_opt_prob} is exactly the same as in \cite{kapsikar2020controlling}, but for the inclusion of different user behaviour in our model. Thus, we need to extend the solution of \cite{kapsikar2020controlling} to the case that includes wi, ws, a, and np-users. Observe that $\delta$ is a design parameter for the OSN.

Before we solve the above problem, we observe the following qualitative behaviour which is true by the virtue of Theorem \ref{thrm_unique_att} -- this behavior is important for further analysis:
\begin{corollary}\label{cor_limits_warning}
For each $u \in \{R, F\}$, the limit $\beta^{o, \infty, u}(w, b)$ strictly increases with $w$ and strictly decreases with $b$. 
\eop
\end{corollary}
The above Corollary intuitively indicates to choose the largest $w$, i.e., $\overline{w}$ and the smallest $b$, i.e., $0$. However, since the optimal $w, b$ needs to satisfy the constraint for the real-post as in \eqref{eqn_opt_prob}, therefore, formal analysis is required.

\begin{theorem}{\bf [Optimal eo-warning design]} \label{thrm_opt} The following statements hold for the optimization problem \eqref{eqn_opt_prob}:
\begin{enumerate}[label=(\roman*)]
    \item If $\beta^{o, \infty, R}(\overline{w}, 0) > \delta$, then 
    the optimizer $(w^*, b^*)$ of \eqref{eqn_opt_prob} is as below and satisfies $\beta^{o, \infty, R}(w^*, b^*) = \delta$:
    \begin{align}\label{eqn_optimal_parameter_old}
    \begin{aligned}
        w^* = \overline{w} \mbox{ and } b^* 
        &= \left(\frac{\delta}{1-\delta}\right)\left(\frac{w^*\eta^R \mu_2(\delta \alpha_x^R + (1-\delta)\alpha_y^R)}{\delta ((\mu_1+\mu_2)\eta^R + \mu_a  \eta_a) - \eta^R (\mu_1 \rho + \mu_2 \gamma) (\delta \alpha_x^R + (1-\delta) \alpha_y^R)} - 1 \right).
    \end{aligned}
    \end{align}
    \item Else, if $\beta^{o, \infty, R}(\overline{w}, 0) \leq \delta$, then $(w^*, b^*) = (\overline{w}, 0)$ and satisfies $\beta^{o, \infty, R}(w^*, b^*) \leq \delta$.  \eop
    \end{enumerate}
\end{theorem}
Thus, as anticipated, $w^* = \overline{w}$. Interestingly, contrary to the expectation, $b^*$ is not always $0$. If $\beta^{o, \infty, R}(\overline{w}, 0) > \delta$, then the optimal choice for $b$ is given by $b^* > 0$. Such $b^*$ is achieved by solving for $\beta^{o, \infty, R}(w^*, b) = \delta$, i.e., relaxing the constraint for the real-post to the maximum $\delta$-level in a bid to achieve the maximum $\beta^{o, \infty, F}$ for fake-post at optimality. In view of Corollary \ref{cor_limits_warning}, it is then easy to see that, $\beta^{o, \infty, F}(w^*, b^*) < \beta^{o, \infty, F}(w^*, 0)$, when $b^* > 0$. For simpler notations, henceforth we will refer to $\beta^{o, \infty, F}(w^*, b^*)$ as $\beta^o$ and $\beta^{o, \infty, R}(w^*, b^*)$ as $\beta^{o, R}$. 

In \cite{kapsikar2020controlling}, the optimization problem \eqref{eqn_opt_prob} is solved partially. Firstly, only the case with the hypothesis of Theorem \ref{thrm_opt}(i) is analyzed. It is shown that the optimal value is achieved for the value of $b$ which satisfies $\beta^{o, \infty, R} = \delta$. However, the optimal choice of $w$ is not derived; rather a projected gradient descent algorithm is suggested to evaluate $w^*$. Furthermore, \cite{kapsikar2020controlling} considers $w \in [0,1]$, while one can allow $w$ to be as large as $\overline{w}$, which can be larger than $1$. As we have proved that $w^* = \overline{w}$, therefore, \textit{our optimal eo-WM should perform better than the optimal WM designed in \cite{kapsikar2020controlling}}. We numerically show this aspect in the next sub-section. 

\subsection{QoS under eo-WM}
\revg{It is clear from Corollary \ref{corollary_ex_wm} and Theorem \ref{thrm_opt} that the QoS under optimal eo-WM, say $Q^o$ equals $\beta^o$. Now, fix any configuration, }
$$
{\cal C} := \bigg\{\{\alpha_i^u\}_{\{i \in \{x, y\}, u \in \{F, R\}\}}, \eta_a, \{\eta^u\}_{\{u \in \{F, R\}\}}, \rho, \gamma, m_f, w, \mu_1, \mu_2 \bigg\},
$$\revg{and let $\mu_a$ vary. Then, we want to investigate how $Q^o$ changes with $\mu_a$. Towards this, define:
\begin{align}
\betana := \beta^o(\mu_a = 0) = Q^o(\mu_a = 0),
\end{align}as the proportion of fake-tags for the fake-post at optimality when there is no adversary. Recall that a-users deliberately tag any post as real. Therefore, 
one can anticipate that the OSN achieves the maximum QoS when there is no adversary, i.e., $\beta^o(\mu_a) = Q^o(\mu_a) < \betana$, when $\mu_a > 0$. We precisely prove this in the next result for an appropriate range of $\delta$.}
\begin{corollary}\label{cor_beta_o_na}
    For given configuration ${\cal C}$, there exists a $\overline{\delta} > 0$ such that $\beta^o(\mu_a) < \betana$ for all $\delta \leq \overline{\delta}$, for each $\mu_a \in (0, 1-\mu_1-\mu_2]$. \eop
\end{corollary}
\revg{Thus, the above corollary confirms our anticipation that the performance degrades with introduction of the a-users in the system, however for a smaller range of $\delta$; observe that the OSN is interested in keeping $\delta$ as small as possible, therefore, such choices of $\delta$ are indeed meaningful. \textit{Henceforth, we consider such $\delta$, i.e., $\delta \leq \overline{\delta}$}. }
 In the next subsection, we will validate this result numerically and reinforce the requirement to design better WMs in the presence of adversaries.

\subsection{Numerical analysis for eo-WM}
\revg{At first, we would like to compare eo-WM with the mechanism in \cite{kapsikar2020controlling} with just a-users added -- in the first example, any user on the OSN can either be a ws-user or an a-user ($\mu_2 + \mu_a = 1$). Thus, there is no wi-user and everyone participates.} Further, let the parameters be as in \cite{kapsikar2020controlling}:
\begin{align}\label{eqn_param_plot1}
\begin{aligned}
    m_f &= 28, \eta^F = 0.08, \eta^R = 0.05, \gamma = 0.1, \eta_a = 0.55, \delta = 0.02,
    \alpha_x^F = 0.85, \alpha_y^F = 0.6375, \alpha_x^R = 0.3 \mbox{ and } \alpha_y^R = 0.09.
\end{aligned}
\end{align}
\begin{wrapfigure}{r}{6.6cm}
\centering
\includegraphics[trim = {1.6cm 6.5cm 0cm 6.5cm}, clip, scale = 0.3]{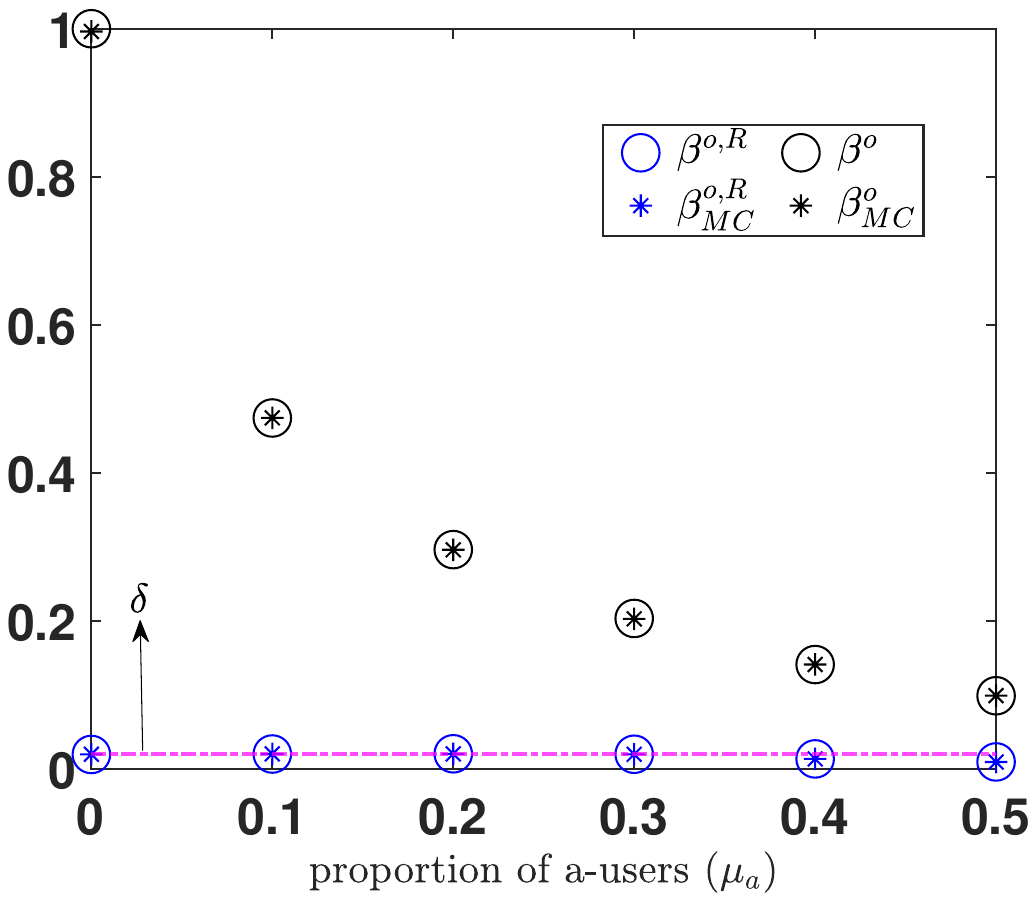}
  \caption{Limits of warning dynamics under eo-WM}
  \label{fig:exWM}
\end{wrapfigure}
For such parameters, we perform Monte-Carlo  (MC) simulation, and also evaluate the zeroes of $g_\beta^{o, u}$ for each $u \in \{R, F\}$. In Figure \ref{fig:exWM}, we plot the outputs of MC simulations and the theoretical limits against $\mu_a$, which can be seen to be close to each other. The constraint for the real-post is satisfied. In fact, the proportion of tags (for fake post) decreases with $\mu_a$, which is intuitive as a-users deliberately real-tag the posts.

Under the optimal eo-WM, $99.981\%$ of users can identify the fake-post as fake; this optimal value is higher than the reported $90\%$ in \cite{kapsikar2020controlling}, as we use $w^* = 1.076$, while algorithm  in \cite{kapsikar2020controlling}  uses $w^* = 1$. 
Now, it is interesting to note that with just $1\%$ and $2\%$ of a-users on the OSN, the performance of the eo-WM decreases to $89.798\%$ and $81.74\%$ respectively (\revr{in fact, there is degradation with respect to the new QoS defined in \eqref{eqn_iqos} which focuses only on non-adversarial users; $99.981\%$ decreases to $95.8\%$ and $92.53\%$ respectively with $1\%$ and $2\%$ of a-users}). This depicts that the original WM is not sufficient to control the fake-post propagation in the presence of adversaries.

\revr{Next, we consider a second example with parameters almost as in \eqref{eqn_param_plot1}, but with 
proportion of ws-users ($\mu_2$) fixed and with $\mu_a$, the proportion of a-users varying. We set $\mu_2 = 0.5$, 
 $\mu_1 = 0$ and let the fraction of non-participants equal $0.5 - \mu_a$.  
For ease of reference, the users of this example are referred to as `\underline{smart users}', as here $\alpha_x^F - \alpha_x^R = 0.55$ and $\alpha_y^F - \alpha_y^R = 0.5475$ indicating that the users are capable of distinguishing the fake posts from real posts to a good extent, even without external aid and irrespective of sender tag.

We compare smart users with users in another example scenario where $\alpha_x^F - \alpha_x^R = 0.18$ and $\alpha_y^F - \alpha_y^R = 0.135$. As the differences between the distinguishing parameters are small, these users are referred to as `\underline{naive users}'. For this example, the remaining parameters are fixed as below (for diversity, we also consider more attractive posts):
\begin{align}\label{eqn_param_plot3}
\begin{aligned}
    \rho &= 0.9, m_f = 30, \eta^F = 0.52, \eta^R = 0.4, \gamma = 0.1, \eta_a = 0.55, \delta = 0.05,\\
    \alpha_x^F &= 0.3, \alpha_y^F = 0.225, \alpha_x^R = 0.12, \alpha_y^R = 0.09, \mu_1 = 0.15 \mbox{ and } \mu_2 = 0.5.
\end{aligned}
\end{align}Typically, the users may be naive -- may not possess sufficient intrinsic ability to distinguish between the posts to the level that smart users can. Interestingly as seen below, the proposed mechanism is effective to guide even naive users.

\begin{figure}[http]
\begin{minipage}{\textwidth}
\centering
\begin{minipage}{.5\textwidth}
  \centering
  \includegraphics[trim = {1cm 6cm 0cm 6cm}, clip, scale = 0.3]{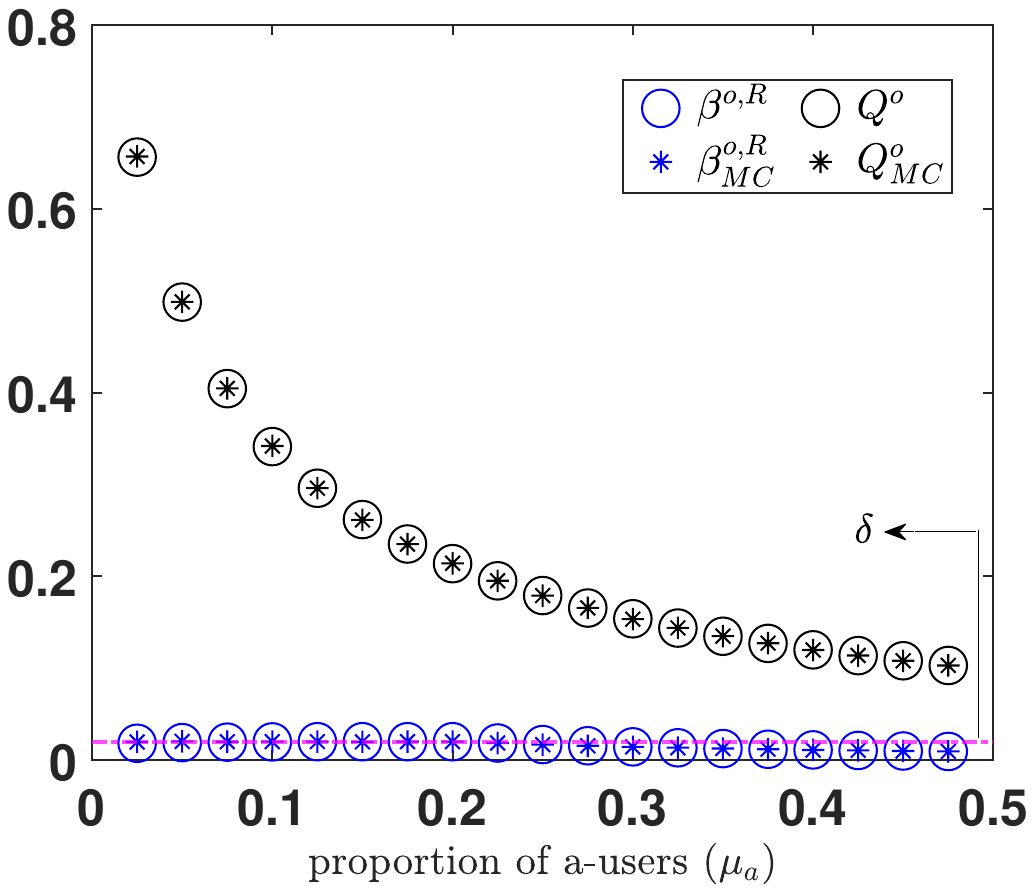}
\end{minipage}%
\begin{minipage}{.5\textwidth}
  \centering
  \includegraphics[trim = {1cm 6cm 0cm 6cm}, clip, scale = 0.3]{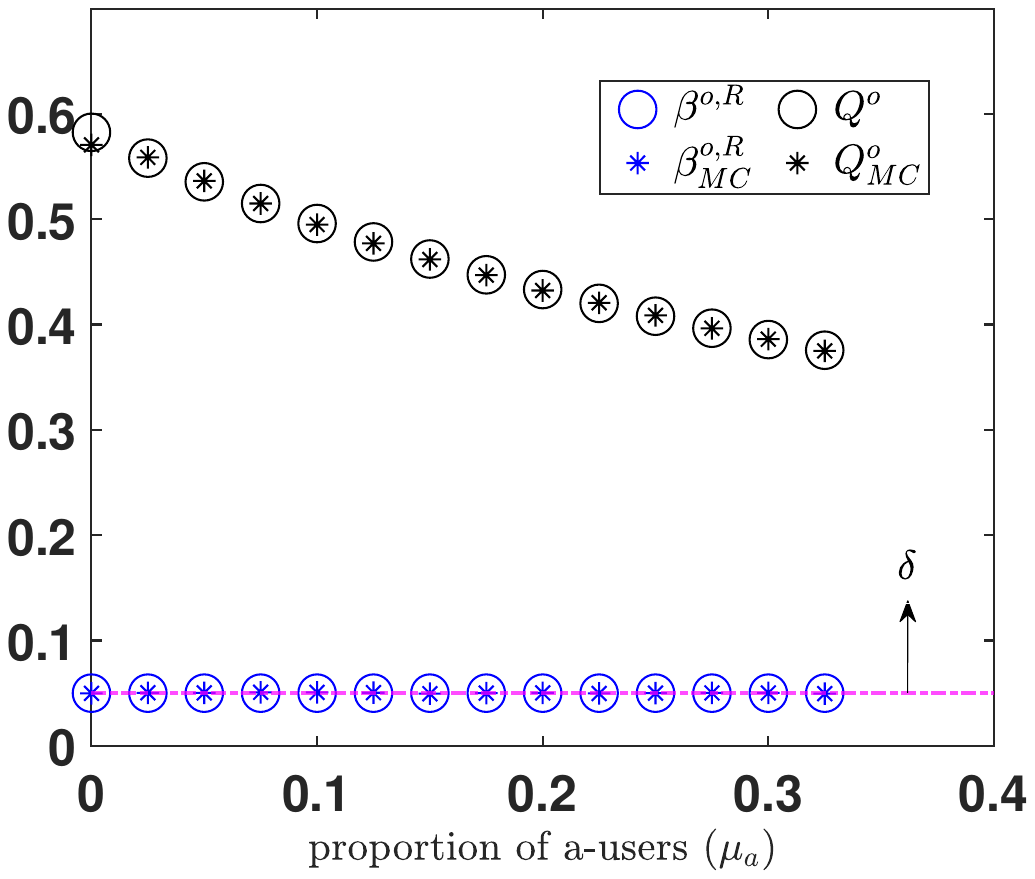}
\end{minipage}
\caption{Limits of warning dynamics under eo-WM with smart (left) and naive (right) users respectively}
\label{fig:eoWM_delta}
\end{minipage}
\begin{minipage}{\textwidth}
\centering
\begin{minipage}{.5\textwidth}
  \centering
  \includegraphics[trim = {1cm 6cm 0cm 6cm}, clip, scale = 0.3]{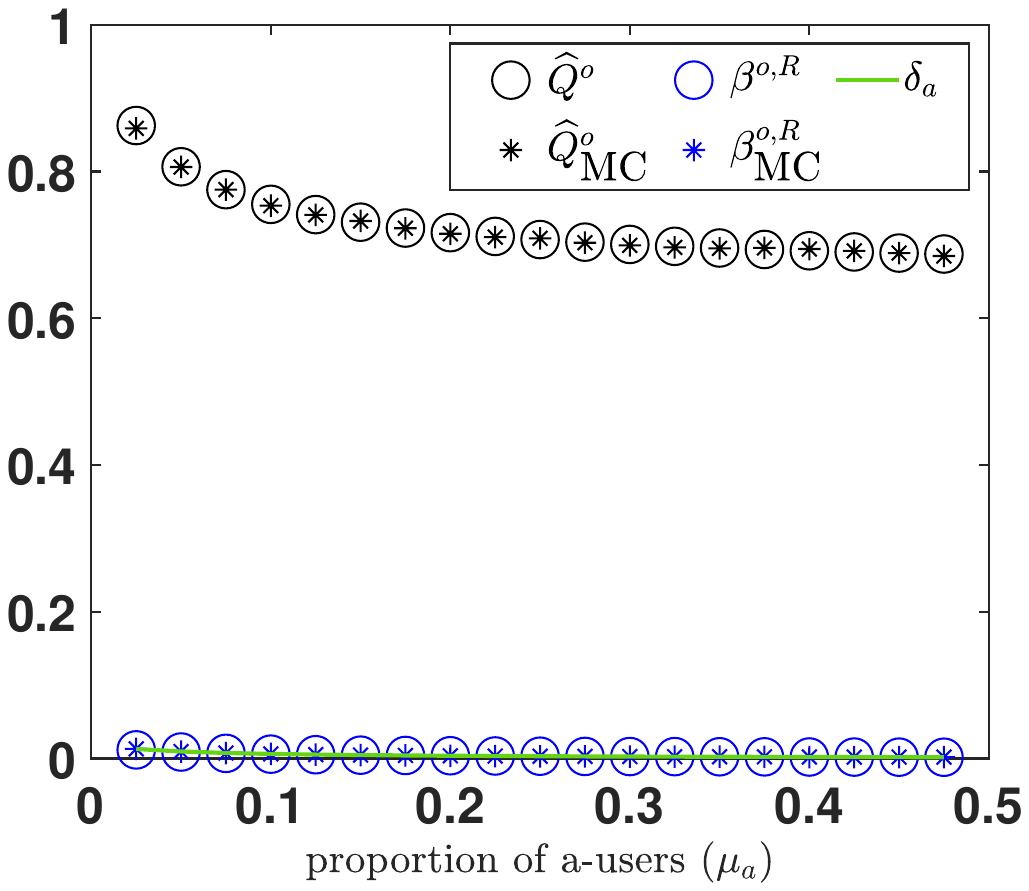}
\end{minipage}%
\begin{minipage}{.5\textwidth}
  \centering
  \includegraphics[trim = {1cm 6cm 0cm 6cm}, clip, scale = 0.3]{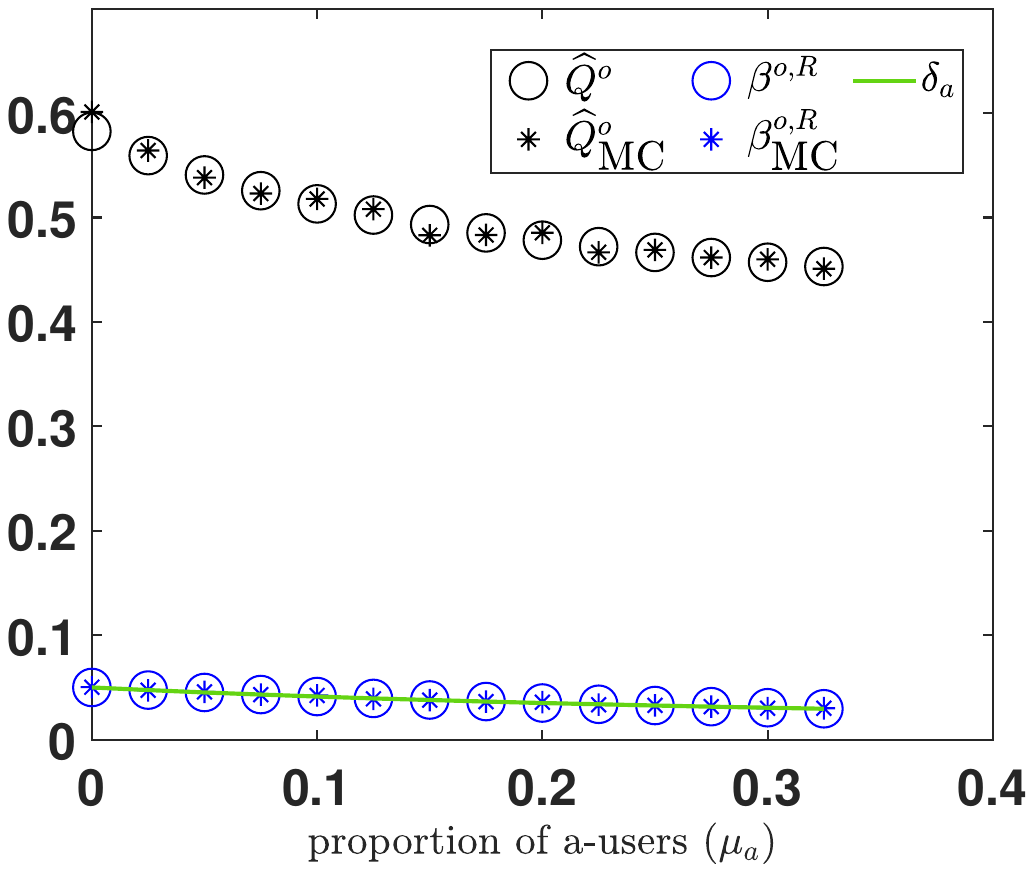}
\end{minipage}
\caption{i-QoS under eo-WM with smart (left) and naive (right) users respectively}
\label{fig:QOS_eo}
\end{minipage}
\end{figure}

In Figure \ref{fig:eoWM_delta}, we illustrate the QoS ($Q$ defined in \ref{eqn_qos}) and the proportion of fake-tags for the real-post for examples with smart and naive users in left and right sub-figures respectively. Many of the observations are similar to that in the first example: QoS decreases with an increase in $\mu_a$, and the proportion of fake-tags for real-post is at most $\delta$. The QoS in the left sub-figure with smart users is also less than that for first example provided in Figure \ref{fig:exWM}, which also considers smart users -- however, for the example in Figure \ref{fig:eoWM_delta}(left), the proportion of ws-users ($\mu_2$) is lesser than that in Figure \ref{fig:exWM} and the number of np-users is non-zero. Furthermore, as one may anticipate, the QoS with naive users is even smaller.}

\revr{\subsection[fdfdf]{Improved QoS -- QoS among non-adversaries\footnote{We would like to thank the  reviewers for their feedback which motivated us to design a new performance metric (named improved QoS (i-QoS)) that helps in explaining the numerical results in a more meaningful way.}}
It is important to note that the OSN can  control/guide the fake-tags only from non-adversarial users. In fact, the aim is also confined to the correct identification of the actuality of the posts by such users. Hence, it is more appropriate to consider a metric/QoS focused on the proportion of fake-tags only from ws and wi-users. We precisely aim to capture the same in this sub-section, and define the appropriate optimization problem. Towards this, let $X_1^u, X_2^u, X_a^u$ be the respective proportion of wi, ws and a-users at limit who fake-tag the $u$-post; observe $X_a^u = 0$ and recall, np-users do not participate. Similarly, define $Y_1^u, Y_2^u, Y_a^u$ as the corresponding proportion of users who real-tag; observe $Y_a^u = 1$. The limit approaches when the number of users that read the post, $t \uparrow \infty$, and consider a large enough $t$. Then, the number of fake-tags by ws-users after $t$-th user reads the post can be approximated by $tX_2^u m_f \eta^u$ (one can anticipate this by law of large numbers and because of \ref{a2_prop}). The other numbers can be approximated in a similar way and as a result, one can re-write the overall proportion of fake-tags as:
\begin{align*}
    \beta^u \approx \frac{(X_1^u + X_2^u)m_f \eta^u}{(X_1^u + X_2^u + Y_1^u + Y_2^u)m_f \eta^u + Y_a^u m_f \eta_a}.
\end{align*}
In a similar manner, the proportion of fake-tags from non a-users represented by $\beta_a^u$, the quantity of actual interest, can be approximated as below:
$$
\beta_a^u \approx \frac{(X_1^u + X_2^u)m_f \eta^u}{(X_1^u + X_2^u + Y_1^u + Y_2^u)m_f \eta^u} = \frac{X_1^u + X_2^u}{X_1^u + X_2^u + Y_1^u + Y_2^u}.
$$
Thus, one can relate the two QoS metrics as follows:
\begin{align}\label{eqn_beta_u_only_wiws}
    \beta_a^u = \left(\frac{(\mu_1 + \mu_2) \eta_u + \mu_a \eta_a}{(\mu_1 + \mu_2) \eta^u}\right)\beta^u.
\end{align}
The above discussion motivates us to define an `improved quality of service (i-QoS)' with respect to any warning-mechanism:
\begin{align}\label{eqn_iqos}
    \widehat{Q} := \inf \left \{ \left(\frac{(\mu_1 + \mu_2) \eta^u + \mu_a \eta_a}{(\mu_1 + \mu_2) \eta^u}\right)\beta : \beta \in \cA_\beta^F \cup \cR_\beta^F\right \} =  \left(\frac{(\mu_1 + \mu_2) \eta^u + \mu_a \eta_a}{(\mu_1 + \mu_2) \eta^u}\right)Q.
\end{align}
One can interpret $\widehat{Q}$ as the almost sure lower bound on the limit proportion of fake-tags for fake-post from non a-users. Henceforth, we also consider comparison of various warning mechanisms using this more relevant metric, i-QoS. Further, we illustrate a lot more improvement when optimization problem \eqref{eqn_opt_prob} is instead designed using i-QoS.
Observe that i-QoS is simply a constant multiple of QoS, and hence by Corollary \ref{corollary_ex_wm} and by \eqref{eqn_iqos}, the i-QoS for eo-WM (represented by $\widehat{Q}^o$) is unique. Thus, the original optimization problem \eqref{eqn_opt_prob} changes to the following, for some $\delta \in (0,1)$:
\begin{align}\label{eqn_new_opt}
    \sup_{w \in \left[0,  \overline{w}\right], b \in [0, \infty)} \widehat{Q}^o(w, b)
    \mbox{ subject to } {\beta}^{o, \infty, R}(w, b) \leq \delta_a := \frac{\delta ((\mu_1 + \mu_2)\eta^R) }{((\mu_1 + \mu_2)\eta^R + \mu_a\eta_a)}.
\end{align}
Observe that the above optimization problem has exactly the same structure as in \eqref{eqn_opt_prob}, except that $\delta$ is replaced by $\delta_a$; hence, $w^*, b^*$ can be derived by Theorem \ref{thrm_opt} directly. The optimal value of the above problem represents the fraction of non a-users (wi and ws-users) who correctly identify the fake-post as fake. 
When $\mu_a > 0$ and is sufficiently large, then QoS is sufficiently small (lesser than $1-\mu_a$), as it includes the effects of a-users real-tagging. But this does not imply that the WM failed; in fact on the contrary, at the extreme end, WM is completely successful in eliminating the effect of adversaries if optimal $\widehat{Q}^o = 1$  (indicating that all the non a-users correctly identify the fake-post).

In Figure \ref{fig:QOS_eo}, we continue with the two examples of Figure \ref{fig:eoWM_delta}, where we plot i-QoS and its MC estimates, and the corresponding quantities for the real-post; the left sub-figure has smart users and right sub-figure has naive users. It is clear that the proportion of fake-tags for the real-post ($\beta^{o, R}$, see blue curves) are within $\delta_a$-threshold for both cases. More interestingly, the results of the said figure for the fake-post indicate that the results of Figure \ref{fig:eoWM_delta} are mis-leading; the latter figure shows extremely high level of degradation in QoS with $\mu_a$, while the same is not the case in the former; this is obviously because  the latter also counts the (intentional) real-tags from a-users. For example, when $\mu_a = 0.3$, the QoS is $15.38\%$ in Figure \ref{fig:eoWM_delta}(left), while the actual fraction of fake tags among the smart non a-users is around $70.06\%$. Thus, the degradation with $\mu_a$ may not be as large as depicted in Figure \ref{fig:eoWM_delta}, nonetheless there is sufficient degradation as $\mu_a$ increases (for example, from $99.981\%$ at $\mu_a = 0$ to $70.06\%$ for $\mu_a = 0.3$). 

The above illustrations motivate us to design better warning mechanisms, which achieve higher performance. In fact, the underlying theme of the entire paper is to optimize/increase the proportion of fake-tags for the fake-post, while still ensuring that the constraint in \eqref{eqn_new_opt} for the real-post is satisfied. In this section, we optimized the performance of the eo-WM for the fake-post, and achieved exactly $\delta$-threshold for the real-post. In the coming sections, we will attempt at designing WMs which increase the performance, without compromising over the real-post. As mentioned before, this goal is achieved by designing appropriate WMs such that the resultant $g_\beta$ of \eqref{eqn_general_g_beta} has zeroes with desirable properties, which in turn dictate the limiting behaviour of WM as confirmed by Theorem \ref{thrm_BP_to_fake}. }  To this end, the first idea is to eliminate the effect of adversaries, which we consider next.

\section{Eliminating Adversarial Effect WM (ea-WM)}
The OSN may not know the exact set of adversarial users, but it knows the proportion of adversarial users ($\mu_a$). We aim to use this knowledge to design a new improved WM which performs better  even when $\mu_a$ is large. \textit{The idea is to construct a WM specific to any given $\mu_a > 0$, namely $\omega^a(\beta)$, such that $g^F_\beta$ under the new WM exactly matches that corresponding to $g^{o, F}_\beta$ with $\mu_a = 0$, at optimality (see \eqref{eqn_beta_ODE_etac1}).} In other words, using the knowledge of $\mu_a$, we are creating a hypothetical situation with no adversaries, and hence we name $\omega^a$ as \underline{eliminating adversarial effect WM} (ea-WM). If that is possible, then one can anticipate that the performance will improve for the fake-post under ea-WM; we will identify such conditions below. Further, one still needs to ensure that the performance of real-post is not degraded beyond $\delta$ as in \eqref{eqn_opt_prob} \revr{(beyond $\delta_a$ as in \eqref{eqn_new_opt} when i-QoS is considered); this is ensured by the WM proposed in this section (and for coming WMs as well).} Towards this, we simply define $\omega^a$ as:
\begin{align}\label{eqn_warning_ea}
    {\omega}^a(\beta) = \omega(\beta) + \frac{\beta \mu_a m_f \eta_a}{ \mu_2 m_f\eta^F \bigg(\beta \alpha_x^F + (1-\beta)\alpha_y^F\bigg) }.
\end{align}
Consider $w, b$ and $\beta$ such that $\min\{\alpha_x^u \omega^a(\beta), 1\} = \alpha_x^u \omega^a(\beta)$. Then $g^F_\beta$ under ea-WM, henceforth denoted as $g^{a, F}_\beta$, matches with $g^{o, F}_\beta(\beta; \mu_a = 0)$, because (see \eqref{eqn_general_g_beta}):
\begin{align}\label{eqn_g_F_eaWM}
\begin{aligned}
    g^{a, F}_\beta(\beta; \mu_a > 0) &= -\beta \mu_2 m_f\eta^F - \beta \mu_1 (1-\alpha_x^F \rho) m_f\eta^F + (1-\beta) \mu_1 \rho  \alpha_y^F m_f\eta^u + \mu_2 \omega^a(\beta) \bigg(\beta \alpha_x^F + (1-\beta)\alpha_y^F\bigg) m_f\eta^F\\
    &= g^F_\beta(\beta; \mu_a = 0).
\end{aligned}
\end{align}
Thus, if $\min\{\alpha_x^u \omega^a(\beta), 1\} = \alpha_x^u \omega^a(\beta)$ is satisfied in a neighborhood of $\betana$, then one can design the required ea-WM, if further the performance of real-post is within $\delta$-threshold \revr{(or $\delta_a$-threshold)}. In view of Theorem \ref{thrm_opt}, we set $w, b$ as follows for the new ea-WM \revr{(similarly, with $\delta_a$)}:
\[
w = \overline{w} \mbox{ and } b = 
\begin{cases}
    b^*|_{\mu_a = 0} = \left(\frac{\delta}{1-\delta}\right)\left(\frac{\overline{w} \mu_2(\delta \alpha_x^R + (1-\delta)\alpha_y^R)}{\delta (\mu_1+\mu_2) - (\mu_1 \rho + \mu_2 \gamma) (\delta \alpha_x^R  + (1-\delta) \alpha_y^R)} - 1 \right), &\mbox{if } \beta^{a, \infty, R}(\overline{w}, 0) > \delta,\\
    0, & \mbox{otherwise}.
\end{cases}
\]
Now, similar to eo-WM, for each $u \in \{R, F\}$, we will first identify the set of attractors ($\cA^{a, u}_\beta$) and the combined set of repellers and saddle points ($\cR^{a, u}_\beta$) for the \revg{ODE \eqref{eqn_general_g_beta} under ea-WM, i.e., $\dot{\beta}^u = g_\beta^{a, u}(\beta)$.}
\begin{theorem}\label{corollary_ea_wm}
    Define 
    \begin{align}\label{eqn_l_threshold_mua_eaWM}
        \Delta_a := \mu_2 \eta^F \left(\frac{1}{\alpha_x^F} - \omega(\betana)\right) \left(\frac{\betana\alpha_x^F + (1-\betana)\alpha_y^F}{\betana \eta_a}\right). 
        \end{align}
    Then, the following statements are true for the fake-post:
    \begin{enumerate}[label=(\roman*)]
        \item If $0 < \mu_a \leq \min\{1-\mu_1-\mu_2, \Delta_a\}$, then ${\beta}^{a} \geq \betana$ for all ${\beta}^{a} \in \cA^{a, F}_\beta \cup \cR^{a, F}_\beta$.
    
        \item Else, i.e., if $\Delta_a < \mu_a < 1-\mu_1-\mu_2$, then $ {\beta}^{a} \in (\beta^o, \betana)$ for all ${\beta}^{a} \in \cA^{a, F}_\beta \cup \cR^{a, F}_\beta$.
    \end{enumerate}
    For the real-post, ${\beta}^{a, R} < \delta$ for all ${\beta}^{a, R} \in \cA^{a, R}_\beta \cup \cR^{a, R}_\beta$. \eop
\end{theorem}
In view of the above and Theorem \ref{thrm_BP_to_fake}, we get that the stochastic iterates $\Ups_n$ under ea-WM for the $u$-post either converge to $\{\mathbf{h}(\beta) : \beta \in \cA^{a, u}_\beta \cup \cR^{a, u}_\beta \} \cup \{\mathbf{0}\}$, or hover around $\{\mathbf{h}(\beta) : \beta \in \cR^{a, u}_\beta \} \cup \{\mathbf{0}\}$. \revr{Unlike eo-WM, above Theorem does not guarantee unique limit for the warning dynamics under ea-WM in the non-extinction paths, but Theorem \ref{thrm_BP_to_fake}(ii) ensures that there exists at least one attractor of the ODE \eqref{eqn_general_g_beta}, as} $\cA^{a, u}_\beta \neq \emptyset$.

\revr{Now, note that ea-WM provides higher warning in comparison to the eo-WM, even for the real-post. Even with such a WM, it is proved above that the proportion of the real-post is maintained\footnote{Some equilibrium points can be saddle points and according to Theorem \ref{thrm_BP_to_fake}, the warning dynamics can hover around such points. But then the warning dynamics moves arbitrarily close to such points and we assume the equilibrium points to be the representative of the limiting behaviour. This leads to a small level of inaccuracy in the sense that the warning dynamics can go above or below the point, in case of hovering around.} within $\delta$-threshold. Further, due to higher warning, we expect a higher QoS under ea-WM; next we discuss the same.
Let the QoS \eqref{eqn_qos} under ea-WM be represented by $Q^a$. In view of Theorem \ref{corollary_ea_wm}, we claim that $Q^a > Q^o$ for the following reasons:

(i) when $\mu_a$ is small, i.e., when  $\mu_a \leq \Delta_a$, we have $Q^a \geq \betana > Q^o$ (by Theorem \ref{corollary_ea_wm}(i) and Corollary \ref{cor_beta_o_na}). Thus, ea-WM with adversaries achieves higher QoS than the original eo-WM with no adversary. Then, one can say that the former eliminated the effect of adversaries completely.

(ii) when $\mu_a$ is larger, i.e., when $\mu_a > \Delta_a$,  ea-WM still improves over eo-WM as $Q^a > Q^o$ by Theorem \ref{corollary_ea_wm}(ii). However, in this case, the QoS under ea-WM is lesser than the QoS under eo-WM with no adversary as $Q^a < \betana $. Thus, in this case, the effect of adversaries is not completely eliminated by ea-WM.

Similar design and  observations follow when one attempts to design ea-WM with i-QoS, i.e., by replacing $\delta$ with $\delta_a$. Recall again that with i-QoS, we consider a more relevant problem that focuses only on the responses from non a-users.}

\subsection{Numerical analysis for ea-WM}\label{num_exp_eaWM}
In this sub-section, we will numerically quantify the improvement achieved by ea-WM, in comparison to eo-WM; we consider only i-QoS based problems and results.
\revr{In Figure \ref{fig:eaWM}, we continue with the two examples considered in Figure \ref{fig:exWM} (i.e., with smart and naive users) for ea-WM. We plot the i-QoS with respect to ea-WM (denoted as $\widehat{Q}^a$) evaluated via the exact zeroes of $g^{a, u}_\beta$ and the corresponding MC estimates for the ea-WM. We once again observe a close match between the theoretical expressions and the corresponding MC estimates. 
\begin{figure}[http]
\centering
\begin{minipage}{.5\textwidth}
  \centering
  \includegraphics[trim = {1cm 6cm 0cm 6cm}, clip, scale = 0.3]{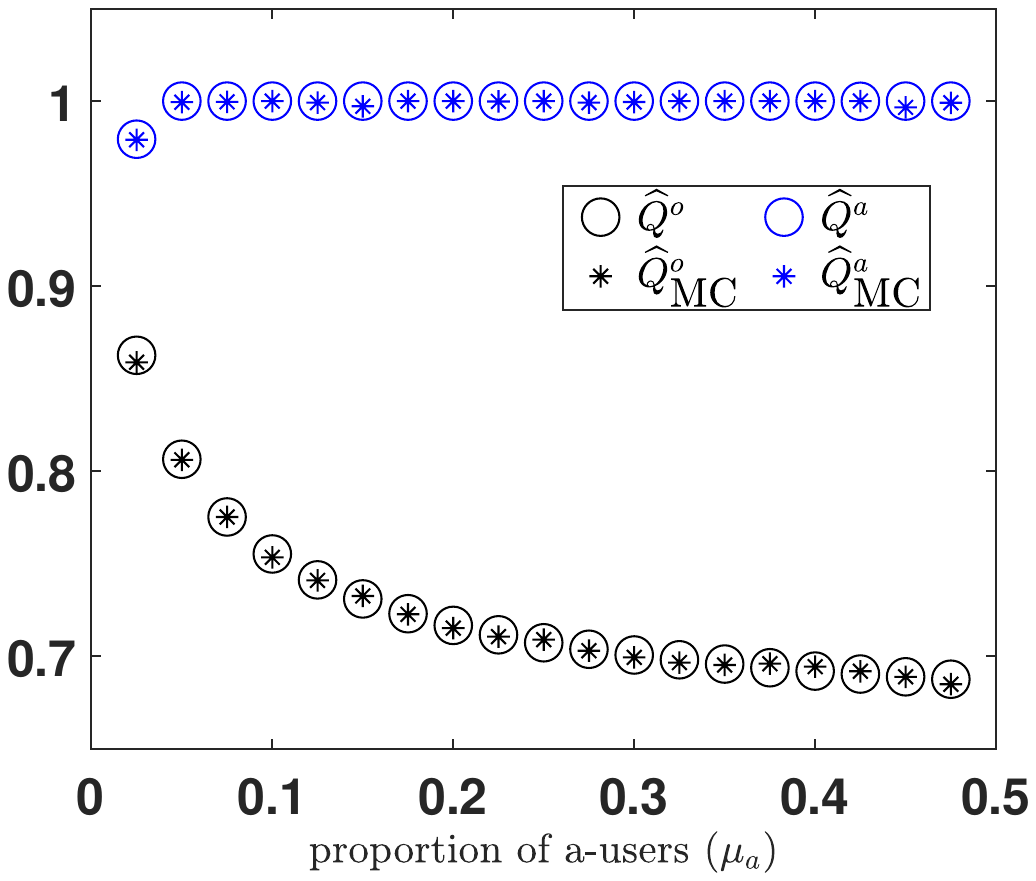}
\end{minipage}%
\begin{minipage}{.5\textwidth}
  \centering
  \includegraphics[trim = {1cm 6cm 0cm 6cm}, clip, scale = 0.3]{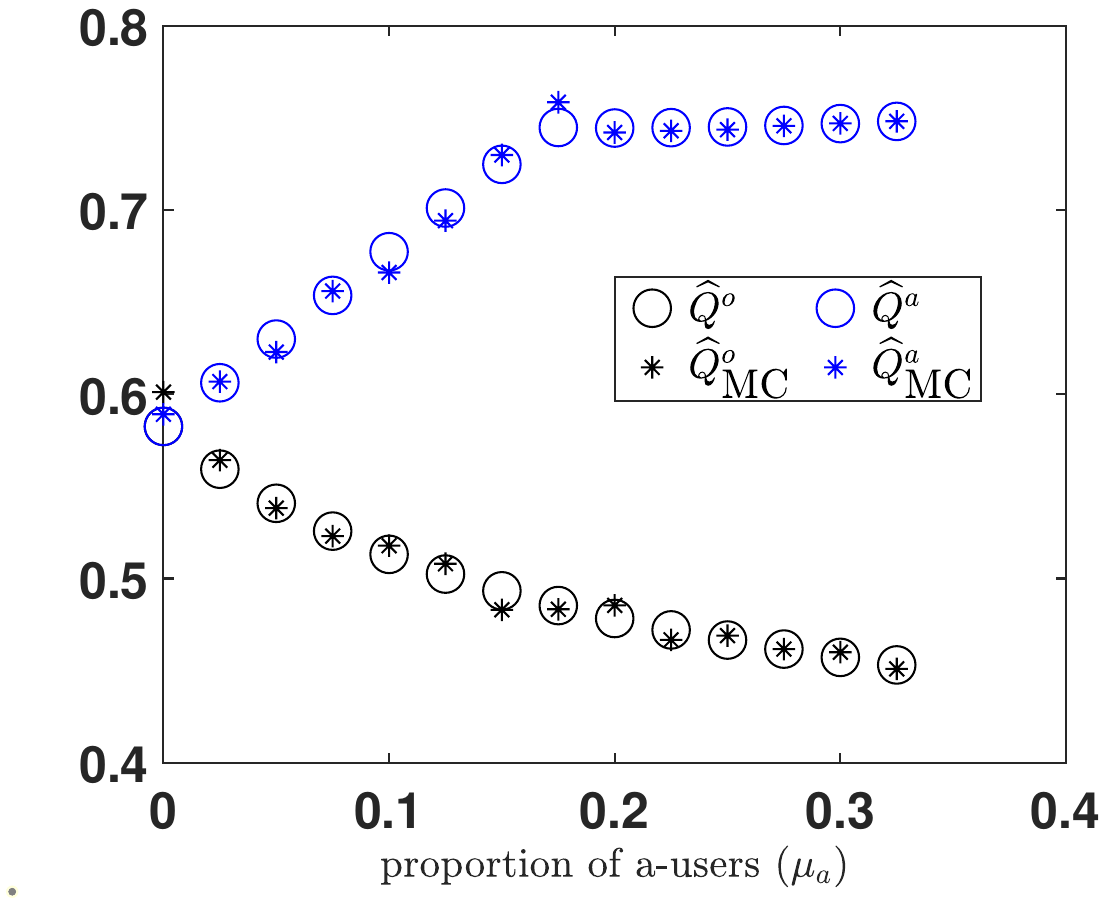}
\end{minipage}
\caption{Comparison of i-QoS under eo-WM and ea-WM, with smart (left) and naive (right) users respectively}
\label{fig:eaWM}
\end{figure}

Further as seen from the figure \eqref{fig:eaWM}, in all the case studies, the i-QoS improves; nonetheless this way of improvement does not degrade the performance of the real-post, as confirmed by Theorem \ref{corollary_ea_wm}) and also as observed in Figure \ref{fig:ehWM_real} which plots the performance for the real-post.
More interestingly,  the i-QoS and the improvement (with respect to eo-WM) both increase sharply with $\mu_a$.  
 Thus,  even in the presence of a larger fraction of a-users confusing the WM, ea-WM is able to nudge the non a-users to correctly identify the fake-post as fake. In view of Theorem \ref{corollary_ea_wm}, this may  be true as  ea-WM provides increasingly high warning levels with increase in $\mu_a$ (see \eqref{eqn_warning_ea}).  
 One probably can design a better WM  that  provides higher warning levels even with smaller value of $\mu_a$ (and which again ensures the required real performance) and the quest further is precisely for the same. 
 }

\revr{From Figure \ref{fig:eaWM}(left), for the case study  with smart users, observe that $\widehat{Q}^a = 1$, the maximum possible i-QoS, for $\mu_a \geq 0.05$. However, ea-WM fails to achieve such high i-QoS with naive users --- i-QoS is less than $0.8$  in right sub-figure of Figure \ref{fig:eaWM}. The quest again is for a better WM  which works well even for naive users, and this is considered in the immediate next. 


}

\section{Enhanced WM (eh-WM)} \label{sec_enWM}
In this section, we design an improved version of ea-WM. The idea is to design a warning $\omega^h$ such that $\omega^a(\beta) < \omega^h(\beta)$ for all $\beta \in [0,1]$. \revg{In lines of Theorem \ref{thrm_unique_att},} such monotonicity of the WM will ensure that the zeroes of the function $g_\beta^F$ (see \eqref{eqn_general_g_beta}) corresponding to the new WM are larger than that of $g_\beta^{a, F}$. However, the design should be such that the performance of the new WM for the real-post is not compromised. Towards this, we design an \underline{enhanced warning mechanism (eh-WM)} as follows:
\begin{align}\label{eqn_warning_eh}
    \omega^h(\beta)  = \phi \omega^a(\beta), \mbox{ for an appropriate choice of } \phi > 1, \mbox{ with } w, b \mbox{ as in ea-WM.}
\end{align}
For given $\phi$, denote the $g_\beta^u$ of \eqref{eqn_general_g_beta} corresponding to the eh-WM as $g_{\beta, \phi}^{h, u}$. Further, define $\beta^{h}_\phi$ as a zero of $g_{\beta, \phi}^{h, F}$ in $[0,1]$ and $\beta^{h, R}_\phi$ as a zero of $g_{\beta, \phi}^{h, R}$ in $[0,1]$. Observe that:

\vspace{-2mm}
{\small
\begin{align*}
    g_{\beta, \phi}^{h, F}(\beta) &= g_\beta^{a, F}(\beta) + \mu_2 m_f \eta^F \bigg\{ \beta \bigg( \min\{1, \phi \omega^a(\beta) \alpha_x^F\} - \min\{1, \omega^a(\beta) \alpha_x^F\} \bigg) + (1-\beta) \bigg( \min\{1, \phi \omega^a(\beta) \alpha_y^F\} - \min\{1, \omega^a(\beta) \alpha_y^F\} \bigg) \bigg\} \\
    &\geq g_\beta^{a, F}(\beta),
\end{align*}}with equality only if $\alpha_j^F \omega^a(\beta) > 1$ for each $j \in \{x, y\}$. This implies that any zero of $g_{\beta, \phi}^{h, F}$ is larger or equal to the smallest zero of $g_\beta^{a, F}$. Thus, it clear that $\beta_\phi^{h} \geq Q^a$ for any $\beta^h_\phi \in \cA^{h, F}_{\beta, \phi} \cup \cR^{h, F}_{\beta, \phi}$. Therefore, we have:
$$
\inf\{\beta : \beta  \in \cA^{h, F}_{\beta, \phi} \cup \cR^{h, F}_{\beta, \phi}\} =: Q^h_\phi \geq Q^a.
$$
That is, the QoS under eh-WM (for any $\phi$) is higher or at par with the QoS corresponding to ea-WM.

Now, one can anticipate that higher the warning level is, the more cautiously users tag the posts. Thus, as $\phi$ increases, the proportion of fake-tags must increase. However, one can not choose an arbitrarily large $\phi$ as then the performance for the real-post is degraded. Thus, we consider the following problem to optimally choose $\phi = \phi^*$ such that $Q^h_\phi$ is maximized, while satisfying constraint in \eqref{eqn_opt_prob}:
        \begin{align}
        \label{eqn_opt_phi}
            \max_\phi & \hspace{2mm} \revg{Q^h_\phi}
             \mbox{ subject to } \beta \leq \delta \mbox{ for each } \beta \in \cA^{h, R}_{\beta, \phi} \cup \cR^{h, R}_{\beta, \phi}.
        \end{align}
We have the following optimal design for the   eh-WM (proof is in  appendix):

\begin{theorem}\label{corollary_eh_WM}
Define the constant 
\begin{align}\label{eqn_phi_bar}
\overline{\phi} := \frac{\delta \bigg( \mu_2 \eta^R + \mu_1 (1-\alpha_x^R \rho) \eta_R + \mu_a \eta_a \bigg) - (1-\delta) \mu_1 \rho \alpha_y^R \eta_R}{\mu_2 \omega^a(\delta) \bigg( \delta \alpha_x^R + (1-\delta) \alpha_y^R \bigg) \eta^R}.
\end{align}
The $\phi^*$ defined below is greater than $1$ and is the optimizer of the problem \eqref{eqn_opt_phi}:
\begin{eqnarray}\label{eqn_phi_star}
    \phi^* :=
    \left \{ \begin{array}{ll}
       \overline{\phi}, & \mbox{ if \hspace{2mm}}  \overline{\phi} < \frac{1}{\alpha_y^R \omega^a(\delta)}, \mbox{ or if  \hspace{2mm}} \overline{\phi} \geq \frac{1}{\alpha_y^R \omega^a(\delta)}, \hspace{2mm} \underline{\beta}^F = 0 \mbox{ and } b = 0,\\
    \frac{1}{\omega^a (\underline{\beta}^F) \alpha^F_y}, & \mbox{ else.} 
     \end{array} \right.  \mbox{\eop}
\end{eqnarray}   
\end{theorem}
\revg{Thus, the choice of $\phi$ which gives the maximum proportion of fake-tags for the fake-post is given by $\phi^*$. Such a $\phi^*$ also ensures that the performance of eh-WM for the real-post is not degraded beyond $\delta$-level. The problem \eqref{eqn_opt_phi} can also be designed and solved in terms of the better metric i-QoS and by replacing $\delta$ by $\delta_a$ analogously. 
Henceforth, when we refer eh-WM, it corresponds to the case with $\phi = \phi^*$ and when $\delta = \delta_a$.
We present the numerical results with respect to eh-WM directly in terms of i-QoS and the correspondingly modified $\delta_a$-threshold.}

\subsection{Numerical analysis for eh-WM}
We now (MC) simulate the warning dynamics under eh-WM \revr{for the two examples with smart and naive users and the  MC-estimates again well match the theoretical values, as seen from  Figure \ref{fig:ehWM_real} (for real-post) and Figure \ref{fig:ehWM} (for fake-post).  
Next, we discuss the qualitative analysis. To begin with, the Figure \ref{fig:ehWM_real} re-affirms the results of Theorem \ref{corollary_eh_WM} with regard to  the real-post --- the proportion of fake-tags for the real-post is at most $\delta_a$.
\begin{figure}[http]
\centering
\begin{minipage}{.5\textwidth}
  \centering
  \includegraphics[trim = {1cm 6cm 0cm 6cm}, clip, scale = 0.3]{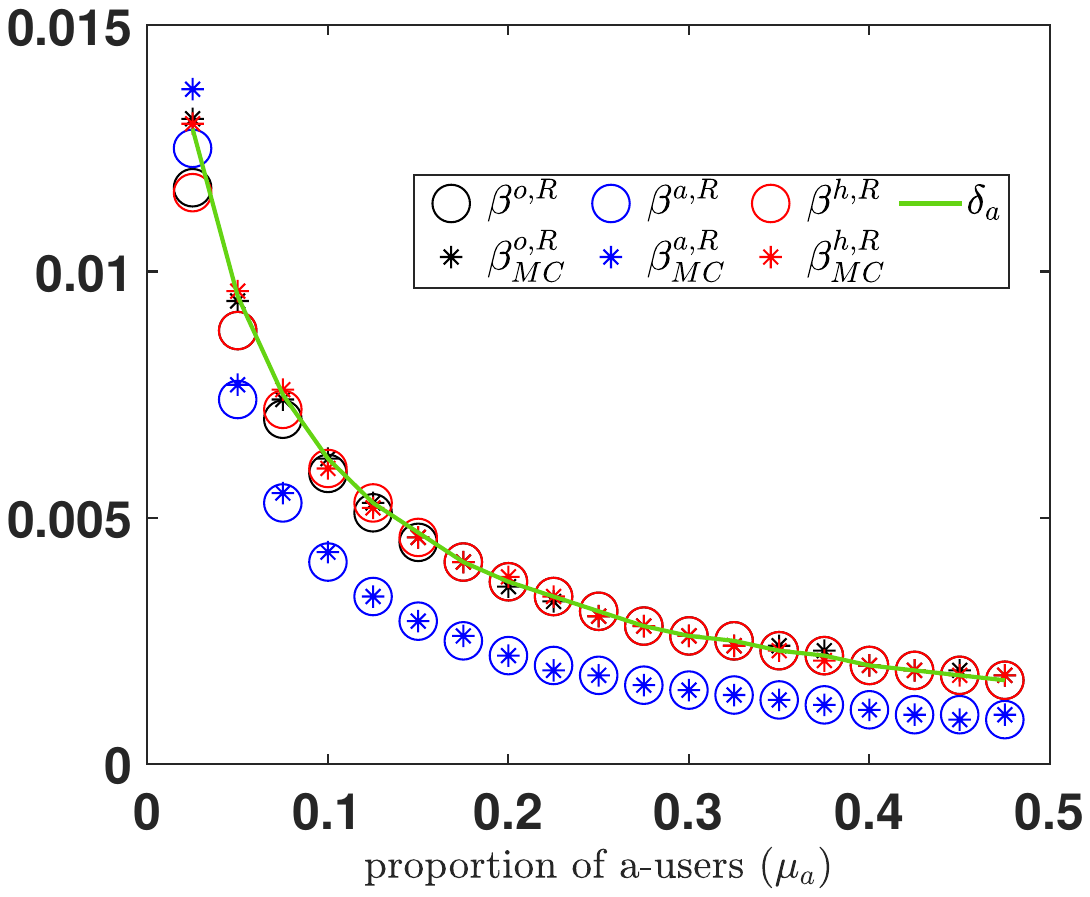}
\end{minipage}%
\begin{minipage}{.5\textwidth}
  \centering
  \includegraphics[trim = {1cm 6cm 0cm 6cm}, clip, scale = 0.3]{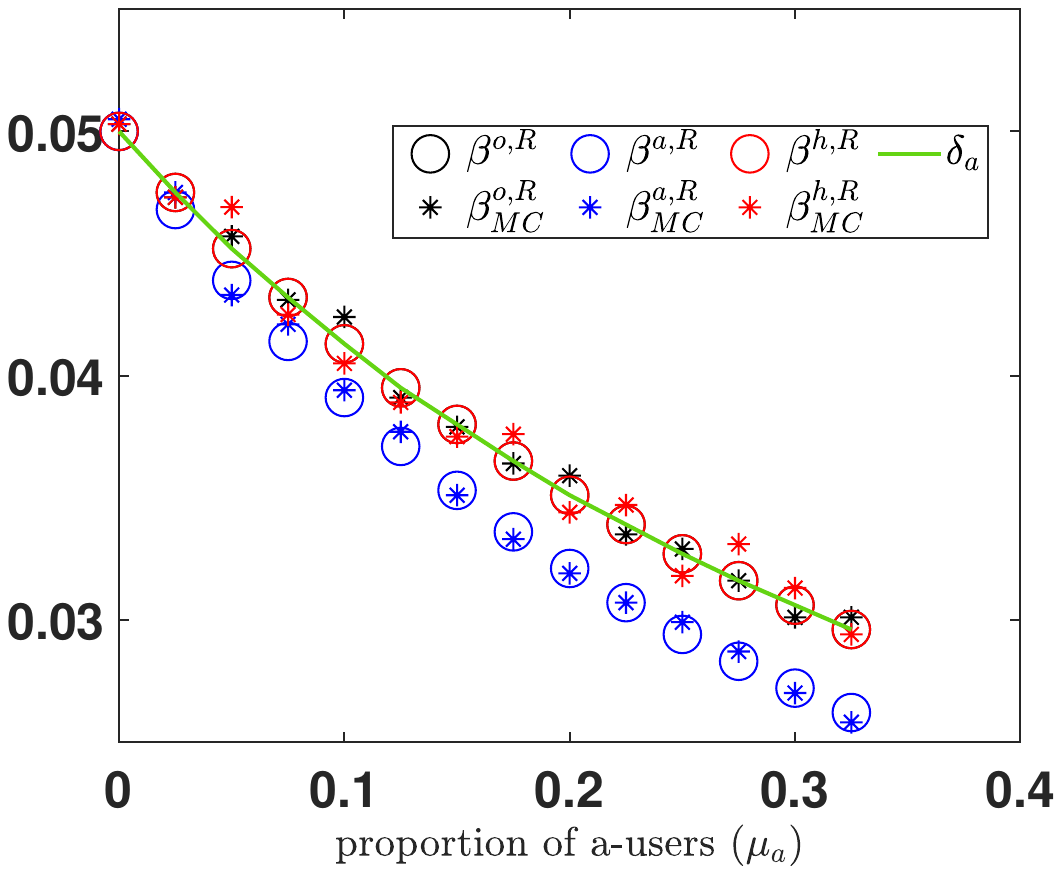}
\end{minipage}
\caption{Limits of warning dynamics for real-post under three WMs with smart (left) and naive (right) users  respectively}
\label{fig:ehWM_real}
\end{figure}

In Figure \ref{fig:ehWM}, we plot the i-QoS under eh-WM $\widehat{Q}^h$ (i.e., with $\phi^*$), along with that corresponding to the previous two WMs. For the example with smart users, eh-WM performs at par with ea-WM; recall, ea-WM almost achieved $\widehat{Q}^a = 1$. However, for the case with naive users, $\widehat{Q}^h \gg \widehat{Q}^a$; thus, eh-WM is more robust against adversaries than ea-WM. Therefore, one can say that eh-WM is able to guide the naive non a-users about the actuality of fake-posts better than ea-WM. }

\revr{As an example, when $10\%$ of a-users are trying to harm the system, the eh-WM ensures that $76.29\%$ of naive non a-users correctly identify the fake-post, while this fraction is only $67.73\%$ under ea-WM (observe, $\widehat{Q}^h - \widehat{Q}^a$ is as large as $0.0856$, for $\mu_a = 0.1$). 

\begin{figure}[http]
\centering
\begin{minipage}{.5\textwidth}
  \centering
  \includegraphics[trim = {1cm 6cm 0cm 6cm}, clip, scale = 0.3]{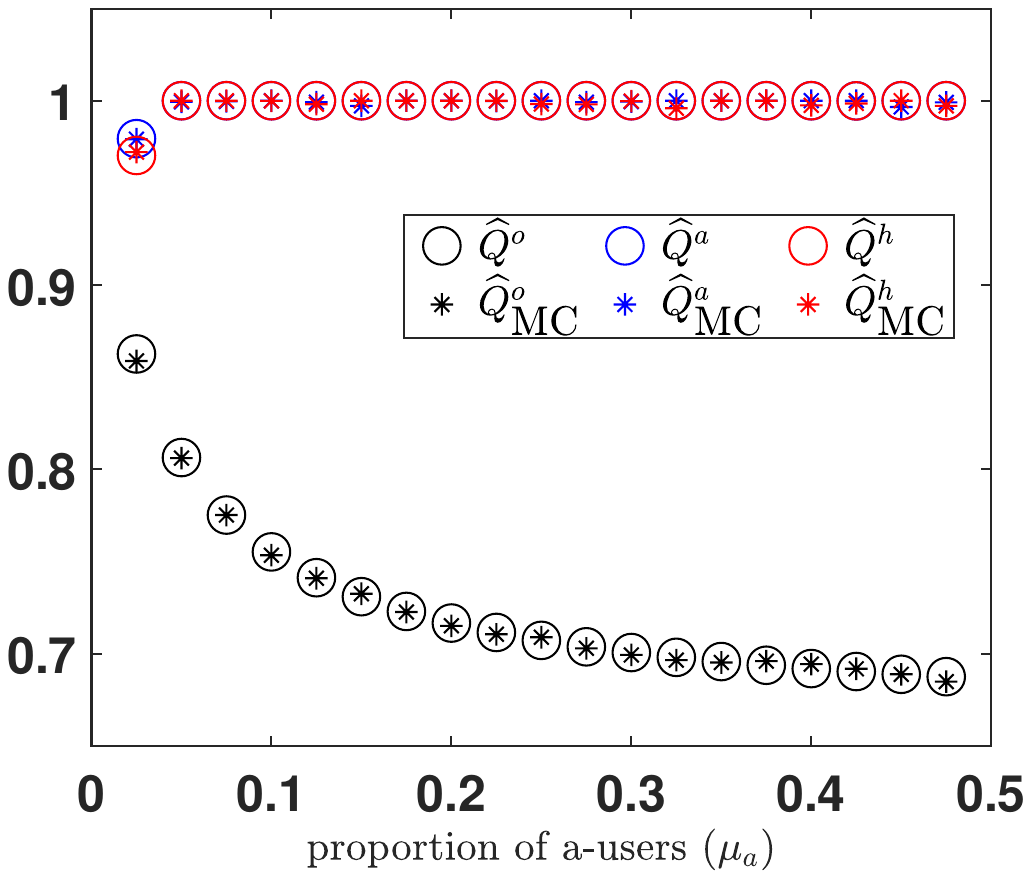}
\end{minipage}%
\begin{minipage}{.5\textwidth}
  \centering
  \includegraphics[trim = {1cm 6cm 0cm 6cm}, clip, scale = 0.3]{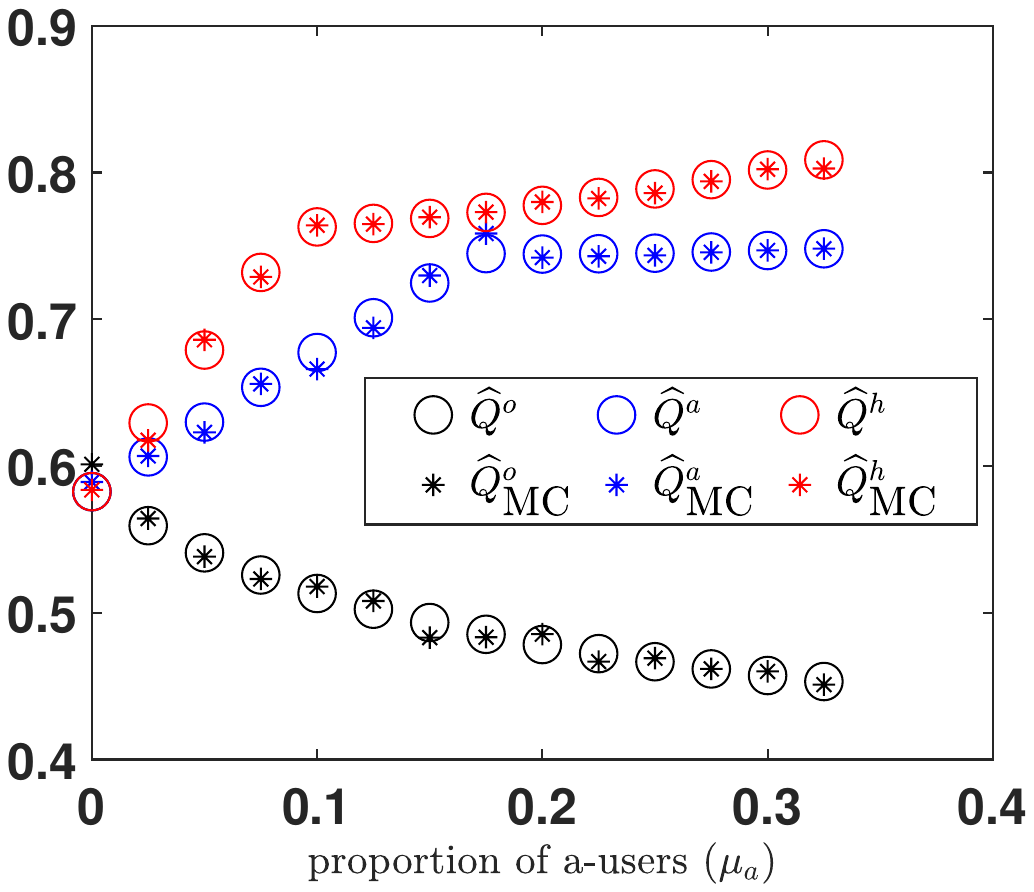}
\end{minipage}
\caption{Comparison of i-QoS under three WMs with smart (left) and naive (right) users  respectively}
\label{fig:ehWM}
\end{figure}
As seen from the example with naive users, eh-WM (red curve) performs significantly better than ea-WM (blue curve). Even then, the i-QoS under eh-WM is much better with higher values of $\mu_a$. This probably calls for a very different design of WM, which can generate high warning levels even for smaller values of $\mu_a$. This is attempted in the immediate next.



}

\revr{

\section[fgfdgfd]{Enhanced-2 WM (eh2-WM) and learning\footnote{We would also like to thank the reviewers of the paper, as their feedback motivated us to design a better warning mechanism (discussed in this section) that does not require the knowledge of system parameters and the proportions.}}
It is intuitive that as warning increases, the users are alarmed rigorously about the actuality of the posts; this should lead to more users correctly identifying the posts, and thus higher QoS; in fact, Theorem \ref{thrm_unique_att} precisely captures this intuition. If one can control the warning such that it does not harm the performance of the real-post beyond $\delta$-threshold, providing higher warning should be effective. We designed ea-WM and eh-WM along these intuitions with higher warning than eo-WM (recall, there is an additive term in \eqref{eqn_warning_ea} and multiplicative term in \eqref{eqn_warning_eh}), and still managed to ensure the performance of the real-post is within the desired level (see Theorem \ref{corollary_ea_wm} and Theorem \ref{corollary_eh_WM}). Further motivated by this, in this section, we aim to design another improved version of eo-WM, named \underline{enhanced-2 WM (eh2-WM)} and denoted by $\omega^{h2}$, which provides higher warning signals to the users (in fact, even for the cases with smaller $\mu_a$); this mechanism also facilitates learning the required parameters $b$ and $w$.

To achieve the same, we again utilize the eo-WM but now with a bigger $w$, and ensure that there is a unique limit proportion for the real-post which satisfies the $\delta$-threshold. From \eqref{eqn_warning}, a bigger $w$ results in higher warning levels, hence, we simply set $w = w^{h2} := \nicefrac{1}{\alpha_x^R} - \gamma$ and choose a corresponding $b$ as in Theorem \ref{thrm_opt}. This value of $w$ ensures that $\alpha_j^R \omega^{h2} (\beta) \le 1$ for all $j \in \{x,y\}$ and  all $\beta \in [0,1]$ for real post (i.e., when $u = R$) and hence using exactly the same logic as in Corollary \ref{corollary_ex_wm}, we have a unique zero/attractor for the real-post; further the choice of $b$ as in Theorem \ref{thrm_opt} ensures the said unique attractor $\beta^{R,h2}$ corresponding to the real-post is within the required threshold $\delta$. 
However, unlike eo-WM, with larger $w$ we may not have   a unique limit proportion for the fake-post under eh2-WM. Nonetheless, the resultant QoS (and hence i-QoS) is bigger than that with eo-WM by Theorem  \ref{thrm_unique_att}, as with bigger $w$, $\omega^{h2}(\beta) > \omega(\beta)$ for all $\beta$.  

It is important to observe here that the new enhanced WM (eh2-WM) generates high levels of warning signals, and its design does not depends on parameters like $\mu_a$. Thus, one can anticipate that it will enhance the performance even for the smaller values of $\mu_a$. To illustrate the same, we tabulate the i-QoS, $\widehat{Q}^{h2}(w^{h2}, b(w^{h2}))$, achieved under eh2-WM for the case with naive users:
\begin{table}[http]
    \centering
    \begin{tabular}{|c|c|c|c|c|}
    \hline
         & $\mu_a = 0$ &   $\mu_a = 0.1$ &  $\mu_a = 0.2$ &  $\mu_a = 0.3$\\ \hline
         $\widehat{Q}^{h2}(w^{h2}, b(w^{h2}))$ & $0.8289$ & $0.8270$ & $0.8257$ & $0.8246$\\ \hline
    \end{tabular}
    \caption{i-QoS under perfect knowledge of user sensitive parameters}\label{table_WM_perfect}
\end{table}
Clearly, the i-QoS under eh2-WM is consistently higher than that with eh-WM (see Figure \ref{fig:ehWM}, where the red curve is below $0.8$ for all $\mu_a$). More importantly, the i-QoS under eh2-WM is almost the same for all values of $\mu_a$.  


%

\noindent \textbf{Learning the parameters:} At this point, it is important to note that all the discussion so far assumed that the user sensitive parameters ($\rho$ and $(\alpha_i^u)$ for each  $i \in \{x, y\}$ and $u \in \{R, F\}$) and proportions of users of different types ($\mu_1, \mu_2$ and $\mu_a$) are known to the OSN. However, such information is not easily accessible to the OSN and the purpose now is to design a WM without such knowledge. Towards this, \textit{we propose an algorithm which directly learns the parameters of the WM, $b$ and $w$}. We only require that there is a non-zero proportion of ws-users\footnote{it can be checked by noticing the users who click on the information button (see Figure \ref{fig_post_design}) }, i.e., $\mu_2 > 0$ and the knowledge of ratio $\alpha_x^R/\alpha_y^R$ (details are given below). The design would only utilize various random quantities that are observed during the post propagation process.

The main idea is to consider a real-post which is known to the OSN and train the parameters $w$ and $b$ using the responses of the users.

Basically, we add a SA-based step which tunes $b$ such that the corresponding $\beta^{o, R}$ eventually approaches $\delta$ - recall, the constraint in optimization problem \eqref{eqn_opt_prob} requires that $\beta^{o, R} \leq \delta$. Further, $w$ is tuned such that $\alpha_y^R \omega^{h2}(1)$ approaches $ 1-\kappa$, where constant $\kappa \geq 1 -~\alpha_y^R/\alpha_x^R$. From \eqref{eqn_warning}, $\omega^{h2}(1; w, b) = w+ \gamma$, and hence such a tuning ensures that $w$ approaches $\nicefrac{(1-\kappa)}{\alpha_y^R} - \gamma$ (and by choice of $\kappa$, eventually $w \leq \nicefrac{(1-\kappa)}{\alpha_y^R} - \gamma$) --- thus, eventually $\alpha_j^R \omega^{h2}(1) \leq 1$ for each $j \in \{x, y\}$, as planned for the real-post. Here, we would like to stress that the tuning of $w$ is done with respect to $\alpha_y^R$, instead of $\alpha_x^R$, as there may not be sufficient estimates corresponding to fake-tags for the real-posts (recall, $\delta$ is typically a small value). Thus, the algorithm requires some idea on the ratio $\alpha_x^R/\alpha_y^R$. In all, if such a tuning (of both $w$ and $b$) is possible, then it would ensure a unique attractor below $\delta$-threshold for the real-post. 

\RestyleAlgo{ruled}
\begin{algorithm}
\caption{Design of learning WM}\label{alg_WM}
(i) Consider a real-post.

(ii) Initialize $\Cx(\tau_0)$ and $\Cy(\tau_0)$; calculate $B_0^{h2, R}$. Fix a large enough $\SampS < \infty$.

(iii) Initialize $b_0$ and $\eta_0$ sufficiently small, and choose a  $w_0 > 1$.

(iv) At $k$-th epoch, $\tau_k$, when $k$-th user reads the post, for $k \in \{1, 2, \dots, \SampS\}$:
\begin{itemize}
\item  set the $w$-update flag, $J_{ws} = 0$
    \item if the reader is a ws-user, then provide warning, $\omega^{h2}$, which is set as below:
        \begin{itemize}
            \item  toss a biased coin such that $P(\mbox{head appears}) = \eta_{k-1} > 0$, let $\eta_{k-1} \to 0$
            \item if head appeared and if the said user received with post with real-tag, 
            \begin{itemize}
                \item set warning corresponding to $\beta = 1$, i.e., set $\omega^{h2}(B^{h2, R}_{k-1}) := w_{k-1} + \gamma$
                \item set the indicator $J_{ws} = 1$
            \end{itemize}
             
            \item else, set warning as per WM, i.e., set $\omega^{h2}(B^{h2, R}_{k-1})$  as in \eqref{eqn_warning_learn}
        \end{itemize}
    \item observe the tag $I_k$ and the number of shares by the said user; accordingly, update proportion of fake tags, $B_k^{h2, R} = \frac{\Cx(\tau_{k-1}^+)}{\Cx(\tau_{k-1}^+) + \Cy(\tau_{k-1}^+)}$
    
    \item update the parameters, using the new estimate $B_k^{h2, R}$ and $I_k$
    \begin{itemize}
        \item if $w$-update flag, $J_{ws} = 1$, then update  $w_k$ as in  \eqref{update_w}
        \item update $b_k$ as in \eqref{update_b}
    \end{itemize}
\end{itemize}

\end{algorithm}
The above tuning for $w$ requires warning levels $\omega^{h2}(1)$, corresponding to $\beta = 1$; however, in the eo-WM, the warning levels were generated according to the then estimates of $\beta$, the proportion of fake-tags. To minimally disrupt the normal functioning of the WM, we propose some special epochs at which such special warning is provided -- at time epoch $k$, if a ws-user who received the post with real-tag clicks on the information button, the OSN generates such a warning with probability $\eta_k$, where $\eta_k \downarrow 0$, as $k \to \infty$. Only such special epochs are used to learn $w$. To summarize, the updates for $w$ at epoch $k$ are as follows: if a ws-user that received the post with real-tag reads the post, then we have:
\begin{align}\label{update_w}
    w_k \gets 
    \max \left\{1, w_{k-1} - \epsilon_k\left( I_{k} - (1-\kappa) \right) \right\}, &\mbox{ with probability }\eta_k, 
\end{align}where $I_{k}$ is the indicator that the user tags the post as fake and $\epsilon_k := c_1(\frac{1}{k+1})^{c_2}$ with some appropriate $c_1 > 0$ and $c_2 \in ( 0.5, 1]$. In all other cases, we set $w_k = w_{k-1}$.



Next, we discuss the updates for $b$.  For each $k \geq 1$, update $b_k$ as below:
\begin{align}\label{update_b}
\begin{aligned}
    b_k &\gets \max \left\{0, b_{k-1} + \epsilon_k (B_k^{h2, R} - \delta) \right\},  \mbox{ where as before } B_k^{h2, R} := \frac{\Cx(\tau_k^-)}{\Cx(\tau_k^-) + \Cy(\tau_k^-)},
\end{aligned}
\end{align}and the post-propagation process updates as in \eqref{eqn_transition_fake_tag} and \eqref{eqn_transition_real_tag} --- the warning shown to the $k$-th user reading the post would have been generated using $(w_k, b_k)$ as below: 
\begin{align}\label{eqn_warning_learn}
\omega^{h2}(B^{h2, R}_k) := \omega(B^{h2, R}_k) = \frac{w_k B^{h2, R}_k}{B^{h2, R}_k + b_k(1-B^{h2, R}_k)} +~\gamma, 
\end{align}
at the normal epochs (when $w_k$ is not updated); for the special epochs, the warning $\omega^{h2}(B^{h2, R}_k) := \omega(1) = w_{k} + \gamma$ is generated.


The brief idea behind such a design is that as is usually the case with SA algorithms, the SA iterates $b_k$ and  $w_k$ converge so as to ensure the expected values of the respective  update-terms $B_k^{h2, R} - \delta$ in \eqref{update_b} and  $I_{k} - (1-\kappa)$ in \eqref{update_w} converge to $0$ as $k \to \infty$. That is, $\beta^{h2,R}_k = E[B_k^{h2, R} ]$ approaches $\delta$ and $w_k$ approaches\footnote{Observe that the the conditional expected value conditioned that the user is a ws-user who received the post with real-tag, $E[I_{ k}] = \alpha_y^R \omega(1) = \alpha_y^R  (w_k + \gamma) $.} $(1-\kappa)/\alpha_y^R - \gamma$. As already mentioned, such a limit of $w$ ensures that the unique limit for the real-post $\beta^{h2,R}$ is near $\delta$; thus, the constraint in \eqref{eqn_opt_prob} is satisfied and the discussion in the beginning of this section also ensures that the QoS is strictly improved in comparison to the eo-WM.

The learning algorithm is summarized in Algorithm \ref{alg_WM}. The analysis of the above learning algorithm  would require rigorous two-time scale (projected) SA-based tools - observe $w_k$ is updated minimally and further probability $\eta_k \downarrow 0$. We skip the analysis here, but validate and illustrate the improved performance of the \underline{learning WM} (referred to as l-eh2-WM)  via numerical examples in the next sub-section.



\subsection{Numerical analysis for l-eh2-WM}\label{subsec_numerical_estimate}
In Table \ref{table:WM}, we continue with the example with naive users to test the learning algorithm. Towards this, we fix   $\kappa = 1 - \nicefrac{\alpha_y^R}{\alpha_x^R} + 10^{-3}$, $\eta_k = 1.5(\nicefrac{1}{k})^{0.8}$, $\eta_0 = 0.008$, $w_0 = 6$ and $b_0 = 10^{-4}$. The choice of $\epsilon_k$ for learning $b$ and $w$ is $2.2(\nicefrac{1}{k})^{0.7}$. We initialize the system such that a real-post is shared by the content provider to $20$ users with the real-tag.

For a given sample-size (number of samples available for learning and represented by $\SampS$), we consider $150$ sample paths for the post-propagation of the real-post under l-eh2-WM; the idea is to  measure  the  efficacy of l-eh2-WM algorithm  via the fraction of times it achieves an i-QoS   within $\pm0.05$ of that corresponding to the case with perfect information (i.e., $\widehat{Q}^{h2}(w^{h2}, b(w^{h2}))$). 
We consider different sample-sizes $\SampS$ in the range $10^4$ to $10^5$.
In each sample-path, l-eh2-WM algorithm  is used to update the estimates of $\{(w_k, b_k)\}$  for $\SampS$ number of epochs and then the i-QoS $\widehat{Q}^{o}_k(w_{\SampS}, b_{\SampS})$  (fake-post) corresponding to eo-WM using the last estimate $(w_{\SampS},  b_{\SampS})$ is computed.

In Table \ref{table:WM}, for different values of $\SampS$, we tabulate $f_\SampS$, the fraction of sample paths for which $|\widehat{Q}^{h2}(w^{h2}, b(w^{h2}))- \widehat{Q}^{o}_\SampS(b_\SampS, w_\SampS)| \leq 0.05$.

\begin{table}[http]
\centering
\begin{tabular}{|c|ccccc|}
\hline
                              & \multicolumn{5}{c|}{$\SampS$ }                                                                                                                                                                           \\ \hline
                                 & \multicolumn{1}{c|}{$10^4$} & \multicolumn{1}{c|}{$2.5* 10^4$} & \multicolumn{1}{c|}{$5*10^4$} & \multicolumn{1}{c|}{$7.5*10^4$} & \multicolumn{1}{c|}{$10^5$} \\ \hline
\multirow{1}{*}{$\mu_a = 0$}  
                               & \multicolumn{1}{c|}{0.73}    & \multicolumn{1}{c|}{0.89}      & \multicolumn{1}{c|}{0.91}    & \multicolumn{1}{c|}{0.95}    &  \multicolumn{1}{c|}{0.93}
                               \\ \hline
\multirow{1}{*}{$\mu_a = 0.1$} 
                               & \multicolumn{1}{c|}{0.41}    & \multicolumn{1}{c|}{0.57}      & \multicolumn{1}{c|}{0.76}    & \multicolumn{1}{c|}{0.84}      & \multicolumn{1}{c|}{0.91}    \\ \hline
\multirow{1}{*}{$\mu_a = 0.2$} 
                               &  \multicolumn{1}{c|}{0.19}    & \multicolumn{1}{c|}{0.44}      & \multicolumn{1}{c|}{0.64}    & \multicolumn{1}{c|}{0.74}      & \multicolumn{1}{c|}{0.79}           \\ \hline
\end{tabular}
\caption{Fraction of sample paths that learnt the parameters ($b, w$) sufficiently well and achieved the desired level of i-QoS under l-eh2-WM}
\label{table:WM}
\end{table}
It can be seen from the table that the fraction of sample paths with the desired property ($f_\SampS$) increases with $\SampS$, thus depicting that the l-eh2-WM is progressively able to achieve the performance close to the case with perfect knowledge. One may anticipate that more iterations/shares should be required to achieve i-QoS of eh2-WM (i.e., with perfect knowledge) as $\mu_a$ increases; the same is evident from the table; for example, when $\SampS = 10^5$, $f_\SampS$ is as large as $0.91$ for $\mu_a = 0.1$, but with $\mu_a = 0.2$, it is much smaller and equals $0.79$. 
 Thus, this example illustrates that l-eh2-WM has managed to learn and tune the WM sufficiently well, when number of samples $ \geq 7.5*10^4$ for proportion of a-users up to $0.2$. 

The performance of the l-eh2-WM is sensitive to the initial conditions and the parameters of the two-timescale algorithm (like, $\epsilon_k$), as is the usual case with SA-based algorithms. Using trial-and-error method, we picked a good enough set of values, while extensive study on better choice of these parameters is outside the scope of this work. 


Next, in Figure \ref{fig:estimate}, we continue with
the two examples considered in Figure \ref{fig:eoWM_delta}. In the left and right sub-figures, we consider the instance with smart and naive users respectively and present the results directly in terms of i-QoS. The learning algorithm is again initialized and tuned appropriately, and now with a large sample-size, $\SampS = 10^6$.

\begin{figure}[http]
\centering
\begin{minipage}{.5\textwidth}
  \centering
  \includegraphics[trim = {1cm 6cm 0cm 6cm}, clip, scale = 0.3]{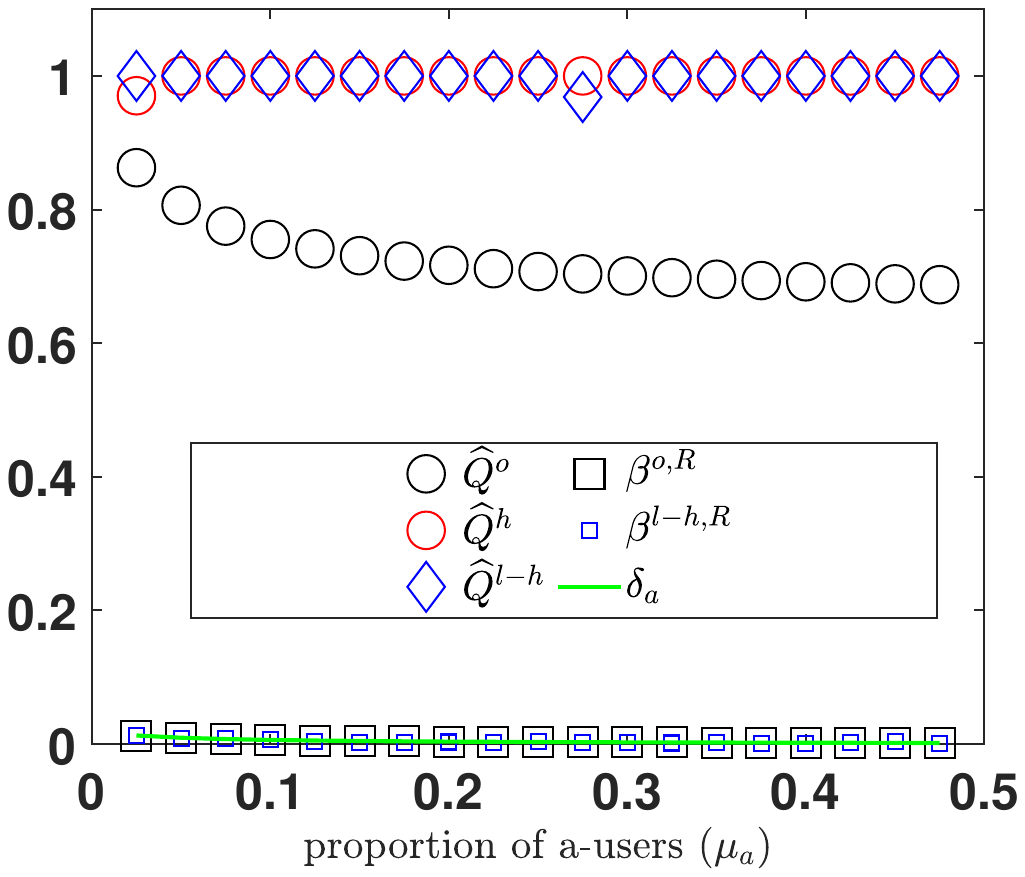}
\end{minipage}%
\begin{minipage}{.5\textwidth}
  \centering
  \includegraphics[trim = {1cm 6cm 0cm 6cm}, clip, scale = 0.3]{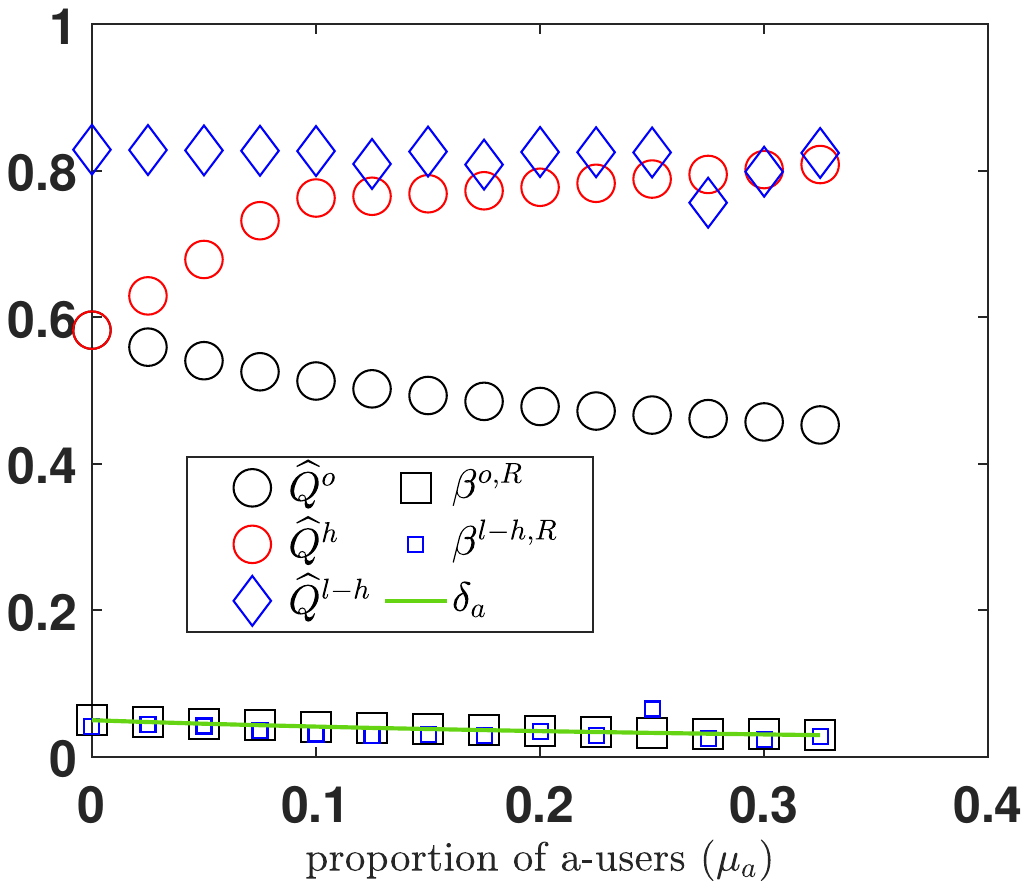}
\end{minipage}
\caption{Comparison of limits of warning dynamics under eo, eh and l-eh2-WM with smart (left) and naive (right) users  respectively}
\label{fig:estimate}
\end{figure}
From the figure, it can be seen that for all values of $\mu_a$, the i-QoS under l-eh2-WM  (marked in diamond) is higher than the eo-WM; in fact, it performs superior to all the previous WMs. Of course, the i-QoS can not be further improved for smart users --- even l-eh2-WM achieves i-QoS close to $1$, as eh-WM. The superior performance of eh2-WM (actually that of l-eh2-WM with large $\SampS$) is clearly depicted in the case with naive users. 
From Table \ref{table_WM_perfect} and Figure \ref{fig:ehWM}, it is clear that the eh2-WM outperforms eh-WM and further has similar performance for all values of $\mu_a$. The l-eh2-WM  with large $\SampS$ has exactly similar performance traits, as can be seen from Figure \ref{fig:estimate}.
Furthermore, the proportion of fake-tags for the real-post is also within the $\delta_a$-threshold, thus satisfying the constraint in \eqref{eqn_new_opt}.
}

\section{Conclusions and Future Work}
There is a huge requirement to identify fake posts on ever active OSNs. Further, any algorithm attempting to identify fake posts faces challenges from adversarial users and users unwilling to participate. Our first aim in this paper is to derive the performance of a promising recently proposed algorithm in the presence of adversaries who always real-tag any post. A severe degradation in performance is observed with just 1\% of adversaries.

The algorithm collects binary signals (fake/real tags) from all previous users, generates a warning based on the fraction of fake tags and compels further users to judge and consume the post cautiously based on the warning level provided. Using new results in branching processes (also derived in this paper), we obtain a one-dimensional ordinary differential equation (ODE) that analyses any generic iterative warning mechanism depending upon the fraction of fake tags. 
This ODE is instrumental in deriving robust adaptations of the previous mechanism -- in particular, we use concepts like eliminating the effects of adversaries,  the inherent monotone characteristics of relevant performance on certain parameters, etc. 
The new mechanisms illustrate significant performance improvement both in the presence and the absence of adversaries compared to the previous method. An algorithm which provides improvement over the existing method, without relying on the usually inaccessible users-specific information, is also proposed.

This paper also contributes towards total-current population-dependent two-type branching processes with population dependent death rates and also considers a variety of unnatural deaths. In particular, we derive all possible limits and limiting behaviours of the population sizes as time progresses.

In future, one can think of several new directions. The one-dimensional ODE can also be utilized to study other types of adversaries, like always fake-tagging adversaries or more informed adversaries that mis-tag both posts (fake-tag authentic post and real-tag the fake post). One can again derive improved algorithms, as we illustrated with real-tagging adversaries. One can also study the influence of users that share but refuse to tag or other important behavioural characteristics. \revr{Further, we designed two types of enhanced warning mechanisms, both of which improved over the existing mechanism. However, the two new mechanisms are not comparable, as one can perform better than the other in some instances. In future, one can attempt to design a combination of the two, which outperforms all of them, and also design the corresponding learning mechanism.}

\newcommand{\ignore}[1]{}

\ignore{

\subsection{Future Work}
\begin{itemize}
    \item We leave for future to capture the impact of users who do not tag but share the post - will lead to three type BP
    \item reluctance factor
    \item Unlike RC-a-users, a mis-guiding adversary (MG-a-user) is the one who knows the actuality of the post and therefore, always tag the post opposite to its true identity. The objective of such users is to completely confuse the system about the actuality of the posts. 
\end{itemize}

\newpage  }
 \appendix 
 
\section{} \label{appendix_journal2}
\begin{pfthrmone}
The proof follows exactly as in \cite{agarwal2021new}, except for some changes due to unnatural deaths. Here, we directly mention the SA based scheme for the new process, and necessary details where ever required. 

From \eqref{evolve_x_up_time_prop}, the embedded process immediately after $n$-th death, when say an $x$-type individual $d$-dies, is given by:
\begin{equation}\label{evolve_x_up_gen_prop}
\begin{aligned}
C_{x, n} &= C_{x, n-1}  + \offs_{xx, d, n}(\Om_{n-1}) - 1, \ \ \  T_{x, n} = T_{x, n-1}  + \offs_{xx, d, n}(\Om_{n-1}), \\
C_{y,n} &= C_{y, n-1} + \offs_{xy, d, n}(\Om_{n-1}),  \ \ \ A_{y, n} = A_{y, n-1} + \offs_{xy, d, n}(\Om_{n-1}).
\end{aligned}
\end{equation}

The ratios in $\Ups_n$ can be re-written as (with $\epsilon_{n-1} := 1/n$): 
\begin{align}\label{eqn_SA_scheme_prop}
\begin{aligned}
\Ups_n &= \Ups_{n-1} + \epsilon_{n-1}\mathbf{L}_{n-1}, \mbox{ where } \mathbf{L}_{n-1} := (L_{n-1}^{\psi, c}, L_{n-1}^{\theta, c}, L_{n-1}^{\psi, a}, L_{n-1}^{\theta, a})^t, \mbox{ with}\\
L_{n-1}^{\psi, c} &:=  \left\{ \sum_{d \in [d_x]} \left( H^x_{n,d} (\offs_{x, d, n}(\Om_{n-1})-1) \right) + \sum_{d \in [d_y]} \left( H^y_{n,d} (\offs_{y, d, n}(\Om_{n-1}) - 1)  \right)  \right\} 1_{\Pc_{n-1} > 0}   -  \Pc_{n-1}, \\
L_{n-1}^{\theta, c} &:= \left\{ \sum_{d \in [d_x]} \left( H^x_{n,d} (\offs_{xx, d, n}(\Om_{n-1})-1) \right) + \sum_{d \in [d_y]} \left( H^y_{n,d} \offs_{yx, d, n}(\Om_{n-1})  \right) \right\}1_{\Pc_{n-1} > 0} - \Tc_{n-1},  \\
L_{n-1}^{\psi, a} &:=  \bigg\{ \sum_{d \in [d_x]} \left( H^x_{n,d} \offs_{x, d, n}(\Om_{n-1}) \right) + \sum_{d \in [d_y]} \left( H^y_{n,d} \offs_{y, d, n}(\Om_{n-1})   \right)  \bigg\}1_{\Pc_{n-1} > 0}  - \Pa_{n-1}, \mbox{ and}\\
L_{n-1}^{\theta, a} &:= \bigg\{ \sum_{d \in [d_x]} \left( H^x_{n,d} \offs_{xx, d, n}(\Om_{n-1}) \right) + \sum_{d \in [d_y]} \left( H^y_{n,d} \offs_{yx, d, n}(\Om_{n-1})  \right) \bigg\}1_{\Pc_{n-1} > 0}  - \Ta_{n-1}, \mbox{ where}\\
\offs_{x, d, k} &:= \offs_{xx, d, k} + \offs_{xy, d, k}, \ \  
    \offs_{y, d, k} := \offs_{yy, d, k} + \offs_{yx, d, k},
\end{aligned}
\end{align}$H^x_{k,d} \in \{0, 1\}$ indicates that an $x$-type individual $d$-dies at $k$-th epoch such that $\sum_{d \in D_x} H^x_{k,d} \in \{0, 1\}$ and $\sum_{d \in D_y} H^y_{k,d} := 1 - \sum_{d \in D_x} H^x_{k,d}$. 

Henceforth, the proof of part (i) has two major steps: (a) to construct a sequence of piece-wise constant interpolated trajectories for almost all sample-paths; (b) to prove that the designed trajectories are equicontinuous in extended sense. We will provide the proof in terms of  $\tc$-component of the vector $\ups$, when the proof for the remaining components goes through in exactly similar manner.

Define $\gna = (\rho_\psi^c, \rho_\theta^c, \rho_\psi^a, \rho_\theta^a)$ as the conditional expectation, $E[\mathbf{L}_n|\mathcal{F}_n] =: \gna(\Ups_n, t_n)$, with respect to the sigma algebra, ${\cal F}_n  := \sigma\{\Om_k : 1 \leq k < n \}$ (see \cite[(16)]{agarwal2021new}).
Let $\Ups^n(\cdot) := (\Psi^{n, c}(\cdot), \Theta^{n, c}(\cdot), \Psi^{n, a}(\cdot), \Theta^{n, a}(\cdot))$  be  the constant piece-wise interpolated trajectory defined as below (see \eqref{eqn_SA_scheme_prop}, and recall $t_n = \sum_{i=1}^n \epsilon_{i-1}$):
\begin{align}\label{eqn_interpolated_traj_1_prop}
    \begin{aligned}
        \Theta^{n, c}(t) &:= \Tc_n + \int_0^t g_\theta^c(\Ups^n(s)) ds +  \sum_{i=n}^{\eta(t_n+t)-1} \epsilon_i L_i^{\theta, c} - \int_0^t g_\theta^c(\Ups^n(s)) ds\\
    &=  \Tc_n + \int_0^t g_\theta^{c}(\Ups^n) ds + M^{n, \theta, c}(t) + R^{n, \theta, c}(t) + D^{n, \theta, c}(t), \mbox{ where}\\
    M^{n, \theta, c}(t) &:=  \sum_{i=n}^{\eta(t_n + t)-1}  \epsilon_i \left(L_i^{\theta, c} - \rho_\theta^{c}(\Ups_i, t_i)\right), \\
    R^{n, \theta, c}(t) &:=  \sum_{i=n}^{\eta(t_n + t)-1}\epsilon_i g_\theta^{c}(\Ups_i) - \int_0^t g_\theta^{c}(\Ups^n) ds,  \\
    D^{n, \theta, c}(t) &:= \sum_{i=n}^{\eta(t_n + t)-1} \epsilon_i D_i^{\theta, c}, \mbox{ where } D_i^{\theta, c}:=\rho_\theta^{c}(\Ups_i, t_i) - g_\theta^{c}(\Ups_i),
    \end{aligned}
\end{align}
$\Psi^{n, c}(t), \Psi^{n, a}(t)$ and $\Theta^{n, a}(t)$ are defined analogously. As in \cite{agarwal2021new}, the  extended equicontinuity can be proved for $M^{n, \theta, c}(\cdot)$, $R^{n, \theta, c}(\cdot)$. For, $D^{n, \theta, c}(\cdot)$ 
 the procedure again follows as in \cite{agarwal2021new}  when $S_n \to 0$; however for sample paths where $S_n \nto 0$, the  arguments for proving the equicontinuity for $D^{n, \theta, c}(\cdot)$ slightly changes as below:
\begin{align}\label{eqn_bound_Di}
\begin{aligned}
|D_i^{\theta, c}| &\leq  |f_{\beta}(\Om_i) (m_{xx}(\Om_i) -1) - f_{\beta}^\infty(\Bc_i) (m_{xx}^\infty(\Bc_i)-1)| + | (1-f_{\beta}(\Om_i)) m_{yx}(\Om_i) - (1-f_{\beta}^\infty(\Bc_i))m_{yx}^\infty(\Bc_i)|\\
&\leq |f_{\beta}(\Om_i) m_{xx}(\Om_i)  - f_{\beta}^\infty(\Bc_i) m_{xx}^\infty(\Bc_i)| + |f_{\beta}(\Om_i) - f_{\beta}^\infty(\Bc_i)|\\
&\hspace{2cm}+ |m_{yx}(\Om_i) - m_{yx}^\infty(\Bc_i)|  + |f_{\beta}(\Om_i) m_{yx}(\Om_i)  - f_{\beta}^\infty(\Bc_i) m_{yx}^\infty(\Bc_i)|
\end{aligned}
\end{align}
In the above, under \ref{a2_prop}, the third term is bounded above by $1/(S_i)^\alpha$. The second term can be bounded above as follows:
\begin{align*}
    |f_{\beta}(\Om_i) - f_{\beta}^\infty(\Bc_i)| &= \Bc_i \left|\frac{\sum_{d \in D_x} \lambda_{x, d}(\Om_i)}{d(\Om_i)} - \frac{\sum_{d \in D_x} \lambda_{x, d}^\infty(\Bc_i)}{d^\infty(\Bc_i)} \right | \\
    &\hspace{-2cm} = \Bc_i \left|\frac{\sum_{d \in D_x} \lambda_{x, d}(\Om_i)}{d(\Om_i)} - \frac{\sum_{d \in D_x}  \lambda_{x, d}^\infty(\Bc_i)}{d(\Om_i)} + \frac{\sum_{d \in D_x}  \lambda_{x, d}^\infty(\Bc_i)}{d(\Om_i)} - \frac{\sum_{d \in D_x} \lambda_{x, d}^\infty(\Bc_i)}{d^\infty(\Bc_i)} \right |\\
    &\hspace{-2cm}\leq \frac{\Bc_i}{d(\Om_i)} \sum_{d \in D_x} |\lambda_{x, d}(\Om_i) - \lambda_{x, d}^\infty(\Bc_i)| + \Bc_i\left|  \sum_{d \in D_x}\lambda_{x, d}^\infty(\Bc_i) \left( \frac{1}{d(\Om_i)} - \frac{1}{d^\infty(\Bc_i)}\right)\right|\\
    &\hspace{-2cm}\leq \frac{\Bc_i}{d(\Om_i)} \left(\frac{|D_x|}{(S_i^c)^\alpha} +  \frac{\left|  \sum_{d \in D_x}\lambda_{x, d}^\infty(\Bc_i)\right|}{d^\infty(\Bc_i)} \left|d^\infty(\Bc_i) - d(\Om_i) \right|  \right)\\
    &\hspace{-2cm}\leq \frac{\Bc_i}{d(\Om_i)} \left(\frac{|D_x|}{(S_i^c)^\alpha} +  \frac{\left|  \sum_{d \in D_x}\lambda_{x, d}^\infty(\Bc_i)\right|}{d^\infty(\Bc_i)} \frac{(|D_x| + |D_y|)}{(S_i^c)^\alpha}  \right) \\
    &\hspace{-2cm}\leq \frac{\Bc_i}{d(\Om_i)} \frac{|D_x| + |D_y|}{(S_i^c)^\alpha} \left(1 + \frac{\left|  \sum_{d \in D_x}\lambda_{x, d}^\infty(\Bc_i)\right|}{d^\infty(\Bc_i)}  \right)\\
    &\hspace{-2cm}\leq \frac{\Bc_i}{d(\Om_i)} \frac{|D_x| + |D_y|}{(S_i^c)^\alpha} \left(1 + \frac{1}{B^c_i} \right) \hspace{5mm}\left(\mbox{ since } d^\infty(\Bc_i) \geq B_c^i  \sum_{d \in D_x}\lambda_{x, d}^\infty(\Bc_i)\right)\\
    &\hspace{-2cm}= \frac{(|D_x| + |D_y|)(B_i^c + 1)}{d(\Om_i) (S_i^c)^\alpha} \leq \frac{2(|D_x| + |D_y|)}{d(\Om_i) (S_i^c)^\alpha}
\end{align*}
Define $\Delta_1 := \min\left\{\inf_{\Om} \lambda_{x, d}(\Om), \inf_{\Om} \lambda_{y, d}(\Om) \right\} > 0$, by \ref{a1_prop}. Then, $d(\Om_i) \geq \Bc_i \inf_{\Om} \lambda_{x, d}(\Om) + (1-\Bc_i) \inf_{\Om} \lambda_{y, d}(\Om) \geq \Delta_1$. Thus, we have:
\begin{align}\label{eqn_sec_term}
    |f_{\beta}(\Om_i) - f_{\beta}^\infty(\Bc_i)| &\leq \frac{2(|D_x| + |D_y|)}{(S_i^c)^\alpha} \frac{1}{\Delta_1}
\end{align}
The first term in \eqref{eqn_bound_Di} can be bounded as follows  under \ref{a2_prop} and \eqref{eqn_sec_term}:
\begin{align*}
    |f_{\beta}(\Om_i) m_{xx}(\Om_i) - f_{\beta}^\infty(\Bc_i) m_{xx}^\infty(\Bc_i)|  &\leq |f_{\beta}(\Om_i)| |m_{xx}(\Om_i) - m_{xx}^\infty(\Bc_i)| + |m_{xx}^\infty(\Bc_i)| |f_{\beta}(\Om_i) - f_{\beta}^\infty(\Bc_i)|\\
    &\leq \frac{1}{(S_i)^\alpha} + \frac{2(|D_x| + |D_y|)}{(S_i^c)^\alpha} \frac{1}{\Delta_1} \left( E[\overline{\Gamma}] + \frac{1}{(S_i)^\alpha}\right).
\end{align*}
Similarly, the fourth term in \eqref{eqn_bound_Di} can be upper bounded as follows:
\begin{align*}
    |f_{\beta}(\Om_i) m_{yx}(\Om_i) - f_{\beta}^\infty(\Bc_i) m_{yx}^\infty(\Bc_i)| 
    &\leq \frac{1}{(S_i)^\alpha} + \frac{2(|D_x| + |D_y|)}{(S_i^c)^\alpha} \frac{1}{\Delta_1} \left( E[\overline{\Gamma}] + \frac{1}{(S_i)^\alpha}\right).
\end{align*}
Thus, $D_i^{\theta, c}$ can be upper bounded as follows for some $K <\infty$ (recall, $\alpha \geq 1$):
\begin{align*}
    D_i^{\theta, c} &\leq 2\left(\frac{1}{(S_i)^\alpha} + \frac{2(|D_x| + |D_y|)}{(S_i^c)^\alpha} \frac{1}{\Delta_1} \left( E[\overline{\Gamma}] + \frac{1}{(S_i)^\alpha}\right)\right) + \frac{1}{(S_i^c)^\alpha} + \frac{2(|D_x| + |D_y|)}{(S_i^c)^\alpha} \frac{1}{\Delta_1} \\
    &\leq \frac{K}{(S_i^c)^\alpha} \leq \frac{K}{S_i^c} = \frac{K}{\Pc_i \eta(t_i)} \leq \frac{K}{\Delta i}.
\end{align*}
This implies that, (recall $\epsilon_i = 1/(i+1)$ and $\alpha \geq 1$)
\begin{align*} 
|D^{n, \theta, c}(t)| = \left|\sum_{i= n}^{\eta(t_n + t) - 1}\epsilon_i D_i^{\theta, c}\right| \leq \sum_{i= n}^{\eta(t_n + t) - 1} \frac{K}{\Delta i (i+1)} \leq \sum_{i= n}^{\infty} \frac{K}{\Delta i (i+1)}, \mbox{ for any }t.
\end{align*}Thus, $D^{n, \theta, c}(t)$ uniformly converges to $0$ as $n \to \infty$.  In all, $(\Theta^{n, c}(\cdot))$ is equicontinuous in the extended sense. 

The proof of part (ii) follows exactly as in \cite{agarwal2021new}. \eop
\end{pfthrmone}



\begin{pfthrmtwo}
Observe that each point $x_i^* \in {\cal I}$ can either be a point of dis-continuity or continuity for $g_\beta$.
In the former case, when $x_i^*$ is either an attractor or repeller of the ODE \eqref{eqn_beta_ode_simple_prop}, the result can be proved exactly as in \cite[Theorem 2]{agarwal2021new}. In fact, when $x_i^*$ is a saddle point of the ODE \eqref{eqn_beta_ode_simple_prop}, the analysis can be easily extended similar to the case when $x_i^*$ is a repeller.

Now consider $x_i^* \in {\cal I}$ such that $g_\beta$ is continuous at $x_i^*$.
Let $\ups(0) \in \cD_I$ with $\pc(0) > 0$.  By \cite[Lemma 5.]{agarwal2021new}, $\pc(t) > 0$ for all $t \geq 0$, thus ODE \eqref{eqn_ODE_prop} simplifies to $\dot{\ups} = \mathbf{h}(\beta(\ups)) - \ups$. Now, we will prove the claim for different possibilities of $x^*$ as in the hypothesis separately. 
Firstly for all cases global solution exists because of Lipschtiz continuity.

\underline{Part (i)} Without loss of generality, let $\beta(0) \in {\cal N}_i^-$. Then, by \cite[Lemma 4(a)(i)]{agarwal2021new}, $\beta(t)$ increases to $x^*_i$ for all $t < \tau:=\inf\{t : \beta(t) = x^*_i\}$. If $t < \infty$, then $\beta(t) = x^*_i$ for all $t \geq \tau$ (as $x^*_i$ is an equilibrium point). Then, clearly, $\beta(t) \to x^*_i$ and $\ups(t) \to \mathbf{h}(x^*_i)$ as $t \to \infty$, as above.

Else say $\tau = \infty$;  then for every $\delta > 0$, there exists a $T_\delta < \infty$ (guaranteed as before by \cite[Lemma 4(a)(i)]{agarwal2021new} because by continuity the RHS of ODE can be uniformly bounded by non-zero values) such that:
\begin{align*}
x^*_i - \delta \leq \beta(t) \leq x^*_i + \delta \mbox{ for all }t \geq T_\delta.
\end{align*}
Thus, $\beta(t) \to x^*_i$ as $t \to \infty$. This also implies that:
\begin{align*}\underline{\mathbf{h}}_\delta(x^*_i) - \ups \leq \dot{\ups} \leq  \overline{\mathbf{h}}_\delta(x^*_i) - \ups \mbox{ for all }t \geq T_\delta, \mbox{ for } \overline{\mathbf{h}}_\delta(x^*_i) := \sup_{x \in \overline{{\cal N}_{\delta}}(x^*_i)}\mathbf{h}(x) \mbox{ and } \underline{\mathbf{h}}_\delta(x^*_i) := \inf_{x \in \overline{{\cal N}_{\delta}}(x^*_i)}\mathbf{h}(x).
\end{align*}
\ignore{Consider $\dot{z} = \overline{\mathbf{h}}_\delta(x^*_i) - z$ with initial condition $z(T_\delta) = \ups(T_\delta)$. Then, using classical methods of ODE, one can obtain the following solution of the ODE:
$$
z(t) = \overline{\mathbf{h}}_\delta(x^*_i) + e^{-t +T_\delta} + e^{-t + T_\delta}(z(T_\delta) - \overline{\mathbf{h}}_\delta(x^*_i)).
$$
Similarly, consider $\dot{p} = \underline{\mathbf{h}}_\delta(x^*_i) - p$ with initial condition $p(T_\delta) = \ups(T_\delta)$. Then, the solution is given by:
$$
p(t) = \underline{\mathbf{h}}_\delta(x^*_i) + e^{-t +T_\delta} + e^{-t + T_\delta}(p(T_\delta) - \underline{\mathbf{h}}_\delta(x^*_i)).
$$}
By Comparison Theorem in \cite{piccini1984ordinary} for ODEs having Lipschitz continuous right hand sides and using classical methods to derive the upper  and lower bounds, we get: 
$$
\underline{\mathbf{h}}_\delta(x^*_i)   + e^{-t + T_\delta}(
\ups(T_\delta) - \underline{\mathbf{h}}_\delta(x^*_i))
\leq \ups(t) \leq
\overline{\mathbf{h}}_\delta(x^*_i)  + e^{-t + T_\delta}( \ups(T_\delta) - \overline{\mathbf{h}}_\delta(x^*_i))
%
\mbox{ for all }t\geq T_\delta.
$$
Then clearly by considering limits $t\to \infty$ we have
\ignore{
This implies that $\limsup_{t\to \infty} \ups(t) \leq \limsup_{t\to \infty} z(t) = \lim_{t\to\infty} z(t) = \overline{\mathbf{h}}_\delta(x^*_i)$ and $\liminf_{t\to \infty} \ups(t) \geq \liminf_{t\to \infty} p(t) = \lim_{t\to\infty} p(t) = \underline{\mathbf{h}}_\delta(x^*_i)$. In all, we get that}:
\begin{align*}
    \underline{\mathbf{h}}_\delta(x^*_i) &\leq \liminf_{t\to \infty} \ups(t) \leq \limsup_{t\to \infty} \ups(t) \leq \overline{\mathbf{h}}_\delta(x^*_i), \mbox{ and now letting } \delta \to 0, 
    \ignore{\\
    \implies  \lim_{\delta \to 0} \underline{\mathbf{h}}_\delta(x^*_i) &\leq \liminf_{t\to \infty} \ups(t) \leq \limsup_{t\to \infty} \ups(t) \leq  \lim_{\delta \to 0} 
 \overline{\mathbf{h}}_\delta(x^*_i)\\
 \implies } \mathbf{h}(x^*_i) 
\leq \liminf_{t\to \infty} \ups(t) \leq \limsup_{t\to \infty} \ups(t) \leq  \mathbf{h}(x^*_i).
\end{align*}
Hence, $\ups(t) \to \mathbf{h}(x^*_i)$ as $t \to \infty$.

\underline{Part (ii)} If $\beta(0) = x^*_i$, then clearly $\beta(t) = x^*_i$ for all $t \geq 0$ and $\ups(t) \to \mathbf{h}(x^*_i)$ as $t \to \infty$. However if $\beta(0) \in {\cal N}_i^-$, then it can be shown as above that $\beta(t) \to y^* := \max\{y \in {\cal I}: y < x^*_i\}$. Similarly, if $\beta(0) \in {\cal N}_i^+$, then $\beta(t) \to y^* := \min\{y \in {\cal I}: y > x^*_i\}$. Thus, $x^*_i$ is a repeller for ODE \eqref{eqn_beta_ode_simple_prop} and $\mathbf{h}(x^*_i)$ is a saddle point for ODE \eqref{eqn_ODE_prop}.

\underline{Part (iii)} If $\beta(0) = x^*_i$, then clearly $\beta(t) = x^*_i$ for all $t \geq 0$ and $\ups(t) \to \mathbf{h}(x^*_i)$ as $t \to \infty$. Say $g(x) > 0$ for all $x \in {\cal N}_i^-$ and $g(x) > 0$ for all $x \in {\cal N}_i^+$.
Then, if $\beta(0) \in {\cal N}_i^-$, $\beta(t) \to x^*_i$, as shown for part 1. While if $\beta(0) \in {\cal N}_i^+$, then $\beta(t) \to y^* := \min\{y \in {\cal I}: y > x^*_i\}$, as shown for part 2. Thus, $x^*_i$ is a saddle point for ODE \eqref{eqn_beta_ode_simple_prop} and $\mathbf{h}(x^*_i)$ is a saddle point for ODE \eqref{eqn_ODE_prop}.

Lastly, consider the initial condition $\ups(0) \in \cD_I$ with $\pc(0) = 0$, then ODE \eqref{eqn_ODE_prop} simplifies to $\dot{\ups} = -\ups$, which clearly has  unique solution and $\ups(t) \to \mathbf{0}$ as $t \to \infty$. We have shown above that whenever $\pc(0) > 0$, $\ups(t) \nto \mathbf{0}$. Therefore, $\mathbf{0} \in \cR$. \eop
\end{pfthrmtwo}


\begin{pfthrmthree} At first, observe that in view of the hypothesis regarding ${\cal F}$ and \eqref{eqn_lambda} with $\sum_i \mu_i = 1$, the assumption \ref{a1_prop} holds. Further, it is clear from \eqref{eqn_tag_wi}-\eqref{eqn_final_shares}, \eqref{eqn_mean_matrix} and \eqref{eqn_mean_matrix_red_coloring} that the assumption \ref{a2_prop} holds. 

We will now prove that $\cA^u_\beta \neq \emptyset$, which will then imply that the assumption \ref{a3_prop} holds,  by Theorem \ref{thrm_beta_ODE_prop}; this would complete Theorem \ref{thm1}(i). Towards proving the claim, note that (recall $\mu_2 > 0$):
\begin{align}\label{eqn_bound_g}
\begin{aligned}
    \underline{g}^u_\beta(\beta)  &< g_\beta^u(\beta) \leq \overline{g}^u_\beta(\beta) \mbox{ for all } \beta \in [0,1] \mbox{ where} \\
    \underline{g}^u_\beta(\beta) &:= \left(-\beta \mu_2 - \beta \mu_1 (1-\alpha_x^u \rho) + (1-\beta) \mu_1 \rho  \alpha_y^u \right) m_f\eta^u - \beta \mu_a m_f \eta_a, 
    \mbox{ and}\\
    \overline{g}^u_\beta(\beta) &:= \left(- \beta \mu_1 (1-\alpha_x^u \rho) + (1-\beta) \mu_1 \rho  \alpha_y^u  + \mu_2 (1-\beta)\right) m_f\eta^u - \beta \mu_a m_f \eta_a.
\end{aligned}
\end{align}
Now, $\underline{g}^u_\beta(0) = \mu_1 \rho  \alpha_y^u m_f\eta^u \geq 0$; thus, $g_\beta^u(0) > 0$. Further, $\overline{g}^u_\beta(1) = -\mu_1(1-\alpha_x^u\rho) m_f\eta^u - \mu_a m_f \eta_a \leq 0$; thus, $g_\beta^u(1) \leq 0$. Since $g_\beta^u(\beta)$ is a continuous function of $\beta$, therefore there exists at least one zero of $g_\beta^u$, say $\beta^{u, \infty}$ such that $g_\beta^u > 0$ in ${\cal N}_\epsilon^-(\beta^{u, \infty})$ and  $g_\beta^u < 0$ in ${\cal N}_\epsilon^+(\beta^{u, \infty})$; ${\cal N}_\epsilon^+(\beta^{u, \infty}) = \emptyset$ if $\beta^{u, \infty} = 1$. Then, by Theorem \ref{thrm_beta_ODE_prop}, $\beta^{u, \infty} \in \cA_\beta^u$; thus, $\cA_\beta^u \neq \emptyset$.

Since $\overline{g}^u_\beta(\beta)$ is a linear function such that $\overline{g}^u_\beta(0) > 0$ and (recall) $\overline{g}^u_\beta(1) \leq 0$, therefore, $\overline{\beta}^u \in (0,1]$, given in \eqref{eqn_beta_bar}, is the unique zero of $\overline{g}^u_\beta(\beta)$. Further, since $g_\beta^u \leq \overline{g}^u_\beta$ and $\overline{g}^u_\beta(\beta) < 0$ for all $\beta \in (\overline{\beta}^u, 1]$ when $\overline{\beta}^u < 1$, therefore, there exists no zero of $g_\beta^u$ in $(\overline{\beta}^u, 1]$; if $\overline{\beta}^u = 1$, then also, any zero of $g_\beta^u$ is atmost $1$. Thus, if at all, there is any zero of $g_\beta^u$, which can be an attractor or repeller or saddle point of \eqref{eqn_general_g_beta}, it is lesser than or equals to $\overline{\beta}^u$. Next, notice that there is a unique zero of the function $\underline{g}^u_\beta$, namely $\underline{\beta}^u \in (0,1)$, as given in \eqref{eqn_beta_bar}. Again using similar arguments as before, we get that $\beta^{u, \infty} > \underline{\beta}^u$. This proves \eqref{eqn_beta_bar}. 

Now, by Theorem \ref{thrm_beta_ODE_prop}, the attractor and saddle sets are as in the hypothesis with subset of the combined domain of attraction as $\cD_I$. 

We will now identify the compact sub-domain of $\cD_I$ for completing the proof using Theorem \ref{thm1}. From \ref{a1_prop} for our case, one can bound $\Pa_n$:
\begin{align*}
    0 \leq \Pa_n \leq \overline{\Psi}_n^a :=  \frac{1}{n}\left(\sum_{k=1}^{\min\{\nu_e, n\}} 2{\cal F}1_{\{\Pc_k > 0\}} + s_0^c \right).
\end{align*}
By strong law of large numbers, $\overline{\Psi}_n^a \to 2E[{\cal F}]$ a.s. in survival paths and $\overline{\Psi}_n^a \to 0$ in extinction paths, as $n \to \infty$. Thus, $\cS :=  \cD_I \cap  \left\{\ups :  \pa \in [0, 2E({\cal F}) ]\right\}$ is the compact subset of $\cD_I$ and $p_{b} := P(\Ups_n \mbox{ visits } \cS \mbox{ i.o.}) = 1$. Hence, by Theorem \ref{thm1}(ii), the claim holds. \eop
\end{pfthrmthree}


\begin{pfthrmfour}
Let all parameters except $\kappa$ be fixed. 
Consider the case when $\nabla^u(\kappa, \kappa + \partial \kappa) = g^u_\beta(\beta^{\infty, u}(\kappa); \kappa+\partial \kappa) > 0$ for some $\partial \kappa > 0$. Since $g^u_\beta(\beta; \kappa+\partial \kappa)$ is either a convex or concave or linear function of $\beta$ with a unique zero in $(0,1)$, therefore, there exists a $\beta^{\infty, u}(\kappa + \partial \kappa) > \beta^{\infty, u}(\kappa)$ such that $g^u_\beta(\beta^{\infty, u}(\kappa + \partial \kappa); \kappa+\partial \kappa) = 0$. One can prove the claim similarly when $\nabla^u(\kappa, \kappa + \partial \kappa) < 0$. Lastly if $\nabla^u(\kappa, \kappa + \partial \kappa) = 0$, then again due to uniqueness, $\beta^{\infty, u}(\kappa + \partial \kappa) = \beta^{\infty, u}(\kappa)$. \eop
\end{pfthrmfour}


\begin{pfcorollaryexWM}
We will first show that the function $g_\beta^{o, u}$ is either convex or concave or linear depending upon warning-specific and user-specific parameters. Towards this, note that for each $u$:
\begin{align}\label{eqn_convex_g_beta}
    \begin{aligned}
    \frac{d g^{o, u}_\beta(\beta)}{d \beta} 
    &= - (\mu_1 + \mu_2)m_f\eta^u + (\alpha_x^u - \alpha_y^u) (\mu_1 \rho + \mu_2 \omega(\beta)) m_f\eta^u + (\beta \alpha_x^u  + (1-\beta)\alpha_y^u)   \frac{bw \mu_2 m_f\eta^u}{(\beta + b (1-\beta))^2} - \mu_a m_f \eta_a\\
    \implies    \frac{d^2 g^{o, u}_\beta(\beta)}{d\beta^2}
    &= \frac{2 m_f\eta^u bw \mu_2}{(\beta + b (1-\beta))^3} \left( b \alpha_x^u - \alpha_y^u \right).
    \end{aligned}
\end{align} 
Thus, if $bw\mu_2(b \alpha_x^u - \alpha_y^u)= 0$ or $< 0$ or $>0$, then $g_\beta^{o, u}$ is a linear, concave or convex function respectively. From \eqref{eqn_beta_ODE_etac1}:
\begin{align*}
    g^{o, u}_\beta(0) &= 
    \left(\mu_1 \rho + \mu_2 \gamma\right)\alpha_y^u m_f\eta^u > 0, \mbox{ and }
    g^{o, u}_\beta(1) 
    = - \bigg( \mu_1 m_f\eta^u \left(1-\alpha_x^u \rho\right) + \mu_2 (1 - \alpha_x^u(w+\gamma))m_f\eta^u + \mu_a m_f \eta_a \bigg) < 0;
\end{align*}the last inequality in above holds as $\alpha_x^u(w+\gamma) \leq 1$ for each $u$ and $\alpha_x^u \rho < \alpha_x^u < 1$. 
Therefore, there exists a unique $\beta^{o, \infty, u} \in (0,1)$ such that $g^u_\beta(\beta^{o, \infty, u}) = 0$, $g^u_\beta(\beta) > 0$ for all $\beta \in [0, \beta^{o, \infty, u})$ and $g^u_\beta(\beta) < 0$ for all $\beta \in (\beta^{o, \infty, u}, 1]$. This implies that for the ODE \eqref{eqn_general_g_beta}, $t \mapsto \beta^{ u}(t)$ is strictly increasing if $\beta^{u}(0) \in [0, \beta^{o, \infty, u})$ and strictly decreasing if $\beta^{u}(0) \in (\beta^{o, \infty, u}, 1]$. Thus, $\cA_\beta^{o,u} = \{\beta^{o, \infty, u}\}$ with  the domain of attraction as $[0,1]$. Lastly, observe that $g_\beta^{o, u}(\beta) \leq \overline{g}_\beta^u(\beta)$ for each $\beta \in [0,1]$, therefore, $\beta^{o, \infty, u} \leq \overline{\beta}^u$, as these two zeroes are unique zeroes of their respective functions (see \eqref{eqn_bound_g}). 
\end{pfcorollaryexWM}


\begin{pflimitswarning}
Recall from Corollary \ref{corollary_ex_wm}, $g_\beta^{o, u}$ has a unique attractor, $\beta^{o, \infty, u} \in (0,1)$, for each $u\in\{R, F\}$. Observe further that $g^{o, u}_\beta(\beta^{o, \infty, u}(w); w) = 0$ and $g^u_\beta(\beta^{o, \infty, u}(b); b)  = 0$. Henceforth, the corollary will be proved using Theorem \ref{thrm_unique_att}. For any $\partial w > 0$ and $\partial b > 0$, we get:
\begin{align*}
        \nabla^u(w, w+\partial w) &= g^{o, u}_\beta(\beta^{o, \infty, u}(w); w+\partial w) \\
        &= g^{o, u}_\beta(\beta^{o, \infty, u}(w); w) + m_f\eta^u \mu_2 \bigg(\alpha_x^u \beta^{o, \infty, u}(w) + \alpha_y^u (1-\beta^{o, \infty, u}(w))\bigg) \left(\frac{\partial w \beta^{o, \infty, u}(w) }{\beta^{o, \infty, u}(w) + (1-\beta^{o, \infty, u}(w))b}\right) > 0 \mbox{ and}\\
        \nabla^u(b, b+\partial b) &= g^{o, u}_\beta(\beta^{o, \infty, u}(b); b+\partial b) \\
        &= g^{o, u}_\beta(\beta^{o, \infty, u}(b); b) - \partial b m_f\eta^u \mu_2 \frac{w \beta^{o, \infty, u}(b) \bigg(\alpha_x^u \beta^{o, \infty, u}(b) + \alpha_y^u (1-\beta^{o, \infty, u}(b)\bigg)}{\bigg(\beta^{o, \infty, u}(b) + (1-\beta^{o, \infty, u}(b))(b+\partial b)\bigg)\bigg(\beta^{o, \infty, u}(b) + (1-\beta^{o, \infty, u}(b))b\bigg)} < 0.
\end{align*}
Thus, by Theorem \ref{thrm_unique_att}, $\beta^{o, \infty, u}(w, b)$ strictly increases with $w$ and strictly decreases with $b$ for any $u \in \{R, F\}$. \eop

\end{pflimitswarning}

\begin{pfthrmfive} 
In this proof, we explicitly show the dependency of zeros of \eqref{eqn_beta_ODE_etac1} on design parameters $(w, b)$.

\noindent \underline{Part (i)} Consider a $\delta > 0$ such that $\beta^{o, \infty, R}(\overline{w}, 0) > \delta$. Then, $w \in [0, \overline{w}] = W_1 \cup W_2$, where $W_1 := \{w : \beta^{o, \infty, R}(w, 0) > \delta \}$ and $W_2 := \{w : \beta^{o, \infty, R}(w, 0) \leq \delta \}$. If $W_2 \neq \emptyset$, by Corollary \ref{corollary_ex_wm}, there exists a $\widetilde{w} > 0$ such that $\beta^{o, \infty, R}(\widetilde{w}, 0) = \delta$,  $W_1 = \{w : w > \widetilde{w} \}$, and $W_2 := \{w : w \leq \widetilde{w} \}$. The proof for case with $W_2 = \emptyset$ is trivially true once the other case is proved. Hence, consider $W_2 \neq \emptyset$.

Consider $w \in W_1$.  Then, by Corollary \ref{cor_limits_warning}, there exists a unique $b(w; \delta) > 0$ such that $\beta^{o, \infty, R}(w, b(w; \delta)) = \delta$ (i.e., the zero of $g_\beta^{o, F}$ equals $\delta$) and hence:
\begin{align}\label{eqn_bw}
     b(w; \delta) := \left(\frac{\delta}{1-\delta}\right)\left(w p(\delta) - 1 \right), \mbox{ where } p(\delta) := \frac{\eta^R \mu_2(\delta \alpha_x^R + (1-\delta)\alpha_y^R)}{\delta ((\mu_1+\mu_2)\eta^R + \mu_a  \eta_a) - \eta^R (\mu_1 \rho + \mu_2 \gamma) (\delta \alpha_x^R + (1-\delta) \alpha_y^R)}.
\end{align}Thus, again by Corollary \ref{cor_limits_warning} and because $[0, \overline{w}] \cap W_1 = (\widetilde{w}, \overline{w}]$ (as said before):
    \begin{align}\label{eqn_opt_prob_reduced}
    \begin{aligned}
        \sup_{w \in [0, \overline{w}] \cap W_1; b \in [0, \infty); \beta^{o, \infty, R}(w, b) \leq \delta}\beta^{o, \infty, F}(w, b) = \sup_{w \in (\widetilde{w}, \overline{w}]}\beta^{o, \infty, F} (w, b(w; \delta)).
    \end{aligned}
    \end{align}

By Lemma \ref{cor_bw}, $\beta^{o, \infty, F}(w, b(w; \delta))$ strictly increases with $w$, for every $\delta > 0$. Then, the optimal value for the problem in \eqref{eqn_opt_prob_reduced} is given by:
\begin{align}\label{eqn_opt_prob_reducedW1}
    \begin{aligned}
 \sup_{w \in (\widetilde{w}, \overline{w}]}\beta^{o, \infty, F} (w, b(w; \delta)) = \beta^{o, \infty, F}(\overline{w}, b(\overline{w}; \delta)).
\end{aligned}
    \end{align}

Now, consider $w \in W_2$. Then, $\beta^{o, \infty, R}(w, 0) \leq \delta$. Further
by Corollary \ref{cor_limits_warning}, for any $w < \widetilde{w}$ and $b > 0$, we have:
\begin{equation*}
     \beta^{o, \infty, F}(\widetilde{w}, 0) > \beta^{o, \infty, F}({w}, 0) > \beta^{o, \infty, F}(w, b),  \mbox{ and } \beta^{o, \infty, F}(\widetilde{w}, 0) > \beta^{o, \infty, F}(\widetilde{w}, b).
\end{equation*}
Thus, we have:
\begin{align}\label{eqn_opt_prob_reduced_W2}
    \begin{aligned}
        \sup_{w \in [0, \overline{w}] \cap W_2; b \in [0, \infty); \beta^{o, \infty, R}(w, b) \leq \delta}\beta^{o, \infty, F}(w, b) = \beta^{o, \infty, F} (\widetilde{w}, 0).
    \end{aligned}
    \end{align}

In all, by \eqref{eqn_opt_prob_reducedW1}, \eqref{eqn_opt_prob_reduced_W2}, we have:
\begin{align}\label{eqn_opt_W1_W2}
    \sup_{w \in [0, \overline{w}]; b \in [0, \infty); \beta^{o, \infty, R}(w, b) \leq \delta}\beta^{o, \infty, F}(w, b)  =  \max\bigg\{\beta^{o, \infty, F}(\overline{w}, b(\overline{w}; \delta)), \beta^{o, \infty, F}(\widetilde{w}, 0)\bigg\}.
 \end{align}
 
Let us now consider a sequence of $w \down \widetilde{w}$ and observe $\frac{\partial b(w; \delta)}{\partial w} = \left(b(w; \delta) + \frac{\delta}{1-\delta}\right)\frac{1}{w} > 0$. Thus, $b(w; \delta)$ decreases as $w$ decreases. We claim that $\lim_{w \down \widetilde{w}} b(w; \delta) = 0$. Let us suppose on the contrary that the limit is positive; note that the limit can not be negative as $b(w; \delta) > 0$. By continuity of $b(w; \delta)$ with respect to $w$ (see \eqref{eqn_bw}), there exists a $w' < \widetilde{w}$ such that $b(w'; \delta) > 0$, and further $\beta^{o, \infty, R}(w', b(w'; \delta)) = \delta$, by definition of $b(w'; \delta)$. However, since $w' \in W_2$, we also have $\beta^{o, \infty, R}(w', 0) \leq \delta$, leading to a contradiction. Thus, the limit is $0$. 

Consider function $L(\beta; w) := \bigg(g_\beta^{o, F} (\beta(w, b(w))\bigg)^2$. Clearly this function is jointly  continuous and has a unique minimum at $\beta^{o, \infty, F}(w, b(w))$ for each $w$  (as it is the unique zero of $g_\beta^o (\cdot)$). Hence  
by Maximum Theorem: 
$$
\beta^{o, \infty, F}(w, b(w)) \to \beta^{o, \infty, F}(\widetilde{w}, 0), \mbox{ as } w \down \widetilde{w}  \  \mbox{ and further by Lemma \ref{cor_bw}, } \  \beta^{o, \infty, F}(w, b(w)) \down \beta^{o, \infty, F}(\widetilde{w}, 0).
$$
Thus, $\beta^{o, \infty, F}(\widetilde{w}, 0) \leq \beta^{o, \infty, F}(w, b(w; \delta)) < \beta^{o, \infty, F}(\overline{w}, b(\overline{w}; \delta))$, where the last inequality is again due to Lemma \ref{cor_bw}. Conclusively, by \eqref{eqn_opt_W1_W2}, we get that 
$$
\sup_{w \in [0, \overline{w}]; b \in [0, \infty); \beta^{o, \infty, R}(w, b) \leq \delta}\beta^{o, \infty, F}(w, b) = \beta^{o, \infty, F}(\overline{w}, b(\overline{w}; \delta)).
$$

\noindent \underline{Part (ii)} Consider $\delta > 0$ such that $\beta^{o, \infty, R}(\overline{w}, 0) \leq \delta$.
Again, by Corollary \ref{cor_limits_warning}, for all $w \in [0, \overline{w}]$ and $b > 0$:
\begin{equation*}
     \beta^{o, \infty, F}(\overline{w}, 0) > \beta^{o, \infty, F}({w}, 0) > \beta^{o, \infty, F}(w, b),  \mbox{ and } \beta^{o, \infty, F}(\overline{w}, 0) > \beta^{o, \infty, F}(\overline{w}, b).
\end{equation*}Thus, the optimal value is achieved at $b = 0$ and $w = \overline{w}$, with $\beta^{o, \infty, R}(\overline{w}, 0) \leq \delta$. \eop
\end{pfthrmfive}

\begin{lem}\label{cor_bw}
    The function $\beta^{o, \infty, F}(w, b(w; \delta))$ strictly increases with $w$, when $w < \overline{w}$, for every $\delta > 0$.
\end{lem}
\begin{proof}
     Fix $w$ and $\partial w > 0$, we have (for simplicity, denote $\beta^{o, \infty, F}(w, b(w; \delta))$ by $\beta_\delta(w)$):

    \vspace{-4mm}
    {\small 
    \begin{align}\label{eqn_nabla_betaF}
    \begin{aligned}
        \nabla^F(w, w+\partial w; b) &= g^{o, F}_\beta(\beta_\delta(w); w+\partial w) - g^{o, F}_\beta(\beta_\delta(w); w)  \\
        &\hspace{-2cm}= m_f\eta^F \mu_2 \bigg(\alpha_x^F \beta_\delta(w) + \alpha_y^F (1-\beta_\delta(w))\bigg) \left( \left(\frac{(w + \partial w) \beta_\delta(w) }{\beta_\delta(w) + (1-\beta_\delta(w))b(w + \partial w; \delta)}\right) - \left(\frac{w \beta_\delta(w) }{\beta_\delta(w) + (1-\beta_\delta(w))b(w; \delta)}\right) \right)\\
        & \hspace{-2cm}= \frac{m_f\eta^F \mu_2 \beta_\delta(w)  \bigg(\alpha_x^F \beta_\delta(w) + \alpha_y^F (1-\beta_\delta(w))\bigg) }{\bigg(\beta_\delta(w) + (1-\beta_\delta(w))b(w + \partial w; \delta)\bigg) \bigg( \beta_\delta(w) + (1-\beta_\delta(w))b(w; \delta) \bigg)} \bigg( \partial w\bigg( (1-\beta_\delta(w)) b(\partial w; \delta) + \beta_\delta(w) \bigg)\bigg), \mbox{ where}
        \end{aligned}
        \end{align}}the last equality follows by simple algebra after substituting for $b(\cdot; \delta)$ from \eqref{eqn_bw}. Since $b(\cdot, \delta) > 0$ and  by Theorem \ref{thrm_BP_to_fake}, $\beta_\delta(w) \in (\underline{\beta}^F, \overline{\beta}^F] \subset [0,1]$, therefore,  $\nabla^F(w, w+\partial w; \delta) > 0$ for any $\delta > 0$. Thus, the proof follows by Theorem \ref{thrm_unique_att}. 
\end{proof}

\begin{pfcorbetaona}
Consider any $\delta > 0$, $\mu_a \in (0, 1-\mu_1-\mu_2]$ and let $b^*, w^*$ be as in Theorem \ref{thrm_opt}.

\noindent \underline{Case 1: when $b^* > 0$:} Let $b^*(\mu_a 
= 0) =: b^*_0$. From \eqref{eqn_optimal_parameter_old}, observe that $b^*$ is a strictly decreasing function of $\mu_a$, therefore, $b^*_0 > b^* > 0$. Further, from \eqref{eqn_beta_ODE_etac1}, we have:
\begin{align}\label{eqn_g_beta_na}
\begin{aligned}
    g_\beta^{o, F}(\betana) &= g_\beta^{o, F}(\betana; \mu_a = 0) - \betana \mu_a m_f \eta_a \\
    &\hspace{1cm}+ \mu_2 m_f \eta^F \bigg(\betana \alpha_x^F + (1-\betana) \alpha_y^F \bigg) \left( \frac{w^*\betana}{\betana + (1-\betana) b^*} - \frac{w^*\betana}{\betana + (1-\betana) b^*_0} \right) \\
    &= 0 + \betana \mu_a m_f \eta_a \left( \frac{\mu_2 \eta^F w^*(1-\betana) (\betana \alpha_x^F + (1-\betana) \alpha_y^F)}{\bigg(\betana + (1-\betana) b^*_0 \bigg) \bigg( \betana + (1-\betana) b^* \bigg) }  \left( \frac{b^*_0 - b^*}{\eta_a \mu_a}\right) - 1 \right).
\end{aligned}
\end{align}
Define $p(\mu_a) := \delta ((\mu_1+\mu_2)\eta^R + \mu_a  \eta_a) - \eta^R (\mu_1 \rho + \mu_2 \gamma) (\delta \alpha_x^R + (1-\delta) \alpha_y^R)$. Then, by \eqref{eqn_optimal_parameter_old}, we have: 
\begin{align*}
\frac{b^*_0 - b^*}{\eta_a \mu_a} = \left(\frac{\delta^2}{1-\delta}\right) \left( \frac{w^* \eta^R \mu_2 (\delta \alpha_x^R + (1-\delta)\alpha_y^R)}{p(\mu_a) p(0)} \right).
\end{align*}
Substitute the above term in  \eqref{eqn_g_beta_na} and consider the following limit to analyse \eqref{eqn_g_beta_na}:
\begin{align*}
\begin{aligned}
    \lim_{\delta \to 0} \left( \frac{\mu_2 \eta^F w^*(1-\betana) (\betana \alpha_x^F + (1-\betana) \alpha_y^F)}{\bigg(\betana + (1-\betana) b^*_0 \bigg) \bigg( \betana + (1-\betana) b^* \bigg) } \right) \lim_{\delta \to 0} \left(\frac{\delta^2}{1-\delta}\right) \lim_{\delta \to 0} \left( \frac{w^* \eta^R \mu_2 (\delta \alpha_x^R + (1-\delta)\alpha_y^R)}{p(\mu_a) p(0)} \right) - 1.
\end{aligned}
\end{align*}In the above, the second limit is clearly $0$ and the rate of convergence is independent of other factors. The first and third limits are finite, and the respective terms can be upper-bounded independent of $\mu_a$ and other factors. Thus, the product of three limits is $0$, and the rate of convergence is uniform in $\mu_a$ and $b^*$, i.e., there exists a $\overline{\delta} > 0$ such that (for example):
$$
 \left( \frac{\mu_2 \eta^F w^*(1-\betana) (\betana \alpha_x^F + (1-\betana) \alpha_y^F)}{\bigg(\betana + (1-\betana) b^*_0 \bigg) \bigg( \betana + (1-\betana) b^* \bigg) } \right)  \left(\frac{\delta^2}{1-\delta}\right)  \left( \frac{w^* \eta^R \mu_2 (\delta \alpha_x^R + (1-\delta)\alpha_y^R)}{p(\mu_a) p(0)} \right) - 1 < -\frac{1}{2} \mbox{ for all } \delta \leq \overline{\delta} \mbox{ and } \mu_a > 0.
$$
Thus, from \eqref{eqn_g_beta_na}, $g_\beta^{o, F}(\betana) < - \betana \mu_a m_f \eta_a/2 < 0$ for any $\mu_a > 0$ and all $\delta \leq \overline{\delta}$.

Recall from the proof of Corollary \ref{cor_limits_warning} that $g_\beta^{o, F}(\cdot)$ is either convex/concave/linear with a unique zero in $(0,1)$. Therefore, the unique zero of $g_\beta^{o, F}(\cdot; \mu_a)$, namely $\beta^o(\mu_a) < \betana$ for all $\delta \leq \overline{\delta}$ and  for all $\mu_a \in (0, 1-\mu_1-\mu_2]$.

\noindent \underline{Case 2: when $b^* = 0$:} Here, again $b^*(\mu_a 
= 0) =: b^*_0 > b^* = 0$.  Then, similar to \eqref{eqn_g_beta_na}, using \eqref{eqn_optimal_parameter_old}:
\begin{align*}
    g_\beta^{o, F}(\betana) &= \betana \mu_a m_f \eta_a \left( \frac{\mu_2 \eta^F w^* (1-\betana) (\betana \alpha_x^F + (1-\betana) \alpha_y^F) b^*_0 }{\eta_a \mu_a  \betana \bigg(\betana + (1-\betana) b^*_0 \bigg)   } \left(\left(\frac{\delta}{1-\delta}\right) \left( \frac{w^* \eta^R \mu_2 (\delta \alpha_x^R + (1-\delta)\alpha_y^R)}{p(0)}  - 1 \right) \right) - 1 \right).
\end{align*}
Hereafter, the proof follows as in Case 1. \eop
\end{pfcorbetaona}

\begin{pfcorollaryeaWM} 
We begin the proof for the fake-post.

\noindent \underline{Part (i)} Consider
$0 < \mu_a \leq \min\{1-\mu_1-\mu_2, \Delta_a\}$. Then, by the definition of upper-bound $\Delta_a$ and \eqref{eqn_warning_ea}, $\alpha_x^F \omega^a(\betana) \leq 1$. Note from \eqref{eqn_warning_ea} that $\omega^a(\beta)$ is a strictly increasing function of $\beta$. Therefore, $\alpha_x^F \omega^a(\beta) \leq 1$ for all $\beta \leq \betana$, for given $\mu_a$.

This implies that for $\beta \le \betana$, we have $g_\beta^{a, F}(\beta) = g_\beta^{o, F}(\beta; \mu_a = 0)$ (see \eqref{eqn_g_F_eaWM}). Further, $\betana$ is a zero of $g_\beta^{a, F}$, as $g_\beta^{a, F}(\betana) = g_\beta^{o, F}(\betana; \mu_a = 0) = 0$. Furthermore, by uniqueness given in Corollary \ref{corollary_ex_wm}, $\betana$ is the unique zero of $g_\beta^{a, F}$ in $[0, \betana]$. Therefore, any $\beta^{a} \in {\cA}^{a, F}_\beta \cup \cR^{a, F}_\beta$ is  in $[\betana, 1]$.

\noindent \underline{Part (ii)} Consider $\mu_a > \Delta_a$. Then, the corresponding $\alpha_x^F \omega^a(\betana) > 1$.
Define the function $h(\beta) := \alpha_x^F \omega^a({\beta}) - 1$. It is easy to see that $h(0) < 0$, $h(1) > 0$ and $h(\cdot)$ is a strictly increasing function. Thus, there exists a unique zero of $h$, denoted by  $\widetilde{\beta} \in (0,1)$, i.e., $\alpha_x^F \omega^a(\widetilde{\beta}) = 1$. As $\beta \mapsto \omega^a (\beta)$ is strictly increasing, we further have $\alpha_x^F \omega^a(\beta) < 1$ for all $\beta < \widetilde{\beta} $; furthermore $\widetilde{\beta} < \betana$ as $\alpha_x^F \omega^a(\betana) > 1$. 

From \eqref{eqn_g_F_eaWM}, we have:
\begin{align}\label{eqn_cor3}
\begin{aligned}
    g^{a, F}_\beta(\beta) &= g^{o, F}_\beta(\beta; \mu_a = 0) + \mu_2 m_f \eta^F \left\{ \beta \bigg(\min\{1, \alpha_x^F \omega^a(\beta)\} - \alpha_x^F \omega^a(\beta) \bigg) + (1-\beta) \bigg( \min\{1, \alpha_y^F \omega^a(\beta)\} - \alpha_y^F\omega^a(\beta) \bigg) \right\}. 
\end{aligned}
\end{align}Thus, $g^{a, F}_\beta(\beta) < g^{o, F}_\beta(\beta; \mu_a = 0)$ if $1 < \alpha_{j}^F \omega^a(\beta)$ for some $j \in \{x, y\}$, and $g^{a, F}_\beta(\beta) = g^{o, F}_\beta(\beta; \mu_a = 0)$ if $\alpha_{j}^F \omega^a(\beta) \leq 1$ for each $j \in \{x, y\}$. As a result, we have:

(a) for $\beta \in [0, \widetilde{\beta}]$, $g_\beta^{a, F}(\beta) = g_\beta^{o, F}(\beta; \mu_a = 0) > 0$, and 

(b) for $\beta \in [\betana, 1]$, $g_\beta^{a, F}(\beta) < g_\beta^{o, F}(\beta; \mu_a = 0) \leq 0$.

By Theorem \ref{thrm_BP_to_fake},  there exists at least one zero of $g_\beta^{a, F}$, say $\beta^a $ and by above arguments, $\beta^a \in (\widetilde{\beta}, \betana)$. We will now claim and show that $\beta^a > \beta^o$, but first observe that $\beta^o < \betana$ by Corollary \ref{cor_beta_o_na}. Towards this, note that for $\beta \in (\widetilde{\beta}, \betana)$, we have:
\begin{align}
\begin{aligned}
    g^{a, F}_\beta(\beta) &= g^{o, F}_\beta(\beta) + \mu_2 m_f \eta^F \left\{ \beta \bigg(\min\{1, \alpha_x^F \omega^a(\beta)\} - \alpha_x^F \omega(\beta) \bigg) + (1-\beta) \bigg( \min\{1, \alpha_y^F \omega^a(\beta)\} - \alpha_y^F\omega(\beta) \bigg) \right\}\\
    &= g^{o, F}_\beta(\beta) + \mu_2 m_f \eta^F \left\{ \beta \bigg(1 - \alpha_x^F \omega(\beta) \bigg) + (1-\beta) \bigg( \min\{1, \alpha_y^F \omega^a(\beta)\} - \alpha_y^F\omega(\beta) \bigg) \right\}.
\end{aligned}
\end{align}
In the above, if $1 > \alpha_y^F \omega^a(\beta)$, then:
\begin{align*}
    g^{a, F}_\beta(\beta) &= g^{o, F}_\beta(\beta) + \mu_2 m_f \eta^F \left\{ \beta \bigg(1 - \alpha_x^F \omega(\beta) \bigg) + (1-\beta)\alpha_y^F  \bigg( \omega^a(\beta) - \omega(\beta) \bigg) \right\}\\
    &= g^{o, F}_\beta(\beta) + \mu_2 m_f \eta^F \left\{ \beta \bigg(1 - \alpha_x^F \omega(\beta) \bigg) + (1-\beta)\alpha_y^F  \left( \frac{\beta \mu_a m_f \eta_a}{ \mu_2 m_f\eta^F \left(\beta \alpha_x^F + (1-\beta)\alpha_y^F\right) } \right) \right\}  > g^{o, F}_\beta(\beta), 
\end{align*}
as $ \omega(\beta)\alpha_x^F < 1$ for all $\beta \in [0,1)$. Further, $ \omega(\beta)\alpha_y^F < 1$ for all $\beta \in [0,1)$, hence even with $1 \leq  \alpha_y^F \omega^a(\beta)$, we have:
\begin{align*}
    g^{a, F}_\beta(\beta) &= g^{o, F}_\beta(\beta) + \mu_2 m_f \eta^F \left\{ \beta \bigg(1 - \alpha_x^F \omega(\beta) \bigg) + (1-\beta)  \bigg( 1- \alpha_y^F\omega(\beta) \bigg) \right\} > g^{o, F}_\beta(\beta).
\end{align*}
Now, for $\beta \in (\widetilde{\beta}, \beta^o]$, $g_\beta^{o, F}(\beta) \geq 0$, and thus, $g^{a, F}_\beta(\beta) > 0$.  
This completes the proof of the claim.

Now, consider the real-post. By Theorem \ref{thrm_BP_to_fake}, $\cA^{a, R}_\beta \neq \emptyset$, therefore, there exists at least one zero of $g_\beta^{a, R}$, say $\beta^{a, R} \in (0,1)$. Now, using arguments as above:
\begin{align*}
    g^{a, R}_\beta(\beta) &= g^{o, R}_\beta(\beta; \mu_a = 0) + \mu_2 m_f \eta^R \left\{ \beta \bigg(\min\{1, \alpha_x^R \omega^a(\beta)\} - \alpha_x^R \omega^a(\beta) \bigg) + (1-\beta) \bigg( \min\{1, \alpha_y^R \omega^a(\beta)\} - \alpha_y^R \omega^a(\beta) \bigg) \right\} \\
    &\hspace{2cm} + \beta \mu_a m_f \eta_a  \left( \frac{\eta^R}{\eta^F}\left( \frac{\beta \alpha_x^R + (1-\beta) \alpha_y^R}{\beta \alpha_x^F + (1-\beta) \alpha_y^F}\right) - 1  \right) < g^{o, R}_\beta(\beta; \mu_a = 0).
\end{align*}
Thus, any zero of $g^{a, R}_\beta(\beta)$ is strictly less than the unique zero of $g^{o, R}_\beta(\beta; \mu_a = 0)$, i.e., $\beta^{a, R} < \beta^{o, R}(\mu_a = 0) \leq \delta$, for any $\beta^{a, R} \in \cA^{a, R}_\beta \cup \cR^{a, R}_\beta$ (see Theorem \ref{thrm_opt}). \eop
\end{pfcorollaryeaWM}

\begin{pfcorollaryehWM}

We divide the proof in two cases.

\textbf{Case 1:} If $\overline{\phi} < \frac{1}{\alpha_y^R \omega^a(\delta)}$. Then $\overline{\phi}$ is the unique zero of $g_{\beta, {\phi}}^{h, R}(\delta) = 0$.  Further, for any $\phi' \in \left( \overline{\phi}, \frac{1}{\alpha_y^R \omega^a(\delta)}\right)$, $g_{\beta, \phi'}^{h, R}(\delta) > 0$. By \eqref{eqn_bound_g}, $g_{\beta, \phi'}^{h, R}(1) < 0$.  Thus, there exists at least one zero of $g_{\beta, \phi'}^{h, R}$ greater than $\delta$. Now, consider any $\phi' \geq \frac{1}{\alpha_y^R \omega^a(\delta)}$. Since the function $\phi \mapsto g_{\beta, {\phi}}^{h, R}(\delta)$ is continuous, therefore, $g_{\beta, {\phi'}}^{h, R}(\delta) > 0$ for such $\phi$. Thus, again as before, there exists at least one zero of $g_{\beta, \phi'}^{h, R}$ greater than $\delta$. Hence, any $\phi$ satisfying the constraint in \eqref{eqn_opt_phi} is  less than or equals to $\overline{\phi}$. Thus, the optimizer of \eqref{eqn_opt_phi} is $\phi^* = \overline{\phi}$.

\textbf{Case 2:} If $\overline{\phi} \geq \frac{1}{\alpha_y^R \omega^a(\delta)}$. Then by monotonicity, for any $\phi \ge \overline{\phi} $:
$$
g_{\beta, \phi}^{h, R} (\beta) \leq q_\phi(\beta) := \left(-\beta \mu_2 - \beta \mu_1 (1-\alpha_x^R \rho) + (1-\beta) \mu_1 \rho  \alpha_y^R + \mu_2 \phi \omega^a(\beta) \bigg(\beta  \alpha_x^R + (1-\beta)\alpha_y^R\bigg) \right) m_f\eta^R - \beta \mu_a m_f \eta_a.
$$
Thus for all such $\phi$, $q_\phi(\delta)$ is a strictly increasing function  of $\phi$  with  $q_1(\delta) < 0$ (by Theorem \ref{corollary_ea_wm}) and $q_{\overline{\phi}}(\delta) = 0$. Thus, $g_{\beta, \phi}^{h, R} (\delta) \le q_\phi(\delta) \leq 0$. 

Further, by strict monotonicity  of $\omega^a(\cdot)$ in $\beta$, 
we have $\phi \alpha_y^R \omega^a(\beta) \geq 1$ for all $\beta > \delta$ whenever   $\phi \geq \overline{\phi}$. Thus, $g_{\beta, \phi}^{h, R} (\beta) $ is linearly (strictly) decreasing in $\beta$, when $\beta > \delta$. As already proved  $g_{\beta, \phi}^{h, R} (\delta) < 0$, and hence  $g_{\beta, \phi}^{h, R} (\beta) < 0$ for all $\beta > \delta$. Hence, the feasibility condition of \eqref{eqn_opt_phi} is satisfied for any $\phi \ge \overline{\phi}$. 

By definition of $\phi^*$ in this case (the second row), we have:
$$
    \phi^* \omega^a(\beta) \alpha_y^F = 1 \mbox{ for all } \beta \geq \underline{\beta}^F.
$$
Further,  $\min\{1,  \phi^* \omega^a(\beta) \alpha_y^F\} = 1 $  for all $\beta > \underline{\beta}^F$,  when  $\phi \geq \phi^*$. Thus, the functions $g^{h, F}_{\beta, \phi}(\beta) = g^{h, F}_{\beta, \phi^*}(\beta)$ for all $\beta \geq \underline{\beta}^F$. Also, by Theorem \ref{thrm_BP_to_fake}, any  zero of $g^{h, F}_{\beta, \phi}$ is larger than $\underline{\beta}^F$. Thus, $\left\{\beta: \beta \in \cA^{h, \phi}_\beta \cup \cR^{h, \phi}_\beta\right\} = \left\{\beta: \beta \in \cA^{h, \phi^*}_\beta \cup \cR^{h, \phi^*}_\beta\right\}$. Now, given any $\beta$, observe that $\phi \mapsto g_{\beta, \phi}^{h, F}(\beta)$ is an increasing  (actually non-decreasing) function. Thus, $\inf\left\{\beta: \beta \in \cA^{h, \phi}_\beta \cup \cR^{h, \phi}_\beta\right\}$ increases with $\phi$. Conclusively,  we get that $\phi^*$ is an  optimizer of \eqref{eqn_opt_phi}. \eop

\end{pfcorollaryehWM}




\end{document}